\documentclass{amsart}

\newtheorem{theorem}{Theorem}[section]
\newtheorem{lemma}[theorem]{Lemma}
\newtheorem{proposition}[theorem]{Proposition}
\newtheorem{corollary}[theorem]{Corollary}

\newtheorem{_conjecture}[theorem]{Conjecture}

\newtheorem{_problem}[theorem]{Problem}

\newtheorem{_claim}[theorem]{Claim}

\newtheorem{_subclaim}[theorem]{Sub-claim}

\newtheorem{defini}[theorem]{Definition}

\newtheorem{rem}[theorem]{\it Remark}
\newenvironment{remark}{\begin{rem}\rm}{\end{rem}}

\newtheorem{_example}[theorem]{Example}

\numberwithin{equation}{section}
\numberwithin{table}{section}
\numberwithin{figure}{section}

\newcommand{\F}{\mathord{\mathbb F}}
\renewcommand{\P}{\mathord{\mathbb  P}}
\newcommand{\Q}{\mathord{\mathbb  Q}}

\newcommand{\Z}{\mathord{\mathbb Z}}

\newcommand{\EEE}{\mathord{\mathcal E}}

\newcommand{\OOO}{\mathord{\mathcal O}}

\newcommand{\RRR}{\mathord{\mathcal R}}
\newcommand{\SSS}{\mathord{\mathcal S}}

\newcommand{\maprightsp}[1]{\; \smash{\mathop{\; \longrightarrow \; }\limits\sp{#1}}\; }

\newcommand{\mapdownsurj}{
\hbox{$\bigm\downarrow$}
\llap{\hbox{\raise 2pt\hbox{$\bigm\downarrow$}}}%
\vstrechmapdown
}

\newcommand{\inj}{\hookrightarrow}

\newcommand{\set}[2]{\{\; {#1} \; \mid \; {#2} \;  \}}

\newcommand{\map}[3]{ #1 \, : \, #2 \, \to \, #3}

\newcommand{\shortmap}[3]{ #1 : #2 \to #3}

\newcommand{\shortmapinj}[3]{ #1 : #2 \inj  #3}

\newcommand{\sm}{\setminus}
\newcommand{\st}{\subset}
\newcommand{\steq}{\subseteq}

\newcommand{\sprime}{\sp\prime}

\newcommand{\spprime}{\sp{\prime\prime}}

\newcommand{\sperp}{\sp{\perp}}
\newcommand{\dual}{\sp{\vee}}

\newcommand{\inv}{\sp{-1}}

\renewcommand{\qed}{\hfill {$\Box$}}

\newcommand{\pr}{\operatorname{\rm pr}\nolimits}

\newcommand{\Aut}{\operatorname{\rm Aut}\nolimits}
\newcommand{\Hom}{\operatorname{\rm Hom}\nolimits}

\newcommand{\Sing}{\operatorname{\rm Sing}\nolimits}

\newcommand{\rank}{\operatorname{\rm rank}\nolimits}
\newcommand{\disc}{\operatorname{\rm disc}\nolimits}

\newcommand{\closure}[1]{\overline{#1}}

\newcommand{\rmand}{\textrm{and}}

\newcommand{\quand}{\quad\rmand\quad}

\newcommand{\Thm}[1]{Theorem~\ref{thm:#1}}
\newcommand{\Prop}[1]{Proposition~\ref{prop:#1}}
\newcommand{\Lem}[1]{Lemma~\ref{lem:#1}}
\newcommand{\Cor}[1]{Corollary~\ref{cor:#1}}

\newcommand{\roots}{{\mathord{\rm roots}}}
\newcommand{\rootsnum}{{\mathord{\varrho}}}
\newcommand{\Roots}{{\mathord{\rm Roots}}}

\newcommand{\NP}{\mathord{\it NP}}
\newcommand{\MW}{\mathord{\it MW\hskip -1pt}}
\newcommand{\ord}{{\mathord{\rm ord}}}
\newcommand{\Id}{{\mathord{\rm Id}}}
\newcommand{\ECD}[1]{\langle #1\rangle}

\begin{document}

\title[Supersingular $K3$ surfaces]{Rational double points 
on supersingular  $K3$ surfaces}

\author{Ichiro Shimada}
\address{
Department of Mathematics,
Faculty of Science,
Hokkaido University,
Sapporo 060-0810,
JAPAN
}
\email{shimada@@math.sci.hokudai.ac.jp
}




\begin{abstract}
We  investigate configurations of rational double points with the total Milnor number $21$
on supersingular $K3$ surfaces. 
The complete list of possible configurations is given.
As an application,
we also give the complete list of extremal (quasi-)elliptic   fibrations
on supersingular $K3$ surfaces.
\end{abstract}

\maketitle
\section{Introduction}\label{sec:intro}
Let $Y$  be a normal projective surface 
defined over an algebraically closed field $k$,
and let $f : X\to Y$ be the minimal resolution.
Suppose that $X$ is a $K3$ surface.
Then the normal surface $Y$ has only rational double points  as its singularities.
(See \cite{Artin62, Artin66, Artin77} for the definition of
rational double points.)
The {\it total Milnor number} $\mu (Y)$
of $Y$ is, by definition, the number of $(-2)$-curves
that are contracted by $f$.
It is obvious that $\mu (Y)$
is less than the Picard number  of $X$.
In particular,
$\mu (Y)$ cannot exceed $19$
in characteristic $0$.
In positive characteristics,
however,
there exist {\it supersingular} $K3$ surfaces (in the sense of Shioda~\cite{Shioda77}),
and  
we have  singular $K3$ surfaces $Y$
with $\mu (Y)\ge 20$.
\par
For example,
Dolgachev and Kondo
\cite{DK}
showed that
a supersingular $K3$ surface $X$ in characteristic $2$
with the Artin  invariant $1$
is birational to 
the quartic surface  in $\P\sp 3$  defined by the equation
\begin{equation*}\label{eq:DK}
x_0^4+x_1^4+x_2^4+x_3^4+x_0^2x_1^2+x_0^2x_2^2+x_1^2x_2^2+x_0x_1x_2(x_0+x_1
+x_2) = 0.
\end{equation*}
This quartic surface  
 has seven rational double points of type $A\sb 3$
so that its total Milnor number  is $21$.
They also showed that $X$ is birational to the purely inseparable double cover
of $\P\sp 2$ defined by
$$
w^2 = x_0x_1x_2(x_0^3+x_1^3+x_2^3), 
$$
which  has $21$ rational double points of type $A\sb 1$.
\par
We say that $Y$ is a {\it supersingular $K3$ surface with maximal rational double points}
if $\mu(Y)$
attains the possible maximum $21$.
It is quite plausible that  $K3$ surfaces
with this property
have many interesting  features
that are peculiar to  algebraic varieties in positive characteristics.
The purpose of this paper is to investigate
the combinatorial aspects of such supersingular $K3$ surfaces 
by means of the lattice theory.
\par
\medskip
An {\it $ADE$}-type is a finite formal sum
$$
R=\sum a\sb l A\sb l +\sum d\sb m D\sb m +\sum e\sb n E\sb n
$$
of symbols $A\sb l \;(l\ge 1)$, $D\sb m \;(m\ge 4)$ and $E\sb n \;(n=6,7,8)$
with non-negative integer coefficients,
and its rank is defined by
$$
\rank (R):=\sum a\sb l  l +\sum d\sb m  m +\sum e\sb n n.
$$
For an $ADE$-type $R$,
we denote by $Q(R)$
the {\it negative} definite integer lattice
whose intersection matrix is  the Cartan matrix of type $R$
multiplied by $-1$.
The rank of $Q(R)$ is therefore equal to $\rank (R)$.
\par
Let $Y$ be a supersingular $K3$ surface with maximal rational double points,
and $f : X\to Y$ its minimal resolution.
We denote by
$R\sb Y$ the $ADE$-type of rational double points on $Y$.
In the Picard  lattice $S\sb X$ of $X$,
we have a negative definite sublattice $T\sb Y$
generated by the classes of $(-2)$-curves
contracted by $f$.
By definition,
we have $\rank (R\sb Y)=\mu (Y)=21$,  and 
$T\sb Y$
is isomorphic to $Q(R\sb Y)$.
The orthogonal complement 
of $T\sb Y$ in $S\sb X$
is therefore
generated by a vector $h\sb Y\in S\sb X$,
and its norm 
$$
n\sb Y:= h\sb Y\sp 2
$$
is uniquely determined.
Note that $n\sb Y$ is a positive even integer.
Since $X$ is a supersingular $K3$ surface,
it is known from \cite{Artin74} that
the discriminant of  $S\sb X$
is equal to $-p\sp{2\sigma\sb X}$,
where $p$
is the characteristic of the base field $k$,
and $\sigma\sb X$ is a positive integer $\le 10$,
which is called the {\it Artin invariant} of $X$.
Thus we obtain  a triple
$$
(R\sb Y, n\sb Y, \sigma\sb X).
$$
\par
For simplicity,
we define  an {\it $RDP$-triple} to be a triple $(R, n, \sigma)$
consisting of an $ADE$-type $R$ with $\rank (R)=21$,
a positive even integer $n$,
and a positive integer $\sigma \le 10$.
We say that an $RDP$-triple 
$(R, n, \sigma)$
is {\it geometrically realizable in characteristic $p$}
if it satisfies the following conditions,
which are equivalent to each other:
\begin{itemize}
\item[(i)]
There exists a supersingular $K3$ surface $Y$ with maximal rational double points in characteristic $p$
such that $(R\sb Y, n\sb Y, \sigma\sb X)=(R, n, \sigma)$.
\item[(ii)]
Every (smooth) supersingular $K3$ surface $X$ with the Artin  invariant $\sigma$
in characteristic $p$  admits a birational morphism  
$f: X\to Y$  
to a supersingular $K3$ surface  $Y$ with maximal rational double points
such that 
$R\sb Y=R$ and $n\sb Y=n$.
\end{itemize}
See \S\ref{sec:preliminaries} for the equivalence of these conditions.
\par
\medskip
Our main results are as follows:
\begin{theorem}
An $RDP$-triple 
$(R, n, \sigma)$
is  geometrically realizable in characteristic $p$
if and only if $(R, n, \sigma)$ is a member of the list  
given in Table {\rm RDP}  at the end of this paper.
\end{theorem}
\begin{corollary}\label{cor:cor1}
There exists a supersingular $K3$ surface with maximal rational double points
in
characteristic $p$ if and only if $p\le 19$.
\qed
\end{corollary}
The appearance  of $(21 A\sb 1, 2, \sigma)$
with $\sigma =1, \dots, 10$ in the list
of $RDP$-triples geometrically realizable in characteristic $2$ implies that 
every supersingular $K3$ surface in characteristic $2$ is birational to a projective surface with $21$ 
ordinary nodes.
In fact, we will prove the following:
\begin{proposition}\label{prop:21A1}
Every supersingular $K3$ surface in characteristic $2$ is birational to 
a purely inseparable double cover of $\P\sp 2$ that has  $21$ ordinary nodes.
\end{proposition}
This proposition
gives us another proof
of the unirationality
of supersingular $K3$ surfaces in characteristic $2$,
which was first 
proved by Rudakov-{\v{S}}afarevi{\v{c}} in 
\cite{RS_char2}.
\par
\medskip
As an application,
we give the  complete list of extremal (quasi-)elliptic   fibrations on 
supersingular $K3$ surfaces.

Let $X$ be a $K3$ surface.
A {\it {\rm(}quasi-{\rm)}elliptic   fibration} on $X$ is, by definition,
a surjective morphism 
$\phi : X\to \P\sp 1$ such that the general fiber $F$ of $\phi$ is 
a reduced irreducible curve of arithmetic genus $1$,
and that there exists a distinguished  section $O : \P\sp 1 \to X$ of $\phi$.
We say that $\phi: X\to \P\sp 1$ is  {\it elliptic}  if $F$ is smooth,
and $\phi$ is {\it quasi-elliptic} if $F$ is singular.
A quasi-elliptic   fibration exists only in characteristic $2$ or $3$ (\cite{Tate, Shimada92}).

Let $\phi : X\to \P\sp 1$ be a (quasi-)elliptic   fibration.
The generic fiber $E$ of $\phi$ is a curve defined over the rational function field $K:=k (\P\sp 1)$
of the base curve, and 
the  set $E(K)$ of $K$-rational points of $E$
is endowed with a structure of the abelian group
such that the point $o\in E(K)$ corresponding to the section $O: \P\sp 1\to X$
is the zero element.
We call $E(K)$ the {\it Mordell-Weil group} of $\phi : X\to \P\sp 1$,
and denote it  by $\MW\sb{\phi}$.
Let $T\sb{\phi}$ be the sublattice of $S\sb X$ spanned by the classes of the general fiber of $\phi$,
the zero section $O$ and 
the irreducible components 
of reducible fibers of $\phi$ that are disjoint from $O$.
Then we have the following famous formula
(\cite{Shioda90}, \cite{Ito92}):
\begin{equation}\label{eq:MW}
\MW\sb{\phi}\;\;\cong\;\; S\sb X/ T\sb{\phi}.
\end{equation}
Let $\RRR\sb{\phi}$ be the set of points $v\in \P\sp 1$ of the base  curve 
such that the fiber $\phi\inv (v)$ is reducible.
It is well-known that, for each $v\in \RRR\sb\phi$, 
the classes of  irreducible components
of $\phi\inv  (v)$ disjoint from the zero section
span an indecomposable negative-definite $ADE$-lattice
in $S\sb X$.
Let $R\sb v$ be the $ADE$-type of this sublattice, 
and  put 
$$
R\sb{\phi} := \sum\sb{v\in \RRR\sb{\phi}} R\sb v.
$$
Since $\rank T\sb\phi=\rank (R\sb{\phi})+2$, 
we have $\rank (R\sb{\phi})\le 20$.
We say that $\phi$ is {\it extremal} if $\rank (R\sb{\phi})$   attains the possible maximum  $20$.
It follows  that, if $\phi$ is an extremal (quasi-)elliptic   fibration, 
then $X$ is supersingular
and $\MW\sb{\phi}$ is  finite.
On the other hand,
if $\phi$ is quasi-elliptic,
then $\phi$ is necessarily extremal (\cite{RS}).

A triple $\ECD{R,  \sigma, \MW}$ consisting of  an $ADE$-type $R$ of rank $20$,
an integer $\sigma$ with $1\le \sigma \le 10$,
and a finite abelian group $\MW$,
is called an {\it elliptic triple}.
An elliptic triple
$\ECD{ R, \sigma, \MW}$ is called a 
 {\it triple of extremal  elliptic
{\rm(}resp.~quasi-elliptic{\rm)}   $K3$ surfaces in characteristic $p$}
if it satisfies the following conditions,
which are equivalent to each other:
\begin{itemize}
\item[(i)]
There exists 
a $K3$ surface $X$
with the Artin invariant $\sigma$ 
in characteristic $p$ 
that has an extremal elliptic (resp.~quasi-elliptic)   fibration   $\phi : X\to \P\sp 1$ 
such that 
$R\sb\phi= R$ and $\MW\sb{\phi}\cong \MW$.
\item[(ii)]
Every  supersingular $K3$ surface $X$ with the Artin  invariant $\sigma$
in characteristic $p$ admits 
 an extremal elliptic (resp.~quasi-elliptic)    fibration  $\phi : X\to \P\sp 1$ 
such that 
$R\sb\phi= R$ and $\MW\sb{\phi}\cong \MW$.
\end{itemize}
\begin{theorem}
The complete lists of  triples of extremal  quasi-elliptic
and  elliptic   $K3$ surfaces
are given by Tables {\rm QE} and {\rm E}.
\end{theorem}

\def\ttskip{0em}
\def\pttskip{0.1pt}
\def\colwidth{185pt}

\def\tablerule{\noalign{\vskip .6pt\hrule\vskip .6pt}}

\begin{table}
\begin{center}
{\bf
 Table QE: 
The complete list of  extremal quasi-elliptic $K3$ surfaces
}
\end{center}
\par
\smallskip
\hbox{
$p=2$, $\MW= (\Z/2\Z)\sp{\oplus r}$.
}
\par
\smallskip
{\scriptsize
%
%
\hbox{
  \vtop{\tabskip=0pt \offinterlineskip
    \halign to \colwidth {\strut\vrule#\tabskip =\ttskip plus\pttskip&#\hfil&\vrule#&\hfil#& \vrule#& #\hfil & \vrule#\tabskip=0pt\cr
    \tablerule
    &\hfil$R$&& \hfil $\sigma$ \hfil   && \hfil $r$ &\cr
    \tablerule
    \tablerule
    &\hfil$  2 E\sb{8} + D\sb{4} $&& $1$  \hfil && \hfil $0$ &\cr
    \tablerule
    &\hfil$ E\sb{8} + E\sb{7} +  5 A\sb{1} $&& $2$  \hfil && \hfil $1$ &\cr
    \tablerule
    &\hfil$ E\sb{8} + D\sb{12} $&& $1$  \hfil && \hfil $0$ &\cr
    \tablerule
    &\hfil$ E\sb{8} + D\sb{8} + D\sb{4} $&& $2$  \hfil && \hfil $0$ &\cr
    \tablerule
    &\hfil$ E\sb{8} + D\sb{6} +  6 A\sb{1} $&& $3$  \hfil && \hfil $1$ &\cr
    \tablerule
    &\hfil$ E\sb{8} +  3 D\sb{4} $&& $3$  \hfil && \hfil $0$ &\cr
    \tablerule
    &\hfil$ E\sb{8} + D\sb{4} +  8 A\sb{1} $&& $4$  \hfil && \hfil $1$ &\cr
    \tablerule
    &\hfil$ E\sb{8} +  12 A\sb{1} $&& $5$  \hfil && \hfil $1$ &\cr
    \tablerule
    &\hfil$  2 E\sb{7} + D\sb{6} $&& $1$  \hfil && \hfil $1$ &\cr
    \tablerule
    &\hfil$  2 E\sb{7} + D\sb{4} +  2 A\sb{1} $&& $2$  \hfil && \hfil $1$ &\cr
    \tablerule
    &\hfil$  2 E\sb{7} +  6 A\sb{1} $&& $3$  \hfil && \hfil $1$ &\cr
    \tablerule
    &\hfil$ E\sb{7} + D\sb{10} +  3 A\sb{1} $&& $1, 2$ \hfil  && \hfil $3-\sigma$ &\cr
    \tablerule
    &\hfil$ E\sb{7} + D\sb{8} +  5 A\sb{1} $&& $2, 3$ \hfil  && \hfil $4-\sigma$ &\cr
    \tablerule
    &\hfil$ E\sb{7} +  2 D\sb{6} + A\sb{1} $&& $2$  \hfil && \hfil $1$ &\cr
    \tablerule
    &\hfil$ E\sb{7} + D\sb{6} + D\sb{4} +  3 A\sb{1} $&& $2, 3$ \hfil  && \hfil $4-\sigma$ &\cr
    \tablerule
    &\hfil$ E\sb{7} + D\sb{6} +  7 A\sb{1} $&& $3, 4$ \hfil  && \hfil $5-\sigma$ &\cr
    \tablerule
    &\hfil$ E\sb{7} +  2 D\sb{4} +  5 A\sb{1} $&& $3, 4$ \hfil  && \hfil $5-\sigma$ &\cr
    \tablerule
    &\hfil$ E\sb{7} + D\sb{4} +  9 A\sb{1} $&& $3, 4, 5$ \hfil  && \hfil $6-\sigma$ &\cr
    \tablerule
    &\hfil$ E\sb{7} +  13 A\sb{1} $&& $4, 5, 6$ \hfil  && \hfil $7-\sigma$ &\cr
    \tablerule
    &\hfil$ D\sb{20} $&& $1$  \hfil && \hfil $0$ &\cr
    \tablerule
    &\hfil$ D\sb{16} + D\sb{4} $&& $1, 2$ \hfil  && \hfil $2-\sigma$ &\cr
    \tablerule
    &\hfil$ D\sb{14} +  6 A\sb{1} $&& $2, 3$ \hfil  && \hfil $4-\sigma$ &\cr
    \tablerule
  }
 }
  \vtop{\tabskip=0pt \offinterlineskip
    \halign to \colwidth {\strut\vrule#\tabskip =\ttskip plus\pttskip&#\hfil&\vrule#&\hfil#& \vrule#& #\hfil & \vrule#\tabskip=0pt\cr
    \tablerule
    &\hfil$R$&& \hfil $\sigma$ \hfil   && \hfil $r$ &\cr
    \tablerule
    \tablerule
    &\hfil$ D\sb{12} + D\sb{8} $&& $1, 2$ \hfil  && \hfil $2-\sigma$ &\cr
    \tablerule
    &\hfil$ D\sb{12} +  2 D\sb{4} $&& $2, 3$ \hfil  && \hfil $3-\sigma$ &\cr
    \tablerule
    &\hfil$ D\sb{12} +  8 A\sb{1} $&& $3, 4$ \hfil  && \hfil $5-\sigma$ &\cr
    \tablerule
    &\hfil$ D\sb{10} + D\sb{6} +  4 A\sb{1} $&& $2, 3$ \hfil  && \hfil $4-\sigma$ &\cr
    \tablerule
    &\hfil$ D\sb{10} + D\sb{4} +  6 A\sb{1} $&& $3, 4$ \hfil  && \hfil $5-\sigma$ &\cr
    \tablerule
    &\hfil$ D\sb{10} +  10 A\sb{1} $&& $4, 5$ \hfil  && \hfil $6-\sigma$ &\cr
    \tablerule
    &\hfil$  2 D\sb{8} + D\sb{4} $&& $1, 2, 3$ \hfil  && \hfil $3-\sigma$ &\cr
    \tablerule
    &\hfil$ D\sb{8} + D\sb{6} +  6 A\sb{1} $&& $2, 3, 4$ \hfil  && \hfil $5-\sigma$ &\cr
    \tablerule
    &\hfil$ D\sb{8} +  3 D\sb{4} $&& $2, 3, 4$ \hfil  && \hfil $4-\sigma$ &\cr
    \tablerule
    &\hfil$ D\sb{8} + D\sb{4} +  8 A\sb{1} $&& $3, 4, 5$ \hfil  && \hfil $6-\sigma$ &\cr
    \tablerule
    &\hfil$ D\sb{8} +  12 A\sb{1} $&& $4, 5, 6$ \hfil  && \hfil $7-\sigma$ &\cr
    \tablerule
    &\hfil$  3 D\sb{6} +  2 A\sb{1} $&& $1, 2, 3$ \hfil  && \hfil $4-\sigma$ &\cr
    \tablerule
    &\hfil$  2 D\sb{6} + D\sb{4} +  4 A\sb{1} $&& $2, 3, 4$ \hfil  && \hfil $5-\sigma$ &\cr
    \tablerule
    &\hfil$  2 D\sb{6} +  8 A\sb{1} $&& $3, 4, 5$ \hfil  && \hfil $6-\sigma$ &\cr
    \tablerule
    &\hfil$ D\sb{6} +  2 D\sb{4} +  6 A\sb{1} $&& $2, 3, 4, 5$ \hfil  && \hfil $6-\sigma$ &\cr
    \tablerule
    &\hfil$ D\sb{6} + D\sb{4} +  10 A\sb{1} $&& $3, 4, 5, 6$ \hfil  && \hfil $7-\sigma$ &\cr
    \tablerule
    &\hfil$ D\sb{6} +  14 A\sb{1} $&& $3, 4, 5, 6, 7$ \hfil  && \hfil $8-\sigma$ &\cr
    \tablerule
    &\hfil$  5 D\sb{4} $&& $1, 2, 3, 4, 5$ \hfil  && \hfil $5-\sigma$ &\cr
    \tablerule
    &\hfil$  3 D\sb{4} +  8 A\sb{1} $&& $2, 3, 4, 5, 6$ \hfil  && \hfil $7-\sigma$ &\cr
    \tablerule
    &\hfil$  2 D\sb{4} +  12 A\sb{1} $&& $3, 4, 5, 6, 7$ \hfil  && \hfil $8-\sigma$ &\cr
    \tablerule
    &\hfil$ D\sb{4} +  16 A\sb{1} $&& $2, 3, 4, 5, 6, 7, 8$ \hfil  && \hfil $9-\sigma$ &\cr
    \tablerule
    &\hfil$  20 A\sb{1} $&& $3, 4, 5, 6, 7, 8, 9$ \hfil  && \hfil $10-\sigma$ &\cr
    \tablerule
  }
 }
}
}
\par
\bigskip
\hbox{
$p=3$, $\MW= (\Z/3\Z)\sp{\oplus r}$.
}
\par
\smallskip
{\scriptsize 
%
%
\hbox{
  \vtop{\tabskip=0pt \offinterlineskip
    \halign to \colwidth {\strut\vrule#\tabskip =\ttskip plus\pttskip&#\hfil&\vrule#&\hfil#& \vrule#& #\hfil & \vrule#\tabskip=0pt\cr
    \tablerule
    &\hfil$R$&& \hfil $\sigma$ \hfil   && \hfil $r$ &\cr
    \tablerule
    \tablerule
    &\hfil$  2 E\sb{8} +  2 A\sb{2} $&&  $1$  \hfil && \hfil $0$ &\cr
    \tablerule
    &\hfil$ E\sb{8} +  2 E\sb{6} $&&  $1$  \hfil && \hfil $0$ &\cr
    \tablerule
    &\hfil$ E\sb{8} + E\sb{6} +  3 A\sb{2} $&&  $2$  \hfil && \hfil $0$ &\cr
    \tablerule
    &\hfil$ E\sb{8} +  6 A\sb{2} $&&  $2, 3$  \hfil && \hfil $3-\sigma$ &\cr
    \tablerule
  }
 }
  \vtop{\tabskip=0pt \offinterlineskip
    \halign to \colwidth {\strut\vrule#\tabskip =\ttskip plus\pttskip&#\hfil&\vrule#&\hfil#& \vrule#& #\hfil & \vrule#\tabskip=0pt\cr
    \tablerule
    &\hfil$R$&& \hfil $\sigma$ \hfil   && \hfil $r$ &\cr
    \tablerule
    \tablerule
    &\hfil$  3 E\sb{6} + A\sb{2} $&&  $1, 2$  \hfil && \hfil $2-\sigma$ &\cr
    \tablerule
    &\hfil$  2 E\sb{6} +  4 A\sb{2} $&&  $1, 2, 3$  \hfil && \hfil $3-\sigma$ &\cr
    \tablerule
    &\hfil$ E\sb{6} +  7 A\sb{2} $&&  $2, 3, 4$  \hfil && \hfil $4-\sigma$ &\cr
    \tablerule
    &\hfil$  10 A\sb{2} $&&  $1, 2, 3, 4, 5$  \hfil && \hfil $5-\sigma$ &\cr
    \tablerule
  }
 }
}
}
\par
\bigskip
\begin{center}
{\bf 
Table E:
The complete list of  extremal elliptic  $K3$ surfaces
}
\end{center}
\par
\smallskip
\hbox{
$[a]= \Z/a\Z, \quad [a, b]= \Z/a\Z\oplus \Z/b\Z$.
}
{\scriptsize
%
%
\hbox{
  \vtop{\tabskip=0pt \offinterlineskip
    \halign to \colwidth {\strut\vrule#\tabskip =\ttskip plus\pttskip&#&\vrule#&#& \vrule#&#& \vrule#& #& \vrule#\tabskip=0pt\cr
    \tablerule
    &\hfil $p$ \hfil&&\hfil$R$\hfil&& \hfil$\sigma$ \hfil && \hfil $\MW$\hfil &\cr
    \tablerule
    \tablerule
 &\hfil  $11$   &&\hfil$  2 A\sb{10} $\hfil && \hfil $1$ \hfil && \hfil $0$\hfil &\cr
    \tablerule
 &\hfil  $7$   &&\hfil$ A\sb{13} + A\sb{6} + A\sb{1} $\hfil && \hfil $1$ \hfil &&  \hfil$[2]$\hfil &\cr
    \tablerule
 &\hfil  $7$   &&\hfil$ E\sb{8} +  2 A\sb{6} $\hfil && \hfil $1$ \hfil && \hfil $0$\hfil &\cr
    \tablerule
 &\hfil  $5$   &&\hfil$ E\sb{7} + A\sb{9} + A\sb{4} $\hfil && \hfil $1$ \hfil &&  \hfil$[2]$\hfil &\cr
    \tablerule
 &\hfil  $5$   &&\hfil$ A\sb{14} + A\sb{4} + A\sb{2} $\hfil && \hfil $1$ \hfil &&  \hfil$[3]$\hfil &\cr
    \tablerule
 &\hfil  $3$   &&\hfil$ D\sb{16} +  2 A\sb{2} $\hfil && \hfil $1$ \hfil &&  \hfil$[2]$\hfil &\cr
    \tablerule
 &\hfil  $3$   &&\hfil$ D\sb{10} +  2 A\sb{5} $\hfil && \hfil $1$ \hfil &&  \hfil $[2, 2]$\hfil &\cr
    \tablerule
  }
 }
  \vtop{\tabskip=0pt \offinterlineskip
    \halign to \colwidth {\strut\vrule#\tabskip =\ttskip plus\pttskip&#&\vrule#&#& \vrule#&#& \vrule#& #& \vrule#\tabskip=0pt\cr
    \tablerule
    &\hfil $p$ \hfil&&\hfil$R$\hfil&& \hfil$\sigma$ \hfil && \hfil $\MW$\hfil &\cr
    \tablerule
    \tablerule
 &\hfil  $3$   &&\hfil$ D\sb{7} + A\sb{11} + A\sb{2} $\hfil && \hfil $1$ \hfil &&  \hfil$[4]$\hfil &\cr
    \tablerule
 &\hfil  $2$   &&\hfil$ A\sb{17} +  3 A\sb{1} $\hfil && \hfil $1$ \hfil &&  \hfil$[6]$\hfil &\cr
    \tablerule
 &\hfil  $2$   &&\hfil$  4 A\sb{5} $\hfil && \hfil $1$ \hfil &&  \hfil $[3, 6]$\hfil &\cr
    \tablerule
 &\hfil  $2$   &&\hfil$  2 A\sb{9} +  2 A\sb{1} $\hfil && \hfil $1$ \hfil &&  \hfil$[10]$\hfil &\cr
    \tablerule
 &\hfil  $2$   &&\hfil$ E\sb{6} + A\sb{11} + A\sb{3} $\hfil && \hfil $1$ \hfil &&  \hfil$[6]$\hfil &\cr
    \tablerule
 &\hfil  $2$   &&\hfil$ D\sb{5} + A\sb{15} $\hfil && \hfil $1$ \hfil &&  \hfil$[4]$\hfil &\cr
    \tablerule
  }
 }
}
}
\end{table}
The complete lists of 
triples of extremal  elliptic
   $K3$ surfaces 
and  of extremal 
 quasi-elliptic   $K3$ surfaces in
characteristic $3$
have been already obtained by Ito~\cite{Ito92, Ito97, Ito02}.
His method is completely different from ours, and 
he gave  explicit defining equations
of these extremal (quasi-)elliptic $K3$ surfaces.
The complete list of 
triples of extremal  quasi-elliptic   $K3$ surfaces in
characteristic $2$
seems to be new.
\par
\medskip
In \cite{Y1} and \cite{Y2},
Yang classified all configurations of rational double points
on complex sextic double planes and complex quartic surfaces.
He used the ideas of Urabe \cite{U},
and reduced the problem of listing up all rational double points on these complex 
$K3$ surfaces to lattice theoretic
calculations via Torelli theorem.
By a similar method,
the complete list of configurations of singular fibers 
on complex elliptic $K3$ surfaces has been obtained in \cite{ShZh} and \cite{Sh}.
In \cite{DMP},
the maximal configurations of ordinary nodes on rational surfaces in
characteristic $\ne 2$  are investigated.
\par
\medskip
The plan of this paper is as follows.
In \S\ref{sec:preliminaries},
we review some known facts
in the theory of  $K3$ surfaces and the lattice theory.
In \S\ref{subsec:rational double points},
we give  lattice theoretic conditions for a $K3$ surface to 
have a given configuration of $(-2)$-curves
and to have a (quasi-)elliptic   fibration with a given $ADE$-type of reducible fibers.
In \S\ref{subsec:disc},
we briefly review the theory of discriminant forms due to Nikulin \cite{N}.
In \S\ref{subsec:Picard},
we quote from Artin~\cite{Artin74}, 
Rudakov-{\v{S}}afarevi{\v{c}}~\cite{RS_char2, RS} and Shioda~\cite{Shioda}
some fundamental facts about the Picard lattices of supersingular $K3$ surfaces.
These facts  play, in positive characteristics,
the same role as the one
Torelli theorem played for complex $K3$ surfaces
in \cite{Sh, ShZh}, \cite{U} and  \cite{Y1, Y2}.
The algorithms
for obtaining the lists
of geometrically realizable $RDP$-triples 
and of triples of extremal (quasi-)elliptic $K3$ surfaces
are presented in \S\ref{sec:algorithm}
and
\S\ref{sec:elliptic_algorithm},
respectively.
In \S\ref{sec:21A1},
we investigate the geometry of supersingular $K3$ surfaces in characteristic $2$
with $21$ ordinary nodes,
and prove
Proposition~\ref{prop:21A1}.
\par
\medskip
The lists in Tables RDP, QE and E 
are also available from the author's homepage:
\begin{center}
\texttt{http://www.math.sci.hokudai.ac.jp/\~{}shimada/ssK3.html}
\end{center}
\par
\medskip
Part of this work was done 
during the author's stay at Korea Institute for Advanced Study
in October 2001.
He would like to thank  Professor Jonghae~Keum for  his warm hospitality.
He also would like to thank
 Professors I.~R.~Dolgachev,  Shigeru Mukai  and  Tetsuji~Shioda
for helpful discussions  and comments.
\section{Preliminaries}\label{sec:preliminaries}
\subsection{Rational double points 
and (quasi-)elliptic   fibrations on a $K3$ surface}\label{subsec:rational double points}
An integer lattice $\Lambda$  is said to be {\it even} if
$v\sp 2\in 2\Z$ for any $v\in \Lambda$.
A vector $v$ of an integer lattice $\Lambda$ is said to be {\it primitive}
if the intersection  of $\Q\cdot v$ and $\Lambda$ 
in $\Lambda\otimes\sb{\Z}\Q$ is generated by $v$.
\par
Let $T$ be a {\it negative} definite even integer lattice.
A vector $v\in T$ is called  a {\it root} if $v\sp 2 =-2$.
Let $\Roots (T)$ be the set of roots of $T$,
and let $T\sb\roots$ be the sublattice of $T$ generated by $\Roots (T)$.
We denote by $\Sigma (T)$ the $ADE$-type of the root lattice $T\sb\roots$;
that is,
$\Sigma  (T)$ is the $ADE$-type such that
$T\sb\roots$ is isomorphic to $Q(\Sigma (T))$.
(See~\cite{B},~\cite{E}.)
\par
Let $X$ be a smooth $K3$ surface defined over an algebraically closed field of  arbitrary characteristic.
The Picard lattice $S\sb X$ of $X$ is an even integer lattice
of signature $(1, \rho\sb{X}-1)$,
where $\rho\sb X$ is the Picard number of $X$.
\par
The purpose of this sub-section is to establish the following propositions:
\begin{proposition}\label{prop:propR}
Let $R$ be an $ADE$-type,
and  $X$  a smooth $K3$ surface.
There exists a birational morphism $f : X\to Y$
with the singularity of $Y$  being rational double points of $ADE$-type $R$
if and only if 
there exists a vector $h\in S\sb X$
such that $h\sp 2 >0$ and $\Sigma (h\sp\perp)=R$,
where $h\sperp$ is the orthogonal complement of $h$ in $S\sb X$.
\end{proposition}
\begin{proposition}\label{prop:ell}
Let $R$ be an $ADE$-type,
$\MW$ an abelian group {\rm(}not necessarily finite{\rm)}, 
and
 $X$   a smooth $K3$ surface.
There exists a {\rm (}quasi-{\rm)}elliptic   fibration  
$\phi : X\to \P\sp 1$
such that $R\sb{\phi}=R$ and $\MW\sb{\phi}\cong\MW$  
if and only if 
there exists an indefinite unimodular sublattice $U\subset S\sb X$
of rank $2$
such that $\Sigma (U\sperp)=R$ and $MW\cong U\sperp/ (U\sperp)\sb{\roots}$,
where $U\sperp$ is the orthogonal complement of $U$ in $S\sb X$.
\end{proposition}
For the proof of these propositions,
we need several lemmas.
A line bundle $L$ on $X$ is said to be {\it numerically effective} ({\it nef}\;)
if $L.C\ge 0$ holds
for any curve $C$ on $X$.
For a line bundle $L$ and a divisor $D$ on $X$,
we denote by $[L]\in S\sb X$ and $[D]\in S\sb X$ 
the corresponding vectors in  $S\sb X$.
\begin{lemma}\label{lem:lem1}
Let $L$ be a nef line bundle on $X$
with $L\sp 2 >0$.
\par
{\rm (1)}
The complete linear system $|2L|$
is fixed-component free and base-point free.
\par
{\rm (2)}
Let $\RRR\sb L$ be the set of reduced irreducible curves $C\sb i$ on $X$
such that $L.C\sb i=0$.
Then $\RRR\sb L$ is an $ADE$-configuration of $(-2)$-curves on $X$
whose $ADE$-type is equal to $\Sigma ([L]\sp\perp)$.
\end{lemma}
\begin{proof}
From Nikulin's proposition~\cite[Proposition 0.1]{N2}, it follows that,
if the complete linear system $|M|$ of a nef line bundle $M$ on $X$ with $M^2>0$  
has a fixed component, then 
there exists a divisor $E$ on $X$ such that $E\sp 2 =0$ and $E.M=1$.
Therefore $|2L|$ is fixed component free.
Then, by Saint-Donat's result (\cite[Corollary 3.2]{SD}),
$|2L|$ is base-point free.
\par
Since $[L]\sperp$
is negative definite
and $C\sp 2 $ is $\ge -2$ for any  reduced irreducible curve $C$ on $X$,
every $C\sb i\in \RRR\sb L$
satisfies $C\sb i\sp 2=-2$
and hence is a $(-2)$-curve.
We put
$$
R\sb L := \set{[C\sb i]}{C\sb i\in \RRR\sb L}.
$$
We will show  that $[L]\sp\perp\sb\roots$
is generated by $R\sb L$.
Let  $v$ be an arbitrary vector in $ \Roots ([L]\sp\perp)$.
Since $v$ or $-v$ is effective,
we can and will assume that $v$ is the class
of an effective divisor $D=\sum\gamma\sb i C\sb i\sprime$,
where $C\sb i\sprime${\,}'s are reduced irreducible curves
and $\gamma\sb i${\,}'s are positive integers.
Since $L$ is nef and $L.D=0$,
it follows that $L.C\sb i\sprime=0$
for each $C\sb i\sprime$, 
and hence $[C\sb i\sprime]$ is contained in $R\sb L$.
\par
The set $R\sb L$ is contained in the set of positive roots
$$
R\sp+:=\set{ r\in \Roots([L]\sperp)}{ r.[H]>0}
$$
in $[L]\sperp$
with respect to the class  $[H]\in S\sb X$
of a hyperplane section $H$  of $X$.
Since every $C\sb i\in \RRR\sb L$ is irreducible,
each $[C\sb i]\in R\sb L$ is indecomposable in $R\sp+$.
Therefore
$R\sb L$
is a fundamental root system of   $[L]\sp\perp\sb{\roots}$.
(See~\cite{B},~\cite{E}.)
The $ADE$-type $\Sigma ([L]\sperp)$ of the root lattice  $[L]\sp\perp\sb\roots$
and that 
of the configuration $\RRR\sb L$
of $(-2)$-curves 
are  thus identical.
\end{proof}
\begin{lemma}\label{lem:lem2}
Let $R$ be an $ADE$-type.
A smooth $K3$ surface $X$ has 
a contraction $f: X\to Y$
of an $ADE$-configuration of  $(-2)$-curves
of type $R$
if and only if
there exists a  nef line bundle $L$
on $X$ such that
$L\sp 2>0$ and 
$\Sigma ([L]\sp\perp)=R$.
\end{lemma}
\begin{proof}
Suppose that $X$
admits a contraction $f: X\to Y$ of an $ADE$-configuration of  $(-2)$-curves
of type $R$.
Let $L$
be the line bundle $\OOO\sb X (f\sp * H)$,
where $H$ is a hyperplane section of $Y$.
Then $L$ is nef and $L\sp 2 >0$.
Since the set 
of reduced irreducible curves on $X$
contracted by $f$
is equal to $\RRR\sb L$, 
$\Sigma ([L]\sp\perp)=R$ follows
from \Lem{lem1}\;(2).
\par
Conversely,
suppose that  there exists  a nef line bundle
$L$ on $X$ with 
$L\sp 2>0$ and 
$\Sigma ([L]\sp\perp)=R$.
By \Lem{lem1}\;(1),
we have a morphism
$$
\map{\Phi\sb{|2L|}}{X}{Z\sprime:=\Phi\sb{|2L|}(X)\;\st\;\P\sp N}
$$
defined by   $|2L|$.
Since
$L\sp 2>0$,
we have $\dim Z\sprime=2$.
Let $Z$ be the normalization of $Z\sprime$,
and let
$$
X\;\maprightsp{f}\; Y \;\maprightsp{g}\; Z
$$
be the Stein factorization of the morphism $X\to Z$
induced from $\Phi\sb{|2L|}$.
Since $g$ is finite,
the set $\RRR\sb L$
coincides with 
the set 
of reduced irreducible curves on $X$
contracted by $f$.
By \Lem{lem1}\;(2),
$f$ is
a contraction 
of an $ADE$-configuration of  $(-2)$-curves
whose type is equal to $\Sigma ([L]\sp\perp)=R$.
\end{proof}
\begin{lemma}\label{lem:lem3}
{\rm (1)}
For any  vector $v\in S\sb X$
with $v\sp 2 >0$,
there exists an isometry
$\gamma$
of $S\sb X$ such that $\gamma (v)$
is the class $[L]$
of a nef line bundle $L$.
\par
{\rm (2)}
For any  primitive vector  $v\in S\sb X$
with $v\sp 2 =0$,
there exists an isometry
$\gamma$ of $S\sb X$ such that $\gamma (v)$
is the class $[F]$
of a reduced irreducible curve $F$ of arithmetic genus $1$.
\end{lemma}
\begin{proof}
Using the scalar multiplication by $-1$ if necessary,
we can  assume that $v$ is the class of an effective divisor.
Then \Lem{lem3} follows from
\cite [Proposition 3 in Section 3]{RS}.
\end{proof}
Let $\phi : X\to\P\sp 1$
be a  {\rm(}quasi-{\rm)}elliptic fibration  on a $K3$ surface $X$.
We denote by  $U\sb\phi$  the sublattice of $S\sb X$
generated by the classes of the general fiber $F$ of $\phi$
and  the zero section.
Then  $U\sb\phi$ is an indefinite unimodular lattice of rank $2$.
Since $U\sb\phi$ is unimodular,
we have an orthogonal decomposition
\begin{equation}\label{eq:orthodecomp}
S\sb X= U\sb\phi \oplus U\sb\phi\sperp.
\end{equation}
It is easy to prove 
$T\sb\phi=U\sb\phi\oplus (U\sb\phi\sperp)\sb{\roots}$
by the same argument as in the proof of Lemma~\ref{lem:lem1}.
(See  also \cite[Lemma 6.1]{Nishiyama}.)
From the definition and~\eqref{eq:MW}, 
we have
\begin{equation}\label{eq:nishiyama}
R\sb\phi=\Sigma (U\sb\phi\sperp)
\quand
\MW\sb\phi\cong U\sb\phi\sperp/ (U\sb\phi\sperp)\sb{\roots}.
\end{equation}
Consider the module
$ [F] \sperp / \Z [F] $,
where $F$ is the general fiber of $\phi$,
$\Z [F]$ is the sublattice of $S\sb X$ generated by $[F]$,
and
$[F] \sperp$ is the orthogonal complement of $ [F]$ in $S\sb X$.
Since $F^2=0$,
we can naturally regard $ [F] \sperp / \Z [F] $
as an integer lattice.
Since the orthogonal complement of $[F]\in U\sb{\phi}$ in $U\sb\phi$ is generated by $[F]$,
we obtain 
from~\eqref{eq:orthodecomp}
the isomorphism of lattices
\begin{equation}\label{eq:kondo}
[F] \sperp / \Z [F] \cong U\sb\phi\sperp.
\end{equation}
\begin{lemma}[\cite{Kondo89}~Lemma 2.1]\label{lem:KN}
Suppose that $S\sb X$ has an indefinite  unimodular sublattice $U\subset S\sb X$
of rank $2$.
Then $X$ has a {\rm(}quasi-{\rm)}elliptic fibration $\phi : X\to \P\sp 1$
such that $  [F] \sperp / \Z [F]  \cong U\sperp$,
where  $U\sperp$ is the orthogonal complement of $U$ in $S\sb X$.
\qed
\end{lemma}
Note that the proof in~\cite{Kondo89} 
is valid  also in positive characteristics
with the only exception
that  $\phi: X\to \P\sp 1$
may be quasi-elliptic.
\begin{proof}[Proof of Propositions~\ref{prop:propR} and~\ref{prop:ell}]
\Prop{propR} follows 
from Lemmas~\ref{lem:lem2} and \ref{lem:lem3}.
\Prop{ell} follows 
from the isomorphisms~\eqref{eq:nishiyama}, \eqref{eq:kondo} and Lemma~\ref{lem:KN}.
\end{proof}
\subsection{The discriminant form of an even integer lattice}\label{subsec:disc}
Let $\Lambda$ be a non-degenerate even  integer lattice.
We denote by $\Lambda\dual$
the dual lattice of $\Lambda$,
which is the free $\Z$-module
$\Hom (\Lambda, \Z)$
equipped with the natural
symmetric bilinear form
$$
\map{\psi\sb{\Lambda}}{\Lambda\dual\times\Lambda\dual}{\Q}.
$$
There exists  a natural  inclusion $\Lambda \inj \Lambda\dual$.
A submodule $\Lambda\sprime$ of $\Lambda\dual$ containing $\Lambda$
is said to be an {\it overlattice} of $\Lambda$
if $\psi\sb{\Lambda}$ restricted to $\Lambda\sprime\times\Lambda\sprime$
takes values in  $\Z$.
\par
The {\it discriminant group} $G\sb{\Lambda}$
of $\Lambda$ is defined to be $\Lambda\dual/\Lambda$.
Note that the order of $G\sb{\Lambda}$ is equal to $|\disc \Lambda|$.
There is a quadratic form
$$
\map{q\sb{\Lambda}}{G\sb{\Lambda}}{\Q/2\Z}
$$
on $G\sb{\Lambda}$ defined by
$q\sb{\Lambda} (\bar v) := \psi\sb{\Lambda} (v, v) \bmod{2\Z}$,
where $\bar v := v + \Lambda\in G\sb{\Lambda}$.
The pair $(G\sb{\Lambda}, q\sb{\Lambda})$
is called the {\it discriminant form} of $\Lambda$.
\par
More generally,
let $(G, q)$ be a pair of a finite abelian group
$G$ and a quadratic form 
$\shortmap{q}{G}{\Q/2\Z}$.
We define a symmetric bilinear form 
$$
\map{b}{G\times G}{\Q/\Z}
$$
by $b (x, y) :=(q (x+y) -q (x) -q (y))/2$.
For a subgroup $S$ of $G$,
we define its orthogonal complement
$S\sp{\perp}$ by
$$
S\sp{\perp}:=\set{x\in G}{ b (x, y) =0 \;\textrm{for all}\;  y \in S}.
$$
\begin{proposition}[\cite{N}~Proposition~1.4.1]\label{prop:nikulin}
Let $\shortmap{\pr\sb{\Lambda}}{\Lambda\dual}{G\sb{\Lambda}}$
be the natural projection.
Then the correspondence
$S\mapsto \pr\sb{\Lambda}\inv (S)$
yields a bijection from the set of  subgroups of $(G\sb{\Lambda}, q\sb{\Lambda})$
to the set of even overlattices of $\Lambda$.
If $S$ is an isotropic subgroup of $(G\sb{\Lambda}, q\sb{\Lambda})$,
then the discriminant group  of the even integer lattice $\pr\sb{\Lambda}\inv (S)$
is isomorphic to $S\sp{\perp}/S$. 
\qed
\end{proposition}
\subsection{The Picard lattice of a supersingular $K3$ surface}\label{subsec:Picard}
Let $p$ be a prime integer.
A non-degenerate even integer lattice $\Lambda$
is  called {\it $p$-elementary}
if its discriminant group $G\sb {\Lambda}$ is $p$-elementary;
that is, $p\cdot \Lambda\dual \subseteq \Lambda$ holds.
A $2$-elementary lattice $\Lambda$ is said to be {\it of type {\rm I}}
if $\psi\sb{\Lambda} (v, v) \in \Z$ holds 
for any $v\in \Lambda\dual$,
where $\psi\sb{\Lambda}$
is the natural symmetric bilinear form on $\Lambda\dual$.
\par
Let $M\sb{\Lambda}$ be the intersection matrix of $\Lambda$
with respect to a certain basis of $\Lambda$.
Then $\Lambda$ is $p$-elementary
if and only if all the coefficients of $p\cdot M\sb{\Lambda}\inv$ are integers.
Suppose that $\Lambda$ is $2$-elementary.
Then $\Lambda$ is of type {\rm I}
if and only if all the diagonal coefficients of $M\sb{\Lambda}\inv$ are integers.
\begin{theorem}[\cite{Artin74}, \cite{RS_char2}]\label{thm:pelm}
Let $X$ be a supersingular $K3$ surface in characteristic $p$.
Then the Picard lattice $S\sb X$ of $X$
is a $p$-elementary lattice of signature $(1, 21)$,
and its discriminant is equal to $-p\sp{2\sigma}$,
where $\sigma$ is an integer satisfying  $1\le \sigma \le 10$.
\par
When $p=2$,
$S\sb X$ is of type {\rm I}.
\qed
\end{theorem}
%
%
%
%
\begin{theorem}[\cite{RS_char2}, \cite{Shioda}]\label{thm:exists}
For each pair $(p, \sigma)$ of a prime integer $p$ and an integer $\sigma$ satisfying  $1\le \sigma \le 10$,
there exists a supersingular $K3$ surface 
with  the  Artin  invariant   $\sigma$ in  characteristic $p$.
\qed
\end{theorem}
For a pair $(p, \sigma)$ of a prime integer $p$ and an integer $\sigma$ with $1\le \sigma \le 10$,
let $\Lambda\sb{p, \sigma}$ be a lattice of rank $22$
satisfying the following conditions:
\par
\medskip
(1)
$\Lambda\sb{p, \sigma}$ is a non-degenerate even integer lattice of signature $(1, 21)$,
\par
(2) $\disc \Lambda\sb {p, \sigma}= -p\sp{2\sigma}$, and 
\par
(3) $\Lambda\sb{p, \sigma}$ is $p$-elementary.
\par
\medskip
\noindent
When $p=2$,
we further impose on $\Lambda\sb{p, \sigma}$ the following condition:
\par
\medskip
(4) $\Lambda\sb{p, \sigma}$ is of type {\rm I}.
\begin{theorem}[\cite{RS}, \cite{CS} Chapter~15 ]\label{thm:unique}
These conditions determine the lattice $\Lambda\sb{p, \sigma}$
uniquely up to isomorphisms.
\qed
\end{theorem}
\begin{corollary}\label{cor:isom}
If $X$ is a supersingular $K3$ surface in characteristic $p$
with the Artin invariant $\sigma$,
then the Picard lattice $S\sb X$ of $X$ is isomorphic to $\Lambda\sb{p, \sigma}$.
\qed
\end{corollary}
By Propositions~\ref{prop:propR},~\ref{prop:ell} and \Cor{isom},
we see that the two conditions
in the definition
of geometric realizability of
 $RDP$-triples
and  elliptic triples
given in Introduction 
are equivalent to each other.
\begin{corollary}\label{cor:comb1}
An $RDP$-triple $(R, n, \sigma)$ is geometrically realizable 
in characteristic $p$ if and only if 
there exists a primitive vector $h\in \Lambda\sb{p, \sigma}$
such that $h\sp 2 =n$ and $\Sigma (h\sperp)=R$,
where $h\sperp$ is the orthogonal complement of $h$ in $\Lambda\sb{p, \sigma}$.
\qed
\end{corollary}
\begin{corollary}\label{cor:ell1}
An elliptic  triple $\ECD{R, \sigma, MW}$ is
a triple of extremal {\rm (}quasi-{\rm)}elliptic
$K3$ surface 
in characteristic $p$ if and only if 
there exists an indefinite unimodular   sublattice $U\subset \Lambda\sb{p, \sigma}$
of rank $2$
such that $\Sigma (U\sperp)=R$ and $MW\cong U\sperp /  (U\sperp)\sb{\roots}$,
where $U\sperp$ is the orthogonal complement of $U$ in $\Lambda\sb{p, \sigma}$.
\qed
\end{corollary}
\begin{proposition}\label{prop:ell2}
An elliptic triple $\ECD{R, \sigma, MW}$ is
a triple of extremal {\rm (}quasi-{\rm)} elliptic
$K3$ surface 
in characteristic $p$ if and only if 
there exist vectors $h$ and $z$ in $\Lambda\sb{p, \sigma}$
satisfying 
the following conditions:
\begin{enumerate}
\renewcommand{\labelenumi}{(\theenumi)\hfil{}}
\renewcommand{\theenumi}{\roman{enumi}}
\item
$h^2=2$, $\Sigma (h\sperp)=R+A\sb 1$, 
\item
$z\in \Roots (h\sperp)$,
\item
$rz=0$ for any $r\in \Roots ( h\sperp)\sm\{\pm z\}$
{\rm(}so that $\Roots (h\sperp)\sm\{\pm z\}$
is a  root system of type $R${\rm)}, 
\item 
$h-z$ is divisible by $2$ in $\Lambda\sb{p, \sigma}$, 
and
\item
the sublattice $U$ of $\Lambda\sb{p, \sigma}$
generated by $f:=(h-z)/2$ and $z$
satisfies  $MW\cong U\sperp/ (U\sperp)\sb{\roots}$.
\end{enumerate}
In particular,
if an elliptic  triple $\ECD{R, \sigma, MW}$ is
a triple of extremal {\rm (}quasi-{\rm)}elliptic
$K3$ surface 
in characteristic $p$, 
then 
the $RDP$-triple $(R+A\sb 1, 2, \sigma)$
is geometrically realizable 
in characteristic $p$.
\end{proposition}
\begin{proof}
Suppose that there exist vectors $h$ and $z$ with the properties (i)-(v).
Since $f^2=0, fz=1$ and $z^2=-2$,
the sublattice $U$ is indefinite and unimodular.
From (iii),
we have $\Sigma (U\sperp)=R$.
Thus the condition in \Cor{ell1} is satisfied.
\par
Conversely,
suppose that
there exists an indefinite  unimodular  sublattice $U\subset \Lambda\sb{p, \sigma}$
of rank $2$
such that $\Sigma (U\sperp)=R$ and $MW\cong U\sperp /  (U\sperp)\sb{\roots}$.
We can find
vectors $f$ and $z$ in $U$
that generate $U$ and 
satisfy
$f^2=0, fz=1$ and $z^2=-2$,
because $U$ is
given by the intersection matrix
$
{\small \begin{pmatrix}
0 & 1 \\
1 & 0
\end{pmatrix}
}
$
with respect to a certain basis of $U$.
The vectors $h:=2f+z$ and $z$ obviously satisfy  $h^2=2$,
and the conditions (ii), (iv)  and (v).
Since $U$ is unimodular,
we have an orthogonal decomposition 
$\Lambda\sb{p, \sigma}=U\oplus U\sperp$.
Therefore we have 
an orthogonal decomposition 
$h\sperp=U\sperp\oplus \langle z\rangle$,
and hence 
the condition (iii) is fulfilled and 
$\Sigma (h\sperp)=R+A\sb 1$ holds.
\par
The last assertion follows from~\Cor{comb1},
because a vector $h$ with $h\sp 2=2$ is necessarily primitive.
\end{proof}
In order to determine whether a (quasi-)elliptic   fibration $\phi : X\to \P\sp 1$
is elliptic or quasi-elliptic,
we use the following:
\begin{theorem}[\cite{RS} Theorem in Section 4 ]\label{thm:qe}
Let $\phi : X\to \P\sp 1$
be a {\rm (}quasi-{\rm)}elliptic   fibration on a $K3$ surface $X$ in characteristic $p$,
where $p=2$ or $3$.
Then $\phi$ is quasi-elliptic if and only if
$\phi$ is extremal and $Q(R\sb{\phi})$ is $p$-elementary.
\qed
\end{theorem}
\section{The algorithm for the list of  $RDP$-triples}\label{sec:algorithm}
Let $(R, n, \sigma)$ be an $RDP$-triple.
We denote by $I(n)$ the lattice of rank $1$
generated by a vector $e\sb n$ with $e\sb n\sp 2=n$,
and by $Q(R, n)$ the lattice $Q(R)\oplus I(n)$
of rank $22$ with signature $(1, 21)$.
By \Thm{unique},
we can rephrase \Cor{comb1}
to the following:
\begin{corollary}\label{cor:comb2}
An $RDP$-triple $(R, n, \sigma)$ is geometrically realizable
in characteristic $p$ if and only if there exists an even overlattice $\Lambda$ of $Q(R, n)$
with the following properties;
\begin{itemize}
\item[(1)]
$\disc \Lambda = - p\sp{2\sigma}$,
\item[(2)]
$\Lambda$ is $p$-elementary,
and,
if $p=2$,  $\Lambda$ is of type {\rm I},
\item[(3)]
the vector $\tilde e\sb n := (0, e\sb n)\in Q(R, n)$
remains primitive in $\Lambda$, and
\item[(4)]
$\Roots (\tilde e\sb n\sperp) =\Roots (Q(R))$,
where $\tilde e\sb n\sperp$ is the orthogonal complement of $\tilde e\sb n$ in $\Lambda$,
and $Q(R)$ is regarded as a sublattice of $\Lambda$.
\qed
\end{itemize}
\end{corollary}
Using \Prop{nikulin},
we can further rephrase \Cor{comb2}
as follows.
For simplicity,
we denote by $(G\sb R, q\sb R)$,
$(G\sb n, q\sb n)$ and $(G\sb{R, n}, q\sb{R, n})$
the discriminant forms of $Q(R)$, $I(n)$ and $Q(R, n)$,
respectively.
There exists a natural isomorphism
$$
(G\sb{R, n}, q\sb{R, n}) \cong (G\sb R, q\sb R) \oplus (G\sb n, q\sb n).
$$
Note that $G\sb n$ is a cyclic group
of order $n$ generated by
$$
\varepsilon\sb n := e\sb n \dual + I(n),
$$
and we have $q\sb n (\varepsilon\sb n) =1/n \bmod 2\Z$.
See \cite[\S6]{Sh} for the structure of  $(G\sb R, q\sb R)$.
Let $S$ be an isotropic subgroup of $(G\sb{R, n}, q\sb{R, n})$.
We denote by $\Lambda\sb{S}$ the even overlattice of $Q(R, n)$
corresponding to $S$ via the bijection given in \Prop{nikulin}.
We regard $Q(R)$ as a sublattice of $\Lambda\sb S$.
Let $(\tilde e\sb n)\sperp\sb S$ denote the orthogonal complement
of $\tilde e\sb n=(0, e\sb n)\in \Lambda\sb S$
in $\Lambda\sb S$,
which is a negative definite even integer lattice.
Since $Q(R)$ is contained in $(\tilde e\sb n)\sperp\sb S$,
we have $\Roots (Q(R))\steq \Roots((\tilde e\sb n)\sperp\sb S)$.
We put
$$
\rootsnum\sb R:=|\Roots (Q(R))|,
\quand 
\rootsnum\sb{R, n} (S):=|\Roots((\tilde e\sb n)\sperp\sb S)|.
$$
See~\cite{ShZh} for $\rootsnum\sb R$ and for the method of calculating 
 $\rootsnum\sb{R, n} (S)$.
\begin{corollary}\label{cor:comb3}
An $RDP$-triple $(R, n, \sigma)$ is geometrically realizable
in characteristic $p$ if and only if there exists an 
isotropic subgroup $S$ of $(G\sb{R, n}, q\sb{R, n})$
with the following properties;
\begin{itemize}
\item[(1)]
$S\sperp/S$ is a $p$-elementary group
of order $p\sp{2\sigma}$,
\item[(2)]
$S\cap G\sb n$ is trivial,
\item[(3)]
$\rootsnum \sb{R, n}(S)=\rootsnum\sb R$, and
\item[(4)]
if $p=2$, then $\Lambda\sb S$ is of type {\rm I}.
\qed
\end{itemize}
\end{corollary}
\subsection{Finiteness of the triples $[R, n, p]$}\label{subsec:finite}
In this sub-section,
we show that there exist  only a finite number of triples $[R, n, p]$
such that $(R, n, \sigma)$ is geometrically realizable in characteristic $p$
for some $\sigma$.
\par
Let $R$ be an $ADE$-type with $\rank (R)=21$,
and $n$ a positive even integer.
We denote by $N\sb R$ the minimal positive integer such that
$$
N\sb R \cdot q\sb R (x) =0 \quad\textrm{in \;\;$\Q/ 2\Z$\;\; for all \;\;$x\in G\sb R$. }
$$
Then $N\sb R$ is the least common multiple of $N\sb X$,
where $X$ runs through the set of indecomposable $ADE$-types
that appear in $R$. 
From the intersection matrix of $q\sb X$ given  in  \cite[\S6, Table 6.1]{Sh},
we obtain  Table~\ref{table:NX}.
\begin{table}
\caption{$N\sb X$ for indecomposable $ADE$-type $X$}
\label{table:NX}
\begin{center}
{\small
\renewcommand{\arraystretch}{1.4}
\begin{tabular}{|c||c| c|c|c|c|c|c|c|}
\hline
 $X$ & \multicolumn{2}{|c|}{$A\sb l$} & \multicolumn{3}{|c|}{$D\sb m$} & \multicolumn{3}{|c|}{$E\sb n$} \\
\hline
 & $l$ : even & $l$ : odd & $m=4k$ & $m=4k+2$ & $m$ : odd & $n=6$ & $n=7$ & $n=8$ \\
\hline
$N\sb X$ & $l+1$ & $2(l+1)$ & $2$ & $4$ & $8$ & $3$ & $4$ & $1$ \\
\hline
$|G\sb X|$ & \multicolumn{2}{|c|}{$l+1$} & \multicolumn{3}{|c|}{$4$} & $3$ & $2$ & $1$ \\
\hline
\end{tabular}
}
\end{center}
\end{table}
We denote by $P\sb R$ the set of prime factors of $|G\sb R|=|\disc Q(R)|$.
\begin{lemma}\label{lem:finite}
Let $p$ be a prime integer.
Suppose that $(G\sb{R, n}, q\sb{R, n})$ contains
an isotropic subgroup $S$
such that $S\sperp/S\cong (\Z/(p))\sp{\oplus 2\sigma}$
for some $\sigma>0$, and
that $S\cap G\sb n$ is trivial.
Then $(n, p)$ is contained in the following finite set $\NP(R)$;
$$
\set{(n, p)}{ p\in P\sb R, \;\; 2n\;|\;(N\sb R\cdot p\sp 2),  \;\; p\sp 2 \;|\; (n\cdot  |G\sb R|)
\;\;\;\textrm{and}\;\;\;
\sqrt{n\cdot  |G\sb R|}\in\Z  }.
$$
\end{lemma}
\begin{proof}
First we fix some notation.
For an abelian group $A$ and a prime integer $l$,
we denote by $A\sb l$ the maximal subgroup of $A$
whose order is a power of $l$.
For a discriminant form $(G, q)$,
we denote by $q\sb l$ the restriction of $q$ to $G\sb l$.
Then we have a natural orthogonal decomposition
$$
(G, q) = \bigoplus\sb l (G\sb l, q\sb l),
$$
where $l$ runs through the set of prime factors of $|G|$.
Let $S$ be an isotropic subgroup of $(G, q)$.
Then the $l$-part $(S\sperp)\sb l$ of $S\sperp$
coincides with the orthogonal complement 
$(S\sb l)\sperp$
of $S\sb l$ in $(G\sb l, q\sb l)$,
and $(S\sperp/S)\sb l$ coincides with $(S\sb l)\sperp/ S\sb l$.
\par
Let $l$ be a prime factor of $n$,
and suppose that $l\notin P\sb R$.
Then $(G\sb R)\sb l$ is trivial.
Since $S\cap G\sb n$ is trivial,
$S\sb l$ is also trivial,
and we have
$$
(S\sperp/S)\sb l= (S\sb l)\sperp=(G\sb{R, n})\sb l=(G\sb n)\sb l,
$$
which is a non-trivial.
Therefore we have $l=p$.
Then we get a contradiction with $(S\sperp/S)\sb p=(\Z/(p))\sp{\oplus 2\sigma}$,
because $(G\sb n)\sb l$ is cyclic.
Thus we have proved  that every prime factor of $n$ is contained in $P\sb R$.
It then follows that $p$ is contained in $ P\sb R$,
because $p$ is a prime factor of $|G\sb{R, n}|=n\cdot |G\sb R|$.
\par
We put
\begin{eqnarray*}
T & :=& \set{t\in G\sb n}{ N\sb R\cdot q\sb n (t) =0 }, \quand\\
k & :=& \min \set{\nu \in \Z\sb{>0}}{N\sb R \cdot \nu\sp 2 =0 \bmod 2n }.
\end{eqnarray*}
Note that  $T$ is a cyclic group of order $n/k$
generated by $k\varepsilon \sb n$.
Let 
$$
\map{\pr\sb n}{G\sb{R, n}}{G\sb n}
$$
be the projection onto the second factor.
Then $\pr\sb n (S)$ is contained in $T$.
Indeed,
since $S$ is isotropic, we have
$q\sb{R, n} (x, y)=0$ for any  $(x, y)\in S$ $(x\in G\sb R, \; y\in G\sb n)$.
By the definition of $N\sb R$,
we have
$$
0= N\sb R\cdot q\sb{R, n} (x, y) = N\sb R\cdot q\sb R (x) + N\sb R \cdot q\sb n (y) = N\sb R \cdot q\sb n (y)
$$
for all $(x, y)\in S$.
Let $\shortmapinj{i\sb n}{G\sb n}{G\sb{R, n}}$
be the inclusion given by $y\mapsto (0, y)$.
Then $i\sb n \inv (S\sperp)$ coincides with the orthogonal complement of $\pr\sb n (S)$ in $G\sb n$,
and hence 
 $i\sb n \inv (S\sperp)$  contains
the orthogonal complement $T\sperp$ of $T$ in $G\sb n$,
which is of order $k$
generated by $(n/k)\cdot \varepsilon\sb n$.
On the other hand,
since $i\sb n\inv (S)$ is trivial by the assumed property
$S\cap G\sb n=\{0\}$
of $S$,
the composite
$$
 i\sb n\inv (S\sperp)
\;
\maprightsp{i\sb n}
\;
S\sperp
\;
\maprightsp{}
\;
S\sperp/S\cong (\Z/(p))\sp{\oplus{2\sigma}}
$$
is injective.
Since  $i\sb n\inv (S\sperp)\st G\sb n$ is cyclic,
it follows that $|i\sb n\inv (S\sperp)|=1$ or $|i\sb n\inv (S\sperp)|=p$.
Therefore $|T\sperp|=k$ is either $1$ or $p$.
If $k=1$,
then $N\sb R =0 \bmod 2n$ holds.
If $k=p$, then 
then $N\sb R \cdot p\sp 2 =0 \bmod 2n$ holds.
\par
Finally,
since $|S\sperp/S|= p\sp{2\sigma}$ is equal to $ n  |G\sb R|/ |S|\sp 2$,
it follows that  $n |G\sb R|$ is a square integer
divisible by $p\sp 2$.
\end{proof}
Lemma~\ref{lem:finite} implies, in particular, that if $(R, n, \sigma)$ is geometrically realizable in
characteristic $p$, then $p$ divides $|G\sb R|$.
The `only if\hskip .9pt' part of Corollary~\ref{cor:cor1}
now follows.
This assertion is   
also derived  from  a theorem of Goto  \cite[Theorem 3.7]{Goto}.
\par\medskip
It is easy to list up all $ADE$-types $R$
with rank  $21$.
For each $R$,
we calculate the set $\NP(R)$, and make the list
$$
\RRR:=\set{[R, n, p]}{\rank(R)=21, \;\; (n, p)\in \NP(R)},
$$
which consists of $20169$ triples.
\subsection{Algorithm I}
We make the list 
$$
\closure{\RRR}:=\set{[R, n]}{\hbox{$[R, n, p]\in \RRR$ \;\;for some $p$}},
$$
which consists of $14487$ pairs.
For each pair $[R, n] \in \closure\RRR$,
we do the following calculations.
\par
Let $\SSS\sb{R, n}$ be the set of isotropic subgroups of $(G\sb{R, n}, q\sb{R, n})$,
and let $\Gamma\sb{R, n}$ be the image of the natural homomorphism 
$$
\Aut (Q(R))\times \Aut(I(n))\;\to\; \Aut(G\sb{R, n}, q\sb{R, n}).
$$
See~\cite{Sh} for the structure of the subgroup $\Gamma\sb{R, n}$ of $\Aut(G\sb{R, n}, q\sb{R, n})$.
The group $\Gamma\sb{R, n}$
acts on $\SSS\sb{R, n}$.
We find a subset $\SSS\sb{R, n}\sprime$ of $\SSS\sb{R, n}$
such that the map
$$
\SSS\sb{R, n}\sprime\;\inj\; \SSS\sb{R, n} \;\to\; \Gamma\sb{R, n} \backslash \SSS\sb{R, n}
$$
is surjective.
For each $S\in \SSS\sprime\sb{R, n}$,
we check the conditions
in \Cor{comb3}.
Note that these conditions 
are invariant under the action of $\Gamma\sb{R, n}$.
If all  the conditions
in \Cor{comb3} are satisfied,
we put $(R, n, \sigma)$
in the list of  $RDP$-triples
geometrically realizable in characteristic $p$.
\subsection{Algorithm II}
 Algorithm I takes impractically long time
when the coefficient  of $A\sb 1$ in $R$ is large.
We improve Algorithm I as follows.
Observe the following trivial facts.
Let $S$ and $T$ be two isotropic subgroups 
of $(G\sb{R, n}, q\sb{R, n})$
such that $T\subset S$.
\begin{itemize}
\item[(a)]
If $T\cap G\sb n$ is non-trivial,
then $S\cap G\sb n$ is non-trivial.
\item[(b)]
If $\rootsnum\sb{R, n} (T) > \rootsnum\sb R$, 
then $\rootsnum\sb{R, n} (S) > \rootsnum\sb R$.
\end{itemize}
We do the following calculations for each $[R, n, p]\in \RRR$.
\par\smallskip\noindent
{\bf Step 1.}
We decompose the set
$P\sb R$ of prime factors of $|G\sb R|$
into the disjoint union
$$
P\sb R=A \;\hbox{$\bigsqcup$}\;B
$$
of certain subsets $A$ and $B$ such that $p\in A$.
Note that every prime factor of $|G\sb{R, n}|$ is contained in $P\sb R$
by the definition of $\NP(R)$.
We then decompose 
$(G, q):=(G\sb{R, n}, q\sb{R, n})$
into the orthogonal direct sum
$$
(G, q)=(G\sb A, q\sb A)\oplus (G\sb B, q\sb B),
$$
where $(G\sb A, q\sb A)$ and $(G\sb B, q\sb B)$
are the orthogonal direct sum
of $(G\sb l, q\sb l)$ $(l\in A)$ and
$(G\sb l, q\sb l)$ $(l\in B)$,
respectively.
\par
\smallskip
Let $\Gamma\sb A\st\Aut(G\sb A, q\sb A)$ and $\Gamma\sb{B}\st\Aut(G\sb{B}, q\sb{B})$
be the intersection of $\Gamma\sb{R, n}\st\Aut(G, q)$
with  the subgroups $\Aut(G\sb{A}, q\sb{A})$ and $\Aut(G\sb{B}, q\sb{B})$ of
$$
\Aut(G, q)\;\cong\; \Aut(G\sb{A}, q\sb{A})\times \Aut(G\sb{B}, q\sb{B}),
$$
respectively.
\par\smallskip\noindent
{\bf Step 2.}
Let 
$\SSS\sb{B}$ be  the set of the isotropic subgroups of $(G\sb{B}, q\sb{B})$,
on which $\Gamma\sb{B}$ acts.
We  find a subset $\SSS\sprime\sb{B}$ of $\SSS\sb{B}$
such that
$$
\SSS\sprime\sb{B}\;\inj\; \SSS\sb{B}\;\to\; \Gamma\sb{B} \backslash \SSS\sb{B}
$$
is surjective.
Then we make the subset $\SSS\spprime\sb{B}$ of $\SSS\sprime\sb{B}$
consisting of all $S\sb{B}\in \SSS\sprime\sb{B}$
with the following properties;
\begin{itemize}
\item
$|S\sb{B}|\sp 2 = |G\sb{B}|$,
\item
$S\sb{B}\cap G\sb n$ is trivial,
where  $S\sb{B}\st G\sb{B}$ is regarded as a subgroup of 
$G$, and 
\item
$\rootsnum\sb{R, n} (S\sb B)=\rootsnum\sb R$.
\end{itemize}
Note that these properties are invariant under the action of $\Gamma\sb{B}$.
\par\smallskip\noindent
If $\SSS\spprime\sb{B}\ne\emptyset$, then we go to the next step.
\par\smallskip\noindent
{\bf Step 3.}
Let 
$\SSS\sb{A}$ be the set of the isotropic subgroups of $(G\sb{A}, q\sb{A})$,
on which $\Gamma\sb{A}$ acts.
We  find a subset $\SSS\sprime\sb{A}$ of $\SSS\sb{A}$
such that
$$
\SSS\sprime\sb{A}\;\inj\; \SSS\sb{A}\;\to\; \Gamma\sb{A}  \backslash \SSS\sb{A}
$$
is surjective.
Then,
for each positive integer $\sigma\le 10$,
we make the subset $\SSS\spprime\sb{A} (\sigma)$ of $\SSS\sprime\sb{A}$
consisting of all $S\sb{A}\in \SSS\sprime\sb{A}$
with the following properties;
\begin{itemize}
\item
$S\sb{A}\sperp / S\sb{A}$
is a $p$-elementary group of order $p\sp{2\sigma}$, 
where $S\sb{A}\sperp$ is the
orthogonal complement of $S\sb{A}$ in $(G\sb{A}, q\sb{A})$, 
\item
$S\sb{A}\cap G\sb n$ is trivial, 
where  $S\sb{A}$ is regarded as a subgroup of 
$G$, and 
\item
$\rootsnum\sb{R, n}(S\sb A)=\rootsnum\sb R$.
\end{itemize}
Note again that these properties are invariant under the action of $\Gamma\sb{A}$.
\par\smallskip\noindent
If $\SSS\spprime\sb{A} (\sigma)\ne\emptyset$, then we go to the next step.
\par\smallskip\noindent
{\bf Step 4.}
For each   pair 
$( S\sb{A}, S\sb{B})\in \SSS\spprime\sb{A}(\sigma) \times \SSS\spprime\sb{B}$,
we make an isotropic subgroup   
$$
S:=S\sb{A}\times  S\sb{B}
$$
of $(G, q)$.
Note that $S\cap G\sb n$ is still trivial.
We check   the condition
$$
\rootsnum\sb{R, n} (S)=\rootsnum\sb R.
$$
When  $p=2$, we further check 
the condition that  the $2$-elementary lattice $\Lambda\sb S$
be  of type {\rm I}.
If we find a pair $( S\sb{A}, S\sb{B})$
satisfying these conditions, 
then we put $(R, n, \sigma)$
in the list of  $RDP$-triples
geometrically realizable in characteristic $p$.
\par
\medskip
In fact,
Algorithm I
is a special case of Algorithm II
where  $B$ is  taken to be an empty set.
\par
\subsection{Remarks} 
For many $[R, n, p]\in \RRR$,
the set 
$$
\SSS\spprime \sb{A}\;\;:=\;\;\bigcup\sb{\sigma= 1}\sp{10} \; \SSS\spprime\sb{A}(\sigma)
$$
or the set $\SSS\spprime\sb{B}$ is empty,
so that
$(R, n, \sigma)$ is not geometrically realizable in characteristic $p$ for any $\sigma$.
\par
\medskip
Let $[R, n, p]$ be an element of $\RRR$.
For a positive integer $k$, 
we denote by $\ord\sb p (k)$ the maximal integer $\nu$ such that $p\sp \nu \;|\; k$.
\begin{lemma}\label{lem:ptwo}
Let $\mu$ be the minimal non-negative integer
such that $p\sp{\mu} x=0$ holds for any $x$
in the $p$-part $(G\sb R)\sb p$ of $G\sb R$.
If $\ord\sb p (n) \ge \mu +2$, then
$\SSS\spprime\sb{A}$ is empty.
\end{lemma}
\begin{proof}
We put $\nu := \ord\sb p (n)$.
Then the $p$-part $(G\sb n)\sb p$ of $G\sb n$
is a cyclic group of order $p\sp\nu$.
Let $\eta$ be a 
generator of  $(G\sb n)\sb p$.
Suppose that $\SSS\spprime\sb{A} (\sigma)\ne\emptyset$.
Then   $\SSS\spprime\sb{\{p\}} (\sigma)$ is not empty, 
because $p\in A$.
Let $S\sb p\sperp$ be the orthogonal complement of $S\sb p$ in $(G\sb p , q\sb p)$.
Since $S\sb p \cap (G\sb n)\sb p$ is trivial,
we have $p\sp \mu y =0$ in $(G\sb n)\sb p$
for any $(x, y)\in S\sb p$,
and hence
$$
b\sb{R, n}((x, y), (0, p\sp{\nu-2}\eta))=b\sb n (p\sp{\nu-2} y, \eta)=0 \bmod{\Z}
$$
holds for any $(x, y)\in S\sb p$ by the assumption $\nu \ge \mu +2$.
Consequently,  the cyclic group
$\langle p\sp{\nu-2}\eta\rangle \subset (G\sb n)\sb p$
of order $p\sp 2 $ is contained in $S\sb p\sperp$.
Since $S\sb p \cap (G\sb n)\sb p$
is trivial,
we obtain an element of order $p\sp 2 $ in $S\sb p\sperp/ S\sb p$.
This  contradicts  the condition that $S\sb p\sperp/ S\sb p$ be $p$-elementary.
\end{proof}
\begin{lemma}\label{lem:pnotthree}
Suppose that $p\ne 3$
and $3\in B$.
Let $a\sb 2$ and $a\sb 5$ be the coefficients
of $A\sb 2$ and $A\sb 5$ in $R$,
respectively.
Then  $\SSS\spprime\sb{B}$ is empty in the following cases:
\begin{itemize}
\item[(i)] $\ord\sb  3 (|G\sb R|)=2$,
$\ord\sb 3 (n)=0$, and $a\sb 2=2$;
\item[(ii)] $\ord\sb 3 (|G\sb R|)=2$,
$\ord\sb 3 (n)=0$, and $a\sb 5=2$;
\item[(iii)] $\ord\sb  3 (|G\sb R|)=1$, 
$\ord\sb 3 (n)=1$, $n/3 = 4 \bmod 6$ and $a\sb 2=1$;
\item[(iv)] $\ord\sb 3 (|G\sb R|)=1$,
$\ord\sb 3 (n)=1$, $n/3 = 2 \bmod 6$ and $a\sb 5=1$.
\end{itemize}
\end{lemma}
\begin{proof}
We decompose $(G\sb{B}, q\sb{B})$
into the orthogonal direct sum of its $3$-part
$(G\sb 3, q\sb 3)$ and the part 
prime to $3$.
In the cases above,
$G\sb 3$ is isomorphic to 
$(\Z/(3))\sp{\oplus 2}$,
and $q\sb 3$ is given by the matrix
$$
\pm  \Bigl[\begin{matrix}
2/3 & 0 \\ 0 & 2/3\\ \end{matrix}\Bigr],
$$
because
$(q\sb{A\sb 2})\sb 3$ is given by the matrix $[-2/3]$,
$(q\sb{A\sb 5})\sb 3$ is given by $[-2^2\cdot 5/6]=[2/3]$,
and 
$(q\sb{n})\sb 3$ is given by $[m^2/n]=[m/3]$, where $m:=n/3$.
Therefore
$(G\sb 3, q\sb 3)$ does not contain any non-zero isotropic vector.
\end{proof}
By the same argument,
we can prove the following:
\begin{lemma}\label{lem:pnotseven}
Suppose that $p\ne 7$,
$\ord\sb 7 (|G\sb R|)=2$
and $\ord\sb 7 (n)=0$.
Suppose also that $a\sb 6=2$.
If $7\in B$, 
then $\SSS\spprime\sb B$ is empty.
\qed
\end{lemma}
By these lemmas,
we can remove $9247$ triples from $\RRR$
and $7722$ pairs from $\closure{\RRR}$,
before we start the calculations.
\subsection{An example}\label{subsec:anexample}
Let $(G, q)$ be the discriminant form of the lattice $Q(21 A\sb 1, 2)$.
Then $(G, q)$ is naturally isomorphic to
the vector space $\F\sb 2 \sp{21}\oplus \F\sb 2$
with the quadratic form
$$
(x\sb 1, \dots, x\sb{21}, y)\;\mapsto\; 
\frac{1}{2}\Bigl(\;-\sum\sb{i=1}\sp{21} x\sb i\sp 2 + y\sp 2\;\Bigr)\;\in\; \Q/2\Z.
$$
In Table~\ref{table:codes1},
we present,
for each $\sigma$,
an  example of an isotropic subgroup $S\sb{\sigma}$ of  $(G, q)$
that yields an even overlattice of $Q(21 A\sb 1, 2)$
isomorphic to $\Lambda\sb{2, \sigma}$.
A vector 
$$
(x\sb 1, \dots, x\sb{21}, y)\in G,\quad\textrm{where $x\sb i$ and  $y$ are $0$ or $1$, }
$$
is expressed by an integer
$$
2\sp{21} x\sb 1+ 2\sp{20} x\sb 2 + \cdots + 2\sp 2 x\sb{20} + 2 x\sb{21} + y
$$
in Table~\ref{table:codes1}.
\begin{table}
\caption{Codes for supersingular $K3$ surfaces with $21$ ordinary  nodes}
\label{table:codes1}
\begin{center}
{\small
\renewcommand{\arraystretch}{1.4}
\begin{tabular}{|c||l|}
\hline
$\sigma$ & Generators of $S\sb{\sigma}\st G$  \\ \hline
$1$ & $2097406, 1050398, 527206, 265642, 134866, 65657, 34069, 18113, 10633, 6693$ \\ \hline          
$2$ & $2097183, 1048803, 525093, 263497, 133521, 69123, 37457, 21637, 14377$ \\ \hline 
$3$ & $2097406, 1048607, 525091, 263493, 132745, 67985, 37457, 25649$ \\ \hline 
$4$ & $2097183, 1048803, 525093, 263497, 137257, 75845, 51459$ \\ \hline 
$5$ & $2097406, 1048607, 525091, 265253, 143401, 114737$ \\ \hline 
$6$ & $2097183, 1048803, 526083, 276483, 245763$ \\ \hline 
$7$ & $2097406, 1050398, 538654, 507905$ \\ \hline 
$8$ & $2097406, 1081088, 1015809$ \\ \hline 
$9$ & $2101246, 2093057$ \\ \hline 
$10$ & $4194303$ \\ \hline 
\end{tabular}
}
\end{center}
\end{table}
\section{The algorithm for the list of elliptic triples}\label{sec:elliptic_algorithm}
We use Proposition~\ref{prop:ell2}.
First we make the list $\EEE$ of  $ADE$-types $R$ of rank $20$
such that
the $RDP$-triple
$(R+A\sb 1, 2, \sigma)$
is geometrically realizable for some $\sigma$ and in some $p$.
This list consists of $95$ elements.
For each $R\in \EEE$,
we make the following calculations.

In the dual lattice
$$
Q(R+A\sb 1, 2)\dual = Q(R)\dual \oplus Q(A\sb 1)\dual \oplus I(2)\dual,
$$
we fix two vectors
$$
h:= (0,0,2)\quand z:=(0,2,0),
$$
both of which are in $Q(R+A\sb 1, 2)\st Q(R+A\sb 1, 2)\dual$.
Let $\SSS\sb{R+A\sb 1, 2}$ be the set of isotropic subgroups of $(G\sb{R+A\sb 1, 2}, q\sb{R+A\sb 1, 2})$,
and let $\Gamma\sb{R+A\sb 1, 2}$ be the image of the natural homomorphism 
$$
\Aut (Q(R))\times \{\Id\sb{Q(A\sb 1)}\}\times \{\Id\sb{I(2)}\}\;\to\; 
\Aut(G\sb{R+A\sb 1, 2}, q\sb{R+A\sb 1, 2}).
$$
We find a subset $\SSS\sb{R+A\sb 1, 2}\sprime$ of $\SSS\sb{R+A\sb 1, 2}$
such that the map
$$
\SSS\sb{R+A\sb 1, 2}\sprime\;\inj\; \SSS\sb{R+A\sb 1, 2} \;\to\; \Gamma\sb{R+A\sb 1, 2} \backslash
\SSS\sb{R+A\sb 1, 2}
$$
is surjective.
For each $S\in \SSS\sprime\sb{R+A\sb 1, 2}$,
we check the conditions
in \Cor{comb3}.
If these conditions are satisfied,
we then check the condition that 
$h-z$ be divisible by $2$ in the overlattice $\Lambda\sb S$.
Suppose that  $h-z$ is divisible by $2$ in $\Lambda\sb S$.
We denote by $U\sb S$ the indefinite unimodular  sublattice
of $\Lambda\sb S$ spanned by $f:=(h-z)/2$ and $z$, and 
calculate
$$
\MW := U\sb S\sperp  / (U\sb S\sperp)\sb{\roots}.
$$
If $p=2$ or $p=3$, 
we determine the quasi-ellipticity 
by Theorem~\ref{thm:qe}; that is, we see whether $Q(R)$ is $p$-elementary or not.
Then we 
put $\ECD{R, \MW, \sigma}$ in the list.
\begin{remark}
Let $\phi: X\to \P\sp 1$ be an \emph{elliptic}  fibration on a $K3$ surface.
Then the $ADE$-types  and the Kodaira types of reducible fibers
are corresponding in the following way:
\begin{eqnarray*}
&&
A\sb 1 \leftrightarrow {\rm I}\sb 2\;\;\text{\rm or}\;\; {\rm III},
\quad
A\sb 2 \leftrightarrow {\rm{I}}\sb 3\;\;\text{\rm or}\;\; {\rm IV},
\quad
A\sb l\; (l>2) \leftrightarrow {\rm I}\sb{l+1},
\\
&&
D\sb m \leftrightarrow {\rm{I}}\sp{*}\sb{m-4},
\quad
E\sb 6 \leftrightarrow {\rm{IV}}\sp{*},
\quad
E\sb 7 \leftrightarrow {\rm{III}}\sp{*},
\quad
E\sb 8 \leftrightarrow {\rm{II}}\sp{*}.
\end{eqnarray*}
When 
 $\phi: X\to \P\sp 1$ is a \emph{quasi-elliptic}   fibration
in characteristic $p$, 
the correspondence becomes one-to-one:
\begin{eqnarray*}
&
A\sb 1 \leftrightarrow {\rm{III}},
\quad
D\sb {2m} \leftrightarrow {\rm{I}}\sp{*}\sb{2m-4},
\quad
E\sb 7 \leftrightarrow {\rm{III}}\sp{*},
\quad
E\sb 8 \leftrightarrow {\rm{II}}\sp{*}
& 
\textrm{in characteristic $2$;}\\
&
A\sb 2 \leftrightarrow  {\rm IV},
\quad
E\sb 6 \leftrightarrow {\rm{IV}}\sp{*},
\quad
E\sb 8 \leftrightarrow {\rm{II}}\sp{*}
& 
\textrm{in characteristic $3$.}
\end{eqnarray*}
Moreover, the Mordell-Weil  group $\MW\sb{\phi}$ is necessarily $p$-elementary,
and the torsion rank $r:=\dim\sb{\F\sb p}\MW\sb{\phi}$ is 
related to
the Artin invariant $\sigma$ of $X$ 
by the following formula,
which is easily derived from the isomorphism~\eqref{eq:MW}:
$$
2 (\sigma +r)=\begin{cases}
2 \sum  \nu ({\rm I}\sb{2m}\sp* ) +\nu ({\rm III}) + \nu ({\rm III}\sp *)
& \text{in characteristic $2$}, \\
\nu ({\rm IV}) + \nu ({\rm IV}\sp *)
& \text{in characteristic $3$},
\end{cases}
$$
where $\nu (\tau)$ is the number of singular fibers of type $\tau$.
See \cite{Ito92, Ito94} for the detail.
\end{remark}
\begin{remark}
In \cite{RS_char2},
it was shown that every supersingular $K3$ surface $X$
in characteristic $2$ has a quasi-elliptic  pencil.
 Table QE shows that,
if the Artin invariant of $X$ is $10$,
then
any quasi-elliptic pencil on $X$ does not have a zero section.
\end{remark}
\begin{remark}
In \cite{Sh},
it was shown that, over the complex number field,
the torsion part of the Mordell-Weil group 
of an elliptic $K3$ surface is isomorphic to one of the following abelian groups:
$$
0, [2], [3], [4], [5], [6], [7], [8], [2,2], [4,2], [6, 2], [3,3], [4,4],
$$
where $[a]=\Z/(a)$ and $[a, b]=\Z/(a)\times \Z/(b)$.
Therefore the appearance of $\ECD{4 A\sb 5, 1, [3, 6]}$ and 
$\ECD{2A\sb 9 + 2 A\sb 1, 1, [10]}$ in Table E
is a so-called pathological phenomenon in characteristic $2$.
\end{remark}
\section{Supersingular $K3$ surfaces with $21$ ordinary nodes}\label{sec:21A1}
We  prove Proposition~\ref{prop:21A1}.
First note the following proposition,
which holds in every characteristic.
\begin{proposition}[\cite{U}~Proposition~1.7]\label{prop:U}
Let $L$ be a nef line bundle with $L\sp 2=2$
on a smooth $K3$ surface $X$.
Then the complete linear system $|L|$
defines
a surjective morphism $X\to \P\sp 2$
if and only if there does \emph{not}
exist
a divisor $E$ such that $E.L=1$ and $E\sp 2 =0$.
\qed
\end{proposition}
The proof in~\cite{U}
is valid in positive characteristics
if we replace Kawamata-Viehweg vanishing theorem
by Nikulin's proposition~\cite[Proposition 0.1]{N2}.
\par
\medskip
From now on,
we will assume that the base field $k$ is of characteristic $2$.
\par
For each $\sigma=1, \dots, 10$,
we  explicitly construct from Table~\ref{table:codes1}
a pair $(\Lambda, h)$
of a lattice $\Lambda$ isomorphic to $\Lambda\sb{2, \sigma}$
and a vector $h\in \Lambda$
with $h^2=2$
such that $\Sigma (h\sperp)=21 A\sb 1$.
It can be checked by direct calculations that the set
of $u\in \Lambda$ satisfying
$u\sp 2=0$ and $uh=1$
is empty.
By  Corollary~\ref{cor:comb1} and Proposition~\ref{prop:U},
it follows that
every  supersingular $K3$ surface
$X$ has a nef line bundle $L$ with  $L^2=2$
such that $|L|$ defines a surjective morphism
$$
\map{\Phi}{X}{\P\sp 2}
$$
that decomposes into a composite
\begin{equation}\label{eq:fpi}
X\;\maprightsp{f}\;\; Y \;\maprightsp{\pi}\; \P\sp 2
\end{equation}
of a contraction $f$ of mutually disjoint twenty-one $(-2)$-curves 
and a finite morphism $\pi$ of degree $2$.
\par
It remains to show that $\pi$ is purely inseparable.
Since $h\sp 0 (X, L\sp{\otimes m})=m^2 + 2 $
for each $m> 0$,
the graded ring $\oplus \sb{m\ge 0} H\sp 0 (X, L\sp{\otimes m})$ is generated by
elements 
$$
x\sb 0, x\sb 1, x\sb 2 \in H\sp 0 (X, L), \quad w\in H\sp 0 (X, L\sp{\otimes 3}),
$$
and the relations are generated by a relation
\begin{equation}\label{eq:sextic}
w^2 + C(x\sb 0, x\sb 1, x\sb 2) w + G(x\sb 0, x\sb 1, x\sb 2)=0
\end{equation}
in degree $6$.
It is enough to show that the cubic homogeneous polynomial 
$C$  is in fact zero.
\par
We will assume that $C$ is non-zero,
and derive a contradiction.
Let $\Gamma\st \P\sp 2$
be the divisor defined by the cubic equation
$$
C (x\sb 0, x\sb 1, x\sb 2)=0.
$$
We write $\Gamma$
as $\sum \gamma\sb i \Gamma\sb i$,
where $\Gamma\sb i$'s are reduced irreducible curves
distinct to each other,
and $\gamma\sb i$'s are positive integers.
Let $Y\sprime$ be the surface defined by~\eqref{eq:sextic} in the weighted projective space
$\P (3,1,1,1)$,
and let
\begin{equation}\label{eq:fpisprime}
X\;\maprightsp{f\sprime}\; Y\sprime \;\maprightsp{\pi\sprime}\; \P\sp 2
\end{equation}
be the natural morphisms.
We have $\pi\sprime\circ f\sprime =\pi\circ f=\Phi$.
The double cover $\pi\sprime$ is \'etale
over $\P\sp 2 \sm \Gamma$.
In particular,
the singular locus $\Sing (Y\sprime)$ of $Y\sprime$ is contained in $\pi\sp{\prime-1} (\Gamma)$.
We will show that $Y\sprime $ is normal,
and hence the decomposition~\eqref{eq:fpisprime}
of $\Phi$ coincides with the Stein factorization~\eqref{eq:fpi}
of $\Phi$.
It suffices to show that
$\dim \Sing (Y\sprime) =0$
by~\cite[Corollary (3.15) in Chapter VII]{AK}.
\par
Suppose that $\dim \Sing (Y\sprime)=1$.
There exists a reduced  irreducible component $\Gamma\sb j$ of $\Gamma$
such that $\pi\sp{\prime-1} (\Gamma\sb j)$
is contained in $\Sing (Y\sprime)$.
Let $\ell$ be a general line on $\P\sp 2$.
The divisor $\tilde\ell:=\Phi\sp * (\ell)$ on $X$  is reduced and irreducible
(\cite[Proposition~2.6]{SD}), but may possibly be singular.
The reduced irreducible curve $\tilde\ell\sprime:=\pi\sp{\prime-1}(\ell)$
on $Y\sprime$
is defined by a homogeneous  equation
$$
w^2 + C\sprime (y\sb 0, y\sb 1) w + G\sprime(y\sb 0, y\sb 1)=0
$$
of degree $6$ in the weighted projective plane $\P (3, 1, 1)$.
Hence its  arithmetic genus $p\sb a (\tilde \ell\sprime)$  is $2$.
Let $Q$ be an intersection point of $\ell$ and $\Gamma\sb j$.
Since $\ell$ is general,
we can find formal local parameters
$(\xi\sb 0, \xi\sb 1)$ of $\P\sp 2$ at $Q$
such that $\ell$ is defined by $\xi\sb 0=0$
and $\Gamma\sb j$ is defined by $\xi\sb 1=0$.
Since $Y\sprime$ is singular along $\pi\sp{\prime-1} (\Gamma\sb j)$,
the surface $Y\sprime$ is defined over $k [[ \xi\sb 0, \xi\sb 1]]$
by an equation of the form
$$
\eta\sp 2 \;+\; \xi\sb 1 \sp{\nu} \;a (\xi\sb 0, \xi\sb 1)\;\eta \;+\;  
\xi\sb 1 \sp\mu \;b (\xi\sb 0, \xi\sb 1)=0
\qquad (\nu\ge 1, \;\;\mu \ge 2), 
$$
where $\eta=w+ \beta (\xi\sb 0, \xi\sb 1)$ for some  $\beta(\xi\sb 0, \xi\sb 1) \in k [[ \xi\sb 0, \xi\sb 1]]$.
Therefore the homomorphism
$\OOO\sp{\wedge}\sb{Y\sprime, Q}\to (f\sprime\sb * \OOO\sb X)\sp\wedge\sb{Y\sprime, Q}$
factors through
$$
\OOO\sp{\wedge}\sb{Y\sprime, Q} \;\to\; \OOO\sp{\wedge}\sb{Y\sprime, Q}[ \eta/\xi\sb 1].
$$
By local calculations,
it follows that the cokernel of the homomorphism
$$
\OOO\sb{\tilde \ell\sprime} \;\to\; (f\sprime|\sb{\tilde \ell})\sb * \OOO\sb{\tilde \ell}
$$
has a non-trivial torsion subsheaf
whose support is on $\pi\sp{\prime-1} (Q)$.
Thus the arithmetic genus $p\sb a (\tilde \ell)$ 
of $\tilde \ell$ is smaller that $p\sb a (\tilde \ell\sprime)=2$,
which contradicts  $p\sb a (L)=2$.
Therefore $Y\sprime$ is normal.
\par
Let $E\sb 1, \dots, E\sb{21}$ be the $(-2)$-curves that are contracted by $f: X\to Y$.
We denote by $\widetilde \Gamma\sb i$ the strict transform
of $\Gamma\sb i$ by $\Phi: X\to \P\sp 2$,
and put
$$
\widetilde \Gamma:=\sum \gamma\sb i \widetilde \Gamma\sb i.
$$
Since $\Phi$ maps each $E\sb i$ to a point of $\Gamma$,
we have
$$
\Phi\sp * (\Gamma) = \widetilde \Gamma + \sum \alpha\sb\nu E\sb \nu \; \in \; |3L|
$$
with $\alpha \sb \nu \ge 1$ for each $\nu=1, \dots, 21$,
and hence
\begin{equation}\label{eq:24}
\widetilde \Gamma\sp 2 = (3L)\sp 2 + \sum \alpha\sb\nu\sp 2 E\sb\nu\sp 2 \le -24.
\end{equation}
Note that each $\widetilde \Gamma\sb i$ is irreducible.
Let $\Delta\sb i$ be the reduced part of $\widetilde \Gamma\sb i$,
and put
$\widetilde \Gamma\sb i=\delta \sb i \Delta\sb i$,
where
$\delta\sb i =1 \;\textrm{or}\; 2$.
Since $\Delta\sb i\sp 2 \ge -2$ for each $i$, we have 
\begin{equation}\label{eq:ge}
\widetilde \Gamma\sp 2\ge -2 \sum (\gamma\sb i\delta\sb i)^2.
\end{equation}
From~\eqref{eq:24} and~\eqref{eq:ge},
it follows that only the following cases can occur:
\par
\medskip
\hbox{
\vbox{
\halign{# \hfil& #  \hfil& # \hfil\cr
Case I: & $\Gamma = 3 \Gamma\sb 1$, & $\widetilde\Gamma=6 \Delta\sb 1$, \cr
Case II-1: & $\Gamma =  2\Gamma\sb 1 + \Gamma\sb 2$, & $\widetilde\Gamma=4 \Delta\sb 1+ \Delta\sb 2$, \cr
Case II-2: & $\Gamma =  2\Gamma\sb 1 + \Gamma\sb 2$, & $\widetilde\Gamma=4 \Delta\sb 1+ 2\Delta\sb 2$, \cr
Case III: &  $\Gamma =  \Gamma\sb 1 + \Gamma\sb 2+\Gamma\sb 3$, &
$\widetilde\Gamma=2 \Delta\sb 1+ 2\Delta\sb 2+2\Delta\sb 3$, \cr
}
}
}
\par
\medskip
\noindent
where, in each cases,
the components 
$\Gamma\sb i$ are mutually distinct lines.
Suppose that $\delta\sb j=2$ for some $\Gamma\sb j$.
We choose affine coordinates $(x, y)$ on $\P\sp 2$
such that $\Gamma\sb j$ is defined by $x=0$,
and let
$$
w\sp 2 \;+\; b (x, y) \, x w \;+\; g (x, y)\;=\;0
$$
be the defining equation of the affine part of $Y$.
Since $\delta\sb  j=2$ for some $\Gamma\sb j$,
the curve defined by 
$$
w\sp 2 + g(0, y)=0
$$
is
non-reduced,
and hence  $g(x, y)$ is of the form $\gamma (y)\sp 2 +  h(x, y)\, x$.
Replacing $w$ by $w+\gamma (y)$,
we can assume that $g(x, y)$ is equal to $ h\sprime (x, y)\, x$,
where $h\sprime (x, y)$ is a polynomial of degree $\le 5$.
The points in  $\Sing (Y) \cap \pi\inv (\Gamma\sb j)$
are therefore mapped by $\pi$ to the intersection points of 
the line $\Gamma\sb j$
and the curve $h\sprime =0$.
Note that,
since $Y$ is normal,
$h\sprime $ is not divisible by $x$.
Hence 
$\Sing (Y) \cap \pi\inv (\Gamma\sb j)$
consists of at most $5$ points.
Therefore Cases I, II-2 and III cannot occur,
because the cardinality of $\Sing (Y)$ is $21$.
In Case II-1,
there are at least $16$ points on 
$\Sing (Y) \cap \pi\inv (\Gamma\sb 2)$.
Hence 
$$
\Delta\sb 2\sp 2 = \widetilde\Gamma\sb 2 \sp 2 \le 2+ 16\cdot(-2)=-30,
$$
which is obviously impossible.
\qed
\vfill\eject
\par\smallskip
%
%
%
%
\hbox{
  \vtop{\tabskip=0pt \offinterlineskip
    \halign to \colwidth {\strut\vrule#\tabskip =\ttskip plus\pttskip&#\hfil&\vrule#&\hfil#& \vrule#& #\hfil & \vrule#\tabskip=0pt\cr
    \tablerule
    &\multispan5 \; $p=19$\hfil &\cr
    \tablerule
    \tablerule
    &\hfil$R$&& \hfil$n$ && \hfil$\sigma$ &\cr
    \tablerule
    \tablerule
    &$ A\sb{18} + A\sb{3} $ &&$76 $ &&$1$ &\cr    \tablerule
  }
 }
  \vtop{\tabskip=0pt \offinterlineskip
    \halign to \colwidth {\strut\vrule#\tabskip =\ttskip plus\pttskip&#\hfil&\vrule#&\hfil#& \vrule#& #\hfil & \vrule#\tabskip=0pt\cr
    \tablerule
    &\multispan5 \; $p=19$\hfil &\cr
    \tablerule
    \tablerule
    &\hfil$R$&& \hfil$n$ && \hfil$\sigma$ &\cr
    \tablerule
    \tablerule
    &$ A\sb{18} + A\sb{2} + A\sb{1} $ &&$114 $ &&$1$ &\cr    \tablerule
  }
 }
}
%
%
\par\smallskip
%
%
%
%
\hbox{
  \vtop{\tabskip=0pt \offinterlineskip
    \halign to \colwidth {\strut\vrule#\tabskip =\ttskip plus\pttskip&#\hfil&\vrule#&\hfil#& \vrule#& #\hfil & \vrule#\tabskip=0pt\cr
    \tablerule
    &\multispan5 \; $p=17$\hfil &\cr
    \tablerule
    \tablerule
    &\hfil$R$&& \hfil$n$ && \hfil$\sigma$ &\cr
    \tablerule
    \tablerule
    &$ A\sb{16} + A\sb{5} $ &&$102 $ &&$1$ &\cr    \tablerule
  }
 }
  \vtop{\tabskip=0pt \offinterlineskip
    \halign to \colwidth {\strut\vrule#\tabskip =\ttskip plus\pttskip&#\hfil&\vrule#&\hfil#& \vrule#& #\hfil & \vrule#\tabskip=0pt\cr
    \tablerule
    &\multispan5 \; $p=17$\hfil &\cr
    \tablerule
    \tablerule
    &\hfil$R$&& \hfil$n$ && \hfil$\sigma$ &\cr
    \tablerule
    \tablerule
    &$ A\sb{16} + A\sb{4} + A\sb{1} $ &&$170 $ &&$1$ &\cr    \tablerule
    &$ A\sb{16} + A\sb{3} + A\sb{2} $ &&$204 $ &&$1$ &\cr    \tablerule
  }
 }
}
%
%
\par\smallskip
%
%
%
%
\hbox{
  \vtop{\tabskip=0pt \offinterlineskip
    \halign to \colwidth {\strut\vrule#\tabskip =\ttskip plus\pttskip&#\hfil&\vrule#&\hfil#& \vrule#& #\hfil & \vrule#\tabskip=0pt\cr
    \tablerule
    &\multispan5 \; $p=13$\hfil &\cr
    \tablerule
    \tablerule
    &\hfil$R$&& \hfil$n$ && \hfil$\sigma$ &\cr
    \tablerule
    \tablerule
    &$ E\sb{8} + A\sb{12} + A\sb{1} $ &&$26 $ &&$1$ &\cr    \tablerule
    &$ E\sb{7} + A\sb{12} + A\sb{2} $ &&$78 $ &&$1$ &\cr    \tablerule
    &$ A\sb{12} + A\sb{8} + A\sb{1} $ &&$234 $ &&$1$ &\cr    \tablerule
  }
 }
  \vtop{\tabskip=0pt \offinterlineskip
    \halign to \colwidth {\strut\vrule#\tabskip =\ttskip plus\pttskip&#\hfil&\vrule#&\hfil#& \vrule#& #\hfil & \vrule#\tabskip=0pt\cr
    \tablerule
    &\multispan5 \; $p=13$\hfil &\cr
    \tablerule
    \tablerule
    &\hfil$R$&& \hfil$n$ && \hfil$\sigma$ &\cr
    \tablerule
    \tablerule
    &$ A\sb{12} + A\sb{7} + A\sb{2} $ &&$312 $ &&$1$ &\cr    \tablerule
    &$ A\sb{12} + A\sb{6} + A\sb{3} $ &&$364 $ &&$1$ &\cr    \tablerule
    &$ A\sb{12} + A\sb{4} + A\sb{3} + A\sb{2} $ &&$780 $ &&$1$ &\cr    \tablerule
  }
 }
}
%
%
\par\smallskip
%
%
%
%
\hbox{
  \vtop{\tabskip=0pt \offinterlineskip
    \halign to \colwidth {\strut\vrule#\tabskip =\ttskip plus\pttskip&#\hfil&\vrule#&\hfil#& \vrule#& #\hfil & \vrule#\tabskip=0pt\cr
    \tablerule
    &\multispan5 \; $p=11$\hfil &\cr
    \tablerule
    \tablerule
    &\hfil$R$&& \hfil$n$ && \hfil$\sigma$ &\cr
    \tablerule
    \tablerule
    &$ E\sb{8} + A\sb{10} + A\sb{3} $ &&$44 $ &&$1$ &\cr    \tablerule
    &$ D\sb{11} + A\sb{10} $ &&$44 $ &&$1$ &\cr    \tablerule
    &$ D\sb{9} + A\sb{10} + A\sb{2} $ &&$132 $ &&$1$ &\cr    \tablerule
    &$ D\sb{7} + A\sb{10} + A\sb{4} $ &&$220 $ &&$1$ &\cr    \tablerule
    &$ A\sb{21} $ &&$22 $ &&$1$ &\cr    \tablerule
  }
 }
  \vtop{\tabskip=0pt \offinterlineskip
    \halign to \colwidth {\strut\vrule#\tabskip =\ttskip plus\pttskip&#\hfil&\vrule#&\hfil#& \vrule#& #\hfil & \vrule#\tabskip=0pt\cr
    \tablerule
    &\multispan5 \; $p=11$\hfil &\cr
    \tablerule
    \tablerule
    &\hfil$R$&& \hfil$n$ && \hfil$\sigma$ &\cr
    \tablerule
    \tablerule
    &$ A\sb{11} + A\sb{10} $ &&$132 $ &&$1$ &\cr    \tablerule
    &$  2 A\sb{10} + A\sb{1} $ &&$2 $ &&$1$ &\cr    \tablerule
    &$ A\sb{10} + A\sb{8} + A\sb{3} $ &&$396 $ &&$1$ &\cr    \tablerule
    &$ A\sb{10} + A\sb{6} + A\sb{5} $ &&$462 $ &&$1$ &\cr    \tablerule
    &$ A\sb{10} + A\sb{6} + A\sb{4} + A\sb{1} $ &&$770 $ &&$1$ &\cr    \tablerule
  }
 }
}
%
%
\par\smallskip
%
%
%
%
\hbox{
  \vtop{\tabskip=0pt \offinterlineskip
    \halign to \colwidth {\strut\vrule#\tabskip =\ttskip plus\pttskip&#\hfil&\vrule#&\hfil#& \vrule#& #\hfil & \vrule#\tabskip=0pt\cr
    \tablerule
    &\multispan5 \; $p=7$\hfil &\cr
    \tablerule
    \tablerule
    &\hfil$R$&& \hfil$n$ && \hfil$\sigma$ &\cr
    \tablerule
    \tablerule
    &$ E\sb{8} + E\sb{7} + A\sb{6} $ &&$14 $ &&$1$ &\cr    \tablerule
    &$ E\sb{8} + D\sb{7} + A\sb{6} $ &&$28 $ &&$1$ &\cr    \tablerule
    &$ E\sb{8} + A\sb{13} $ &&$14 $ &&$1$ &\cr    \tablerule
    &$ E\sb{8} + A\sb{7} + A\sb{6} $ &&$56 $ &&$1$ &\cr    \tablerule
    &$ E\sb{8} +  2 A\sb{6} + A\sb{1} $ &&$2 $ &&$1$ &\cr    \tablerule
    &$ E\sb{7} + A\sb{13} + A\sb{1} $ &&$14 $ &&$1$ &\cr    \tablerule
    &$ E\sb{7} + A\sb{8} + A\sb{6} $ &&$126 $ &&$1$ &\cr    \tablerule
    &$ E\sb{7} +  2 A\sb{6} + A\sb{2} $ &&$6 $ &&$1$ &\cr    \tablerule
    &$ E\sb{6} + A\sb{9} + A\sb{6} $ &&$210 $ &&$1$ &\cr    \tablerule
    &$ E\sb{6} +  2 A\sb{6} + A\sb{3} $ &&$12 $ &&$1$ &\cr    \tablerule
    &$ D\sb{15} + A\sb{6} $ &&$28 $ &&$1$ &\cr    \tablerule
    &$ D\sb{14} + A\sb{6} + A\sb{1} $ &&$14 $ &&$1$ &\cr    \tablerule
    &$ D\sb{12} + A\sb{6} + A\sb{3} $ &&$28 $ &&$1$ &\cr    \tablerule
    &$ D\sb{9} +  2 A\sb{6} $ &&$4 $ &&$1$ &\cr    \tablerule
    &$ D\sb{9} + A\sb{6} + A\sb{4} + A\sb{2} $ &&$420 $ &&$1$ &\cr    \tablerule
    &$ D\sb{8} + A\sb{7} + A\sb{6} $ &&$56 $ &&$1$ &\cr    \tablerule
    &$ D\sb{7} + A\sb{13} + A\sb{1} $ &&$28 $ &&$1$ &\cr    \tablerule
    &$ D\sb{5} + A\sb{10} + A\sb{6} $ &&$308 $ &&$1$ &\cr    \tablerule
  }
 }
  \vtop{\tabskip=0pt \offinterlineskip
    \halign to \colwidth {\strut\vrule#\tabskip =\ttskip plus\pttskip&#\hfil&\vrule#&\hfil#& \vrule#& #\hfil & \vrule#\tabskip=0pt\cr
    \tablerule
    &\multispan5 \; $p=7$\hfil &\cr
    \tablerule
    \tablerule
    &\hfil$R$&& \hfil$n$ && \hfil$\sigma$ &\cr
    \tablerule
    \tablerule
    &$ D\sb{5} +  2 A\sb{6} + A\sb{4} $ &&$20 $ &&$1$ &\cr    \tablerule
    &$ A\sb{20} + A\sb{1} $ &&$42 $ &&$1$ &\cr    \tablerule
    &$ A\sb{15} + A\sb{6} $ &&$112 $ &&$1$ &\cr    \tablerule
    &$ A\sb{15} + A\sb{6} $ &&$28 $ &&$1$ &\cr    \tablerule
    &$ A\sb{14} + A\sb{6} + A\sb{1} $ &&$210 $ &&$1$ &\cr    \tablerule
    &$ A\sb{13} + A\sb{8} $ &&$126 $ &&$1$ &\cr    \tablerule
    &$ A\sb{13} + A\sb{7} + A\sb{1} $ &&$56 $ &&$1$ &\cr    \tablerule
    &$ A\sb{13} + A\sb{6} + A\sb{2} $ &&$6 $ &&$1$ &\cr    \tablerule
    &$ A\sb{13} + A\sb{6} +  2 A\sb{1} $ &&$2 $ &&$1$ &\cr    \tablerule
    &$ A\sb{12} + A\sb{6} + A\sb{2} + A\sb{1} $ &&$546 $ &&$1$ &\cr    \tablerule
    &$ A\sb{11} + A\sb{6} + A\sb{4} $ &&$420 $ &&$1$ &\cr    \tablerule
    &$ A\sb{9} +  2 A\sb{6} $ &&$10 $ &&$1$ &\cr    \tablerule
    &$ A\sb{9} + A\sb{6} + A\sb{5} + A\sb{1} $ &&$210 $ &&$1$ &\cr    \tablerule
    &$ A\sb{8} + A\sb{7} + A\sb{6} $ &&$504 $ &&$1$ &\cr    \tablerule
    &$ A\sb{8} +  2 A\sb{6} + A\sb{1} $ &&$18 $ &&$1$ &\cr    \tablerule
    &$ A\sb{8} + A\sb{6} + A\sb{5} + A\sb{2} $ &&$126 $ &&$1$ &\cr    \tablerule
    &$  3 A\sb{6} + A\sb{2} + A\sb{1} $ &&$42 $ &&$1, 2$ &\cr    \tablerule
    &$  2 A\sb{6} + A\sb{5} + A\sb{4} $ &&$30 $ &&$1$ &\cr    \tablerule
  }
 }
}
%
%
\par\smallskip
%
%
%
%
\hbox{
  \vtop{\tabskip=0pt \offinterlineskip
    \halign to \colwidth {\strut\vrule#\tabskip =\ttskip plus\pttskip&#\hfil&\vrule#&\hfil#& \vrule#& #\hfil & \vrule#\tabskip=0pt\cr
    \tablerule
    &\multispan5 \; $p=5$\hfil &\cr
    \tablerule
    \tablerule
    &\hfil$R$&& \hfil$n$ && \hfil$\sigma$ &\cr
    \tablerule
    \tablerule
    &$  2 E\sb{8} + A\sb{4} + A\sb{1} $ &&$10 $ &&$1$ &\cr    \tablerule
    &$ E\sb{8} + D\sb{7} + A\sb{4} + A\sb{2} $ &&$60 $ &&$1$ &\cr    \tablerule
    &$ E\sb{8} + A\sb{9} + A\sb{4} $ &&$2 $ &&$1$ &\cr    \tablerule
    &$ E\sb{8} + A\sb{9} + A\sb{3} + A\sb{1} $ &&$20 $ &&$1$ &\cr    \tablerule
    &$ E\sb{8} + A\sb{6} + A\sb{4} + A\sb{3} $ &&$140 $ &&$1$ &\cr    \tablerule
    &$ E\sb{8} +  3 A\sb{4} + A\sb{1} $ &&$10 $ &&$1, 2$ &\cr    \tablerule
    &$ E\sb{7} + D\sb{10} + A\sb{4} $ &&$10 $ &&$1$ &\cr    \tablerule
    &$ E\sb{7} + D\sb{5} + A\sb{9} $ &&$20 $ &&$1$ &\cr    \tablerule
    &$ E\sb{7} + A\sb{14} $ &&$30 $ &&$1$ &\cr    \tablerule
    &$ E\sb{7} + A\sb{10} + A\sb{4} $ &&$110 $ &&$1$ &\cr    \tablerule
    &$ E\sb{7} + A\sb{9} + A\sb{5} $ &&$30 $ &&$1$ &\cr    \tablerule
    &$ E\sb{7} + A\sb{9} + A\sb{4} + A\sb{1} $ &&$2 $ &&$1$ &\cr    \tablerule
    &$ E\sb{6} + D\sb{6} + A\sb{9} $ &&$30 $ &&$1$ &\cr    \tablerule
    &$ D\sb{16} + A\sb{4} + A\sb{1} $ &&$10 $ &&$1$ &\cr    \tablerule
    &$ D\sb{15} + A\sb{4} + A\sb{2} $ &&$60 $ &&$1$ &\cr    \tablerule
    &$ D\sb{12} + A\sb{4} + A\sb{3} + A\sb{2} $ &&$60 $ &&$1$ &\cr    \tablerule
    &$ D\sb{11} + A\sb{9} + A\sb{1} $ &&$20 $ &&$1$ &\cr    \tablerule
    &$ D\sb{11} + A\sb{6} + A\sb{4} $ &&$140 $ &&$1$ &\cr    \tablerule
    &$ D\sb{7} + A\sb{9} + A\sb{4} + A\sb{1} $ &&$4 $ &&$1$ &\cr    \tablerule
    &$ D\sb{7} + A\sb{7} + A\sb{4} + A\sb{3} $ &&$40 $ &&$1$ &\cr    \tablerule
    &$ D\sb{7} +  3 A\sb{4} + A\sb{2} $ &&$60 $ &&$1, 2$ &\cr    \tablerule
    &$ D\sb{6} + A\sb{11} + A\sb{4} $ &&$60 $ &&$1$ &\cr    \tablerule
    &$ D\sb{6} + A\sb{9} + A\sb{6} $ &&$70 $ &&$1$ &\cr    \tablerule
    &$ D\sb{6} + A\sb{9} + A\sb{4} + A\sb{2} $ &&$6 $ &&$1$ &\cr    \tablerule
    &$ D\sb{5} + A\sb{14} + A\sb{2} $ &&$20 $ &&$1$ &\cr    \tablerule
  }
 }
  \vtop{\tabskip=0pt \offinterlineskip
    \halign to \colwidth {\strut\vrule#\tabskip =\ttskip plus\pttskip&#\hfil&\vrule#&\hfil#& \vrule#& #\hfil & \vrule#\tabskip=0pt\cr
    \tablerule
    &\multispan5 \; $p=5$\hfil &\cr
    \tablerule
    \tablerule
    &\hfil$R$&& \hfil$n$ && \hfil$\sigma$ &\cr
    \tablerule
    \tablerule
    &$ A\sb{19} + A\sb{2} $ &&$60 $ &&$1$ &\cr    \tablerule
    &$ A\sb{17} + A\sb{4} $ &&$90 $ &&$1$ &\cr    \tablerule
    &$ A\sb{17} + A\sb{4} $ &&$10 $ &&$1$ &\cr    \tablerule
    &$ A\sb{15} + A\sb{4} + A\sb{2} $ &&$240 $ &&$1$ &\cr    \tablerule
    &$ A\sb{15} + A\sb{4} + A\sb{2} $ &&$60 $ &&$1$ &\cr    \tablerule
    &$ A\sb{14} + A\sb{7} $ &&$120 $ &&$1$ &\cr    \tablerule
    &$ A\sb{14} + A\sb{5} + A\sb{2} $ &&$30 $ &&$1$ &\cr    \tablerule
    &$ A\sb{14} + A\sb{4} + A\sb{3} $ &&$12 $ &&$1$ &\cr    \tablerule
    &$ A\sb{14} + A\sb{4} + A\sb{2} + A\sb{1} $ &&$2 $ &&$1$ &\cr    \tablerule
    &$ A\sb{13} + A\sb{4} + A\sb{3} + A\sb{1} $ &&$140 $ &&$1$ &\cr    \tablerule
    &$ A\sb{12} + A\sb{9} $ &&$130 $ &&$1$ &\cr    \tablerule
    &$ A\sb{12} + A\sb{5} + A\sb{4} $ &&$390 $ &&$1$ &\cr    \tablerule
    &$ A\sb{11} + A\sb{4} +  2 A\sb{3} $ &&$60 $ &&$1$ &\cr    \tablerule
    &$ A\sb{10} + A\sb{9} + A\sb{2} $ &&$330 $ &&$1$ &\cr    \tablerule
    &$ A\sb{10} + A\sb{7} + A\sb{4} $ &&$440 $ &&$1$ &\cr    \tablerule
    &$ A\sb{9} + A\sb{8} + A\sb{4} $ &&$18 $ &&$1$ &\cr    \tablerule
    &$ A\sb{9} + A\sb{8} + A\sb{3} + A\sb{1} $ &&$180 $ &&$1$ &\cr    \tablerule
    &$ A\sb{9} + A\sb{7} + A\sb{5} $ &&$120 $ &&$1$ &\cr    \tablerule
    &$ A\sb{9} + A\sb{6} + A\sb{4} + A\sb{2} $ &&$42 $ &&$1$ &\cr    \tablerule
    &$ A\sb{9} + A\sb{5} + A\sb{4} + A\sb{3} $ &&$12 $ &&$1$ &\cr    \tablerule
    &$ A\sb{9} +  3 A\sb{4} $ &&$2 $ &&$1, 2$ &\cr    \tablerule
    &$ A\sb{9} +  2 A\sb{4} + A\sb{3} + A\sb{1} $ &&$20 $ &&$1, 2$ &\cr    \tablerule
    &$  2 A\sb{8} + A\sb{4} + A\sb{1} $ &&$90 $ &&$1$ &\cr    \tablerule
    &$ A\sb{8} + A\sb{6} + A\sb{4} + A\sb{3} $ &&$1260 $ &&$1$ &\cr    \tablerule
    &$ A\sb{6} +  3 A\sb{4} + A\sb{3} $ &&$140 $ &&$1, 2$ &\cr    \tablerule
    &$  5 A\sb{4} + A\sb{1} $ &&$10 $ &&$1, 2, 3$ &\cr    \tablerule
  }
 }
}
%
%
\par\smallskip
%
%
%
%
\hbox{
  \vtop{\tabskip=0pt \offinterlineskip
    \halign to \colwidth {\strut\vrule#\tabskip =\ttskip plus\pttskip&#\hfil&\vrule#&\hfil#& \vrule#& #\hfil & \vrule#\tabskip=0pt\cr
    \tablerule
    &\multispan5 \; $p=3$\hfil &\cr
    \tablerule
    \tablerule
    &\hfil$R$&& \hfil$n$ && \hfil$\sigma$ &\cr
    \tablerule
    \tablerule
    &$  2 E\sb{8} + A\sb{5} $ &&$6 $ &&$1$ &\cr    \tablerule
    &$  2 E\sb{8} + A\sb{3} + A\sb{2} $ &&$12 $ &&$1$ &\cr    \tablerule
    &$  2 E\sb{8} +  2 A\sb{2} + A\sb{1} $ &&$2 $ &&$1$ &\cr    \tablerule
    &$ E\sb{8} + E\sb{7} + E\sb{6} $ &&$6 $ &&$1$ &\cr    \tablerule
    &$ E\sb{8} + E\sb{7} + A\sb{4} + A\sb{2} $ &&$30 $ &&$1$ &\cr    \tablerule
    &$ E\sb{8} + E\sb{7} +  3 A\sb{2} $ &&$6 $ &&$1, 2$ &\cr    \tablerule
    &$ E\sb{8} +  2 E\sb{6} + A\sb{1} $ &&$2 $ &&$1$ &\cr    \tablerule
    &$ E\sb{8} + E\sb{6} + A\sb{6} + A\sb{1} $ &&$42 $ &&$1$ &\cr    \tablerule
    &$ E\sb{8} + E\sb{6} + A\sb{5} + A\sb{2} $ &&$6 $ &&$1, 2$ &\cr    \tablerule
    &$ E\sb{8} + E\sb{6} + A\sb{3} +  2 A\sb{2} $ &&$12 $ &&$1, 2$ &\cr    \tablerule
    &$ E\sb{8} + E\sb{6} +  3 A\sb{2} + A\sb{1} $ &&$2 $ &&$2$ &\cr    \tablerule
    &$ E\sb{8} + D\sb{11} + A\sb{2} $ &&$12 $ &&$1$ &\cr    \tablerule
    &$ E\sb{8} + D\sb{9} +  2 A\sb{2} $ &&$4 $ &&$1$ &\cr    \tablerule
    &$ E\sb{8} + D\sb{5} + A\sb{6} + A\sb{2} $ &&$84 $ &&$1$ &\cr    \tablerule
    &$ E\sb{8} + D\sb{5} + A\sb{4} +  2 A\sb{2} $ &&$20 $ &&$1$ &\cr    \tablerule
    &$ E\sb{8} + A\sb{11} + A\sb{2} $ &&$4 $ &&$1$ &\cr    \tablerule
    &$ E\sb{8} + A\sb{11} +  2 A\sb{1} $ &&$12 $ &&$1$ &\cr    \tablerule
    &$ E\sb{8} + A\sb{10} + A\sb{2} + A\sb{1} $ &&$66 $ &&$1$ &\cr    \tablerule
    &$ E\sb{8} + A\sb{9} +  2 A\sb{2} $ &&$10 $ &&$1$ &\cr    \tablerule
    &$ E\sb{8} + A\sb{8} +  2 A\sb{2} + A\sb{1} $ &&$18 $ &&$1, 2$ &\cr    \tablerule
  }
 }
  \vtop{\tabskip=0pt \offinterlineskip
    \halign to \colwidth {\strut\vrule#\tabskip =\ttskip plus\pttskip&#\hfil&\vrule#&\hfil#& \vrule#& #\hfil & \vrule#\tabskip=0pt\cr
    \tablerule
    &\multispan5 \; $p=3$\hfil &\cr
    \tablerule
    \tablerule
    &\hfil$R$&& \hfil$n$ && \hfil$\sigma$ &\cr
    \tablerule
    \tablerule
    &$ E\sb{8} + A\sb{7} + A\sb{5} + A\sb{1} $ &&$24 $ &&$1$ &\cr    \tablerule
    &$ E\sb{8} + A\sb{7} + A\sb{4} + A\sb{2} $ &&$120 $ &&$1$ &\cr    \tablerule
    &$ E\sb{8} + A\sb{6} +  3 A\sb{2} + A\sb{1} $ &&$42 $ &&$1, 2$ &\cr    \tablerule
    &$ E\sb{8} +  2 A\sb{5} + A\sb{3} $ &&$4 $ &&$1$ &\cr    \tablerule
    &$ E\sb{8} + A\sb{5} + A\sb{4} +  2 A\sb{2} $ &&$30 $ &&$1, 2$ &\cr    \tablerule
    &$ E\sb{8} + A\sb{5} +  4 A\sb{2} $ &&$6 $ &&$1, 2, 3$ &\cr    \tablerule
    &$ E\sb{8} + A\sb{3} +  5 A\sb{2} $ &&$12 $ &&$2, 3$ &\cr    \tablerule
    &$ E\sb{8} +  6 A\sb{2} + A\sb{1} $ &&$2 $ &&$2, 3$ &\cr    \tablerule
    &$  2 E\sb{7} + D\sb{5} + A\sb{2} $ &&$12 $ &&$1$ &\cr    \tablerule
    &$  2 E\sb{7} + A\sb{3} +  2 A\sb{2} $ &&$4 $ &&$1$ &\cr    \tablerule
    &$ E\sb{7} +  2 E\sb{6} + A\sb{2} $ &&$6 $ &&$1, 2$ &\cr    \tablerule
    &$ E\sb{7} + E\sb{6} + A\sb{7} + A\sb{1} $ &&$24 $ &&$1$ &\cr    \tablerule
    &$ E\sb{7} + E\sb{6} + A\sb{5} + A\sb{3} $ &&$4 $ &&$1$ &\cr    \tablerule
    &$ E\sb{7} + E\sb{6} + A\sb{4} +  2 A\sb{2} $ &&$30 $ &&$1, 2$ &\cr    \tablerule
    &$ E\sb{7} + E\sb{6} +  4 A\sb{2} $ &&$6 $ &&$1, 2, 3$ &\cr    \tablerule
    &$ E\sb{7} + D\sb{10} +  2 A\sb{2} $ &&$2 $ &&$1$ &\cr    \tablerule
    &$ E\sb{7} + D\sb{8} + A\sb{5} + A\sb{1} $ &&$6 $ &&$1$ &\cr    \tablerule
    &$ E\sb{7} + D\sb{5} + A\sb{5} +  2 A\sb{2} $ &&$12 $ &&$1, 2$ &\cr    \tablerule
    &$ E\sb{7} + D\sb{4} +  2 A\sb{5} $ &&$2 $ &&$1$ &\cr    \tablerule
    &$ E\sb{7} + A\sb{10} +  2 A\sb{2} $ &&$22 $ &&$1$ &\cr    \tablerule
  }
 }
}
%
%
\hbox{
  \vtop{\tabskip=0pt \offinterlineskip
    \halign to \colwidth {\strut\vrule#\tabskip =\ttskip plus\pttskip&#\hfil&\vrule#&\hfil#& \vrule#& #\hfil & \vrule#\tabskip=0pt\cr
    \tablerule
    &\multispan5 \; $p=3$\hfil &\cr
    \tablerule
    \tablerule
    &\hfil$R$&& \hfil$n$ && \hfil$\sigma$ &\cr
    \tablerule
    \tablerule
    &$ E\sb{7} + A\sb{9} +  2 A\sb{2} + A\sb{1} $ &&$10 $ &&$1$ &\cr    \tablerule
    &$ E\sb{7} + A\sb{7} +  3 A\sb{2} + A\sb{1} $ &&$24 $ &&$1, 2$ &\cr    \tablerule
    &$ E\sb{7} + A\sb{6} + A\sb{4} +  2 A\sb{2} $ &&$70 $ &&$1$ &\cr    \tablerule
    &$ E\sb{7} + A\sb{5} + A\sb{3} +  3 A\sb{2} $ &&$4 $ &&$2$ &\cr    \tablerule
    &$ E\sb{7} + A\sb{4} +  5 A\sb{2} $ &&$30 $ &&$2, 3$ &\cr    \tablerule
    &$ E\sb{7} +  7 A\sb{2} $ &&$6 $ &&$2, 3, 4$ &\cr    \tablerule
    &$  3 E\sb{6} + A\sb{3} $ &&$12 $ &&$1, 2$ &\cr    \tablerule
    &$  3 E\sb{6} + A\sb{2} + A\sb{1} $ &&$2 $ &&$1, 2$ &\cr    \tablerule
    &$  2 E\sb{6} + D\sb{9} $ &&$4 $ &&$1$ &\cr    \tablerule
    &$  2 E\sb{6} + D\sb{5} + A\sb{4} $ &&$20 $ &&$1$ &\cr    \tablerule
    &$  2 E\sb{6} + A\sb{9} $ &&$10 $ &&$1$ &\cr    \tablerule
    &$  2 E\sb{6} + A\sb{8} + A\sb{1} $ &&$18 $ &&$1, 2$ &\cr    \tablerule
    &$  2 E\sb{6} + A\sb{6} + A\sb{2} + A\sb{1} $ &&$42 $ &&$1, 2$ &\cr    \tablerule
    &$  2 E\sb{6} + A\sb{5} + A\sb{4} $ &&$30 $ &&$1, 2$ &\cr    \tablerule
    &$  2 E\sb{6} + A\sb{5} +  2 A\sb{2} $ &&$6 $ &&$1, 2, 3$ &\cr    \tablerule
    &$  2 E\sb{6} + A\sb{3} +  3 A\sb{2} $ &&$12 $ &&$1, 2, 3$ &\cr    \tablerule
    &$  2 E\sb{6} +  4 A\sb{2} + A\sb{1} $ &&$2 $ &&$1, 2, 3$ &\cr    \tablerule
    &$ E\sb{6} + D\sb{14} + A\sb{1} $ &&$6 $ &&$1$ &\cr    \tablerule
    &$ E\sb{6} + D\sb{11} +  2 A\sb{2} $ &&$12 $ &&$1, 2$ &\cr    \tablerule
    &$ E\sb{6} + D\sb{10} + A\sb{5} $ &&$2 $ &&$1$ &\cr    \tablerule
    &$ E\sb{6} + D\sb{9} +  3 A\sb{2} $ &&$4 $ &&$2$ &\cr    \tablerule
    &$ E\sb{6} + D\sb{5} + A\sb{9} + A\sb{1} $ &&$60 $ &&$1$ &\cr    \tablerule
    &$ E\sb{6} + D\sb{5} + A\sb{6} +  2 A\sb{2} $ &&$84 $ &&$1, 2$ &\cr    \tablerule
    &$ E\sb{6} + D\sb{5} +  2 A\sb{5} $ &&$12 $ &&$1, 2$ &\cr    \tablerule
    &$ E\sb{6} + D\sb{5} + A\sb{4} +  3 A\sb{2} $ &&$20 $ &&$2$ &\cr    \tablerule
    &$ E\sb{6} + A\sb{14} + A\sb{1} $ &&$10 $ &&$1$ &\cr    \tablerule
    &$ E\sb{6} + A\sb{13} +  2 A\sb{1} $ &&$42 $ &&$1$ &\cr    \tablerule
    &$ E\sb{6} + A\sb{11} +  2 A\sb{2} $ &&$4 $ &&$1, 2$ &\cr    \tablerule
    &$ E\sb{6} + A\sb{11} + A\sb{2} +  2 A\sb{1} $ &&$12 $ &&$1, 2$ &\cr    \tablerule
    &$ E\sb{6} + A\sb{10} + A\sb{5} $ &&$22 $ &&$1$ &\cr    \tablerule
    &$ E\sb{6} + A\sb{10} + A\sb{4} + A\sb{1} $ &&$330 $ &&$1$ &\cr    \tablerule
    &$ E\sb{6} + A\sb{10} +  2 A\sb{2} + A\sb{1} $ &&$66 $ &&$1, 2$ &\cr    \tablerule
    &$ E\sb{6} + A\sb{9} + A\sb{5} + A\sb{1} $ &&$10 $ &&$1$ &\cr    \tablerule
    &$ E\sb{6} + A\sb{9} +  3 A\sb{2} $ &&$10 $ &&$2$ &\cr    \tablerule
    &$ E\sb{6} + A\sb{8} +  3 A\sb{2} + A\sb{1} $ &&$18 $ &&$1, 2, 3$ &\cr    \tablerule
    &$ E\sb{6} + A\sb{7} + A\sb{5} + A\sb{2} + A\sb{1} $ &&$24 $ &&$1, 2$ &\cr    \tablerule
    &$ E\sb{6} + A\sb{7} + A\sb{4} +  2 A\sb{2} $ &&$120 $ &&$1, 2$ &\cr    \tablerule
    &$ E\sb{6} + A\sb{6} + A\sb{5} + A\sb{4} $ &&$70 $ &&$1$ &\cr    \tablerule
    &$ E\sb{6} + A\sb{6} +  4 A\sb{2} + A\sb{1} $ &&$42 $ &&$1, 2, 3$ &\cr    \tablerule
    &$ E\sb{6} +  2 A\sb{5} + A\sb{3} + A\sb{2} $ &&$4 $ &&$1, 2$ &\cr    \tablerule
    &$ E\sb{6} + A\sb{5} + A\sb{4} +  3 A\sb{2} $ &&$30 $ &&$1, 2, 3$ &\cr    \tablerule
    &$ E\sb{6} + A\sb{5} +  5 A\sb{2} $ &&$6 $ &&$1, 2, 3, 4$ &\cr    \tablerule
    &$ E\sb{6} + A\sb{3} +  6 A\sb{2} $ &&$12 $ &&$1, 2, 3, 4$ &\cr    \tablerule
    &$ E\sb{6} +  7 A\sb{2} + A\sb{1} $ &&$2 $ &&$2, 3, 4$ &\cr    \tablerule
    &$ D\sb{19} + A\sb{2} $ &&$12 $ &&$1$ &\cr    \tablerule
    &$ D\sb{17} +  2 A\sb{2} $ &&$4 $ &&$1$ &\cr    \tablerule
    &$ D\sb{16} + A\sb{5} $ &&$6 $ &&$1$ &\cr    \tablerule
    &$ D\sb{16} + A\sb{3} + A\sb{2} $ &&$12 $ &&$1$ &\cr    \tablerule
    &$ D\sb{16} +  2 A\sb{2} + A\sb{1} $ &&$2 $ &&$1$ &\cr    \tablerule
    &$ D\sb{14} + A\sb{4} + A\sb{2} + A\sb{1} $ &&$30 $ &&$1$ &\cr    \tablerule
  }
 }
  \vtop{\tabskip=0pt \offinterlineskip
    \halign to \colwidth {\strut\vrule#\tabskip =\ttskip plus\pttskip&#\hfil&\vrule#&\hfil#& \vrule#& #\hfil & \vrule#\tabskip=0pt\cr
    \tablerule
    &\multispan5 \; $p=3$\hfil &\cr
    \tablerule
    \tablerule
    &\hfil$R$&& \hfil$n$ && \hfil$\sigma$ &\cr
    \tablerule
    \tablerule
    &$ D\sb{14} +  3 A\sb{2} + A\sb{1} $ &&$6 $ &&$1, 2$ &\cr    \tablerule
    &$ D\sb{13} + A\sb{6} + A\sb{2} $ &&$84 $ &&$1$ &\cr    \tablerule
    &$ D\sb{13} + A\sb{4} +  2 A\sb{2} $ &&$20 $ &&$1$ &\cr    \tablerule
    &$ D\sb{12} + D\sb{7} + A\sb{2} $ &&$12 $ &&$1$ &\cr    \tablerule
    &$ D\sb{12} + D\sb{5} +  2 A\sb{2} $ &&$4 $ &&$1$ &\cr    \tablerule
    &$ D\sb{11} +  2 A\sb{5} $ &&$4 $ &&$1$ &\cr    \tablerule
    &$ D\sb{11} +  5 A\sb{2} $ &&$12 $ &&$2, 3$ &\cr    \tablerule
    &$ D\sb{10} + D\sb{6} + A\sb{5} $ &&$6 $ &&$1$ &\cr    \tablerule
    &$ D\sb{10} + A\sb{11} $ &&$12 $ &&$1$ &\cr    \tablerule
    &$ D\sb{10} + A\sb{6} + A\sb{5} $ &&$42 $ &&$1$ &\cr    \tablerule
    &$ D\sb{10} +  2 A\sb{5} + A\sb{1} $ &&$2 $ &&$1$ &\cr    \tablerule
    &$ D\sb{10} + A\sb{5} +  3 A\sb{2} $ &&$2 $ &&$2$ &\cr    \tablerule
    &$ D\sb{9} +  6 A\sb{2} $ &&$4 $ &&$2, 3$ &\cr    \tablerule
    &$  2 D\sb{8} + A\sb{3} + A\sb{2} $ &&$12 $ &&$1$ &\cr    \tablerule
    &$ D\sb{8} + A\sb{11} + A\sb{2} $ &&$4 $ &&$1$ &\cr    \tablerule
    &$ D\sb{8} + A\sb{7} + A\sb{5} + A\sb{1} $ &&$24 $ &&$1$ &\cr    \tablerule
    &$ D\sb{8} + A\sb{7} + A\sb{4} + A\sb{2} $ &&$120 $ &&$1$ &\cr    \tablerule
    &$ D\sb{8} +  2 A\sb{5} + A\sb{2} + A\sb{1} $ &&$6 $ &&$1, 2$ &\cr    \tablerule
    &$ D\sb{7} + A\sb{14} $ &&$60 $ &&$1$ &\cr    \tablerule
    &$ D\sb{7} + A\sb{12} + A\sb{2} $ &&$156 $ &&$1$ &\cr    \tablerule
    &$ D\sb{7} + A\sb{11} + A\sb{3} $ &&$12 $ &&$1$ &\cr    \tablerule
    &$ D\sb{7} + A\sb{11} + A\sb{2} + A\sb{1} $ &&$2 $ &&$1$ &\cr    \tablerule
    &$ D\sb{7} + A\sb{9} + A\sb{5} $ &&$60 $ &&$1$ &\cr    \tablerule
    &$ D\sb{7} +  2 A\sb{5} + A\sb{4} $ &&$20 $ &&$1$ &\cr    \tablerule
    &$ D\sb{6} + D\sb{5} +  2 A\sb{5} $ &&$4 $ &&$1$ &\cr    \tablerule
    &$  2 D\sb{5} + A\sb{7} +  2 A\sb{2} $ &&$8 $ &&$1$ &\cr    \tablerule
    &$ D\sb{5} + A\sb{13} + A\sb{2} + A\sb{1} $ &&$84 $ &&$1$ &\cr    \tablerule
    &$ D\sb{5} + A\sb{12} +  2 A\sb{2} $ &&$52 $ &&$1$ &\cr    \tablerule
    &$ D\sb{5} + A\sb{11} + A\sb{4} + A\sb{1} $ &&$30 $ &&$1$ &\cr    \tablerule
    &$ D\sb{5} + A\sb{11} + A\sb{3} + A\sb{2} $ &&$4 $ &&$1$ &\cr    \tablerule
    &$ D\sb{5} + A\sb{11} +  2 A\sb{2} + A\sb{1} $ &&$6 $ &&$1, 2$ &\cr    \tablerule
    &$ D\sb{5} + A\sb{9} +  3 A\sb{2} + A\sb{1} $ &&$60 $ &&$1, 2$ &\cr    \tablerule
    &$ D\sb{5} + A\sb{8} + A\sb{4} +  2 A\sb{2} $ &&$180 $ &&$1, 2$ &\cr    \tablerule
    &$ D\sb{5} +  2 A\sb{7} + A\sb{2} $ &&$12 $ &&$1$ &\cr    \tablerule
    &$ D\sb{5} + A\sb{6} +  5 A\sb{2} $ &&$84 $ &&$2, 3$ &\cr    \tablerule
    &$ D\sb{5} +  2 A\sb{5} +  3 A\sb{2} $ &&$12 $ &&$1, 2, 3$ &\cr    \tablerule
    &$ D\sb{5} + A\sb{4} +  6 A\sb{2} $ &&$20 $ &&$2, 3$ &\cr    \tablerule
    &$ D\sb{4} + A\sb{15} + A\sb{2} $ &&$12 $ &&$1$ &\cr    \tablerule
    &$ D\sb{4} + A\sb{11} + A\sb{6} $ &&$84 $ &&$1$ &\cr    \tablerule
    &$ D\sb{4} + A\sb{11} + A\sb{4} + A\sb{2} $ &&$20 $ &&$1$ &\cr    \tablerule
    &$ D\sb{4} +  3 A\sb{5} + A\sb{2} $ &&$2 $ &&$1, 2$ &\cr    \tablerule
    &$ A\sb{17} +  2 A\sb{2} $ &&$2 $ &&$1$ &\cr    \tablerule
    &$ A\sb{16} +  2 A\sb{2} + A\sb{1} $ &&$34 $ &&$1$ &\cr    \tablerule
    &$ A\sb{15} +  2 A\sb{2} +  2 A\sb{1} $ &&$4 $ &&$1$ &\cr    \tablerule
    &$ A\sb{14} +  3 A\sb{2} + A\sb{1} $ &&$10 $ &&$1, 2$ &\cr    \tablerule
    &$ A\sb{13} + A\sb{4} +  2 A\sb{2} $ &&$70 $ &&$1$ &\cr    \tablerule
    &$ A\sb{13} +  3 A\sb{2} +  2 A\sb{1} $ &&$42 $ &&$1, 2$ &\cr    \tablerule
    &$ A\sb{12} + A\sb{4} +  2 A\sb{2} + A\sb{1} $ &&$130 $ &&$1$ &\cr    \tablerule
    &$ A\sb{11} +  5 A\sb{2} $ &&$4 $ &&$2, 3$ &\cr    \tablerule
    &$ A\sb{11} +  4 A\sb{2} +  2 A\sb{1} $ &&$12 $ &&$1, 2, 3$ &\cr    \tablerule
  }
 }
}
%
%
\hbox{
  \vtop{\tabskip=0pt \offinterlineskip
    \halign to \colwidth {\strut\vrule#\tabskip =\ttskip plus\pttskip&#\hfil&\vrule#&\hfil#& \vrule#& #\hfil & \vrule#\tabskip=0pt\cr
    \tablerule
    &\multispan5 \; $p=3$\hfil &\cr
    \tablerule
    \tablerule
    &\hfil$R$&& \hfil$n$ && \hfil$\sigma$ &\cr
    \tablerule
    \tablerule
    &$ A\sb{10} + A\sb{5} +  3 A\sb{2} $ &&$22 $ &&$2$ &\cr    \tablerule
    &$ A\sb{10} + A\sb{4} +  3 A\sb{2} + A\sb{1} $ &&$330 $ &&$1, 2$ &\cr    \tablerule
    &$ A\sb{10} +  5 A\sb{2} + A\sb{1} $ &&$66 $ &&$2, 3$ &\cr    \tablerule
    &$ A\sb{9} + A\sb{5} +  3 A\sb{2} + A\sb{1} $ &&$10 $ &&$2$ &\cr    \tablerule
    &$ A\sb{9} +  2 A\sb{4} +  2 A\sb{2} $ &&$10 $ &&$1$ &\cr    \tablerule
    &$ A\sb{9} +  6 A\sb{2} $ &&$10 $ &&$2, 3$ &\cr    \tablerule
    &$ A\sb{8} +  6 A\sb{2} + A\sb{1} $ &&$18 $ &&$1, 2, 3, 4$ &\cr    \tablerule
    &$ A\sb{7} + A\sb{5} +  4 A\sb{2} + A\sb{1} $ &&$24 $ &&$1, 2, 3$ &\cr    \tablerule
  }
 }
  \vtop{\tabskip=0pt \offinterlineskip
    \halign to \colwidth {\strut\vrule#\tabskip =\ttskip plus\pttskip&#\hfil&\vrule#&\hfil#& \vrule#& #\hfil & \vrule#\tabskip=0pt\cr
    \tablerule
    &\multispan5 \; $p=3$\hfil &\cr
    \tablerule
    \tablerule
    &\hfil$R$&& \hfil$n$ && \hfil$\sigma$ &\cr
    \tablerule
    \tablerule
    &$ A\sb{7} + A\sb{4} +  5 A\sb{2} $ &&$120 $ &&$2, 3$ &\cr    \tablerule
    &$ A\sb{6} + A\sb{5} + A\sb{4} +  3 A\sb{2} $ &&$70 $ &&$2$ &\cr    \tablerule
    &$ A\sb{6} +  7 A\sb{2} + A\sb{1} $ &&$42 $ &&$2, 3, 4$ &\cr    \tablerule
    &$  2 A\sb{5} + A\sb{3} +  4 A\sb{2} $ &&$4 $ &&$1, 2, 3$ &\cr    \tablerule
    &$ A\sb{5} + A\sb{4} +  6 A\sb{2} $ &&$30 $ &&$1, 2, 3, 4$ &\cr    \tablerule
    &$ A\sb{5} +  8 A\sb{2} $ &&$6 $ &&$1, 2, 3, 4, 5$ &\cr    \tablerule
    &$ A\sb{3} +  9 A\sb{2} $ &&$12 $ &&$1, 2, 3, 4, 5$ &\cr    \tablerule
    &$  10 A\sb{2} + A\sb{1} $ &&$2 $ &&$1, 2, 3, 4, 5$ &\cr    \tablerule
  }
 }
}
%
%
\par\smallskip
%
%
%
%
\hbox{
  \vtop{\tabskip=0pt \offinterlineskip
    \halign to \colwidth {\strut\vrule#\tabskip =\ttskip plus\pttskip&#\hfil&\vrule#&\hfil#& \vrule#& #\hfil & \vrule#\tabskip=0pt\cr
    \tablerule
    &\multispan5 \; $p=2$\hfil &\cr
    \tablerule
    \tablerule
    &\hfil$R$&& \hfil$n$ && \hfil$\sigma$ &\cr
    \tablerule
    \tablerule
    &$  2 E\sb{8} + D\sb{5} $ &&$4 $ &&$1$ &\cr    \tablerule
    &$  2 E\sb{8} + D\sb{4} + A\sb{1} $ &&$2 $ &&$1$ &\cr    \tablerule
    &$  2 E\sb{8} + A\sb{2} +  3 A\sb{1} $ &&$6 $ &&$1$ &\cr    \tablerule
    &$  2 E\sb{8} +  5 A\sb{1} $ &&$2 $ &&$2$ &\cr    \tablerule
    &$ E\sb{8} + E\sb{7} + D\sb{6} $ &&$2 $ &&$1$ &\cr    \tablerule
    &$ E\sb{8} + E\sb{7} + D\sb{4} + A\sb{2} $ &&$6 $ &&$1$ &\cr    \tablerule
    &$ E\sb{8} + E\sb{7} + D\sb{4} +  2 A\sb{1} $ &&$2 $ &&$2$ &\cr    \tablerule
    &$ E\sb{8} + E\sb{7} + A\sb{3} +  3 A\sb{1} $ &&$4 $ &&$1, 2$ &\cr    \tablerule
    &$ E\sb{8} + E\sb{7} + A\sb{2} +  4 A\sb{1} $ &&$6 $ &&$2$ &\cr    \tablerule
    &$ E\sb{8} + E\sb{7} +  6 A\sb{1} $ &&$2 $ &&$2, 3$ &\cr    \tablerule
    &$ E\sb{8} + E\sb{6} + D\sb{7} $ &&$12 $ &&$1$ &\cr    \tablerule
    &$ E\sb{8} + E\sb{6} + D\sb{4} + A\sb{3} $ &&$12 $ &&$1, 2$ &\cr    \tablerule
    &$ E\sb{8} + E\sb{6} + A\sb{4} +  3 A\sb{1} $ &&$30 $ &&$1$ &\cr    \tablerule
    &$ E\sb{8} + D\sb{13} $ &&$4 $ &&$1$ &\cr    \tablerule
    &$ E\sb{8} + D\sb{12} + A\sb{1} $ &&$2 $ &&$1$ &\cr    \tablerule
    &$ E\sb{8} + D\sb{10} + A\sb{2} + A\sb{1} $ &&$6 $ &&$1$ &\cr    \tablerule
    &$ E\sb{8} + D\sb{10} +  3 A\sb{1} $ &&$2 $ &&$1, 2$ &\cr    \tablerule
    &$ E\sb{8} + D\sb{9} + D\sb{4} $ &&$4 $ &&$1, 2$ &\cr    \tablerule
    &$ E\sb{8} + D\sb{9} + A\sb{4} $ &&$20 $ &&$1$ &\cr    \tablerule
    &$ E\sb{8} + D\sb{8} + D\sb{5} $ &&$4 $ &&$1, 2$ &\cr    \tablerule
    &$ E\sb{8} + D\sb{8} + D\sb{4} + A\sb{1} $ &&$2 $ &&$2$ &\cr    \tablerule
    &$ E\sb{8} + D\sb{8} + A\sb{2} +  3 A\sb{1} $ &&$6 $ &&$2$ &\cr    \tablerule
    &$ E\sb{8} + D\sb{8} +  5 A\sb{1} $ &&$2 $ &&$2, 3$ &\cr    \tablerule
    &$ E\sb{8} + D\sb{7} +  6 A\sb{1} $ &&$4 $ &&$2, 3$ &\cr    \tablerule
    &$ E\sb{8} +  2 D\sb{6} + A\sb{1} $ &&$2 $ &&$2$ &\cr    \tablerule
    &$ E\sb{8} + D\sb{6} + D\sb{4} + A\sb{2} + A\sb{1} $ &&$6 $ &&$2$ &\cr    \tablerule
    &$ E\sb{8} + D\sb{6} + D\sb{4} +  3 A\sb{1} $ &&$2 $ &&$2, 3$ &\cr    \tablerule
    &$ E\sb{8} + D\sb{6} + A\sb{3} +  4 A\sb{1} $ &&$4 $ &&$2, 3$ &\cr    \tablerule
    &$ E\sb{8} + D\sb{6} + A\sb{2} +  5 A\sb{1} $ &&$6 $ &&$2, 3$ &\cr    \tablerule
    &$ E\sb{8} + D\sb{6} +  7 A\sb{1} $ &&$2 $ &&$3, 4$ &\cr    \tablerule
    &$ E\sb{8} + D\sb{5} +  2 D\sb{4} $ &&$4 $ &&$2, 3$ &\cr    \tablerule
    &$ E\sb{8} + D\sb{5} + D\sb{4} + A\sb{4} $ &&$20 $ &&$1, 2$ &\cr    \tablerule
    &$ E\sb{8} + D\sb{5} + A\sb{5} +  3 A\sb{1} $ &&$12 $ &&$1, 2$ &\cr    \tablerule
    &$ E\sb{8} + D\sb{5} +  8 A\sb{1} $ &&$4 $ &&$3, 4$ &\cr    \tablerule
    &$ E\sb{8} +  3 D\sb{4} + A\sb{1} $ &&$2 $ &&$3$ &\cr    \tablerule
    &$ E\sb{8} +  2 D\sb{4} + A\sb{2} +  3 A\sb{1} $ &&$6 $ &&$3$ &\cr    \tablerule
    &$ E\sb{8} +  2 D\sb{4} +  5 A\sb{1} $ &&$2 $ &&$3, 4$ &\cr    \tablerule
    &$ E\sb{8} + D\sb{4} + A\sb{9} $ &&$10 $ &&$1$ &\cr    \tablerule
  }
 }
  \vtop{\tabskip=0pt \offinterlineskip
    \halign to \colwidth {\strut\vrule#\tabskip =\ttskip plus\pttskip&#\hfil&\vrule#&\hfil#& \vrule#& #\hfil & \vrule#\tabskip=0pt\cr
    \tablerule
    &\multispan5 \; $p=2$\hfil &\cr
    \tablerule
    \tablerule
    &\hfil$R$&& \hfil$n$ && \hfil$\sigma$ &\cr
    \tablerule
    \tablerule
    &$ E\sb{8} + D\sb{4} + A\sb{8} + A\sb{1} $ &&$18 $ &&$1$ &\cr    \tablerule
    &$ E\sb{8} + D\sb{4} + A\sb{6} + A\sb{2} + A\sb{1} $ &&$42 $ &&$1$ &\cr    \tablerule
    &$ E\sb{8} + D\sb{4} + A\sb{5} + A\sb{4} $ &&$30 $ &&$1$ &\cr    \tablerule
    &$ E\sb{8} + D\sb{4} + A\sb{3} +  6 A\sb{1} $ &&$4 $ &&$3, 4$ &\cr    \tablerule
    &$ E\sb{8} + D\sb{4} + A\sb{2} +  7 A\sb{1} $ &&$6 $ &&$3, 4$ &\cr    \tablerule
    &$ E\sb{8} + D\sb{4} +  9 A\sb{1} $ &&$2 $ &&$3, 4, 5$ &\cr    \tablerule
    &$ E\sb{8} + A\sb{10} +  3 A\sb{1} $ &&$22 $ &&$1$ &\cr    \tablerule
    &$ E\sb{8} + A\sb{9} +  4 A\sb{1} $ &&$10 $ &&$1, 2$ &\cr    \tablerule
    &$ E\sb{8} + A\sb{8} +  5 A\sb{1} $ &&$18 $ &&$2$ &\cr    \tablerule
    &$ E\sb{8} + A\sb{7} + A\sb{2} +  4 A\sb{1} $ &&$24 $ &&$1, 2$ &\cr    \tablerule
    &$ E\sb{8} + A\sb{6} + A\sb{4} +  3 A\sb{1} $ &&$70 $ &&$1$ &\cr    \tablerule
    &$ E\sb{8} + A\sb{6} + A\sb{2} +  5 A\sb{1} $ &&$42 $ &&$2$ &\cr    \tablerule
    &$ E\sb{8} + A\sb{5} + A\sb{4} +  4 A\sb{1} $ &&$30 $ &&$2$ &\cr    \tablerule
    &$ E\sb{8} + A\sb{5} + A\sb{3} +  5 A\sb{1} $ &&$12 $ &&$2, 3$ &\cr    \tablerule
    &$ E\sb{8} + A\sb{4} + A\sb{3} +  6 A\sb{1} $ &&$20 $ &&$2, 3$ &\cr    \tablerule
    &$ E\sb{8} + A\sb{3} +  10 A\sb{1} $ &&$4 $ &&$4, 5$ &\cr    \tablerule
    &$ E\sb{8} + A\sb{2} +  11 A\sb{1} $ &&$6 $ &&$4, 5$ &\cr    \tablerule
    &$ E\sb{8} +  13 A\sb{1} $ &&$2 $ &&$4, 5, 6$ &\cr    \tablerule
    &$  3 E\sb{7} $ &&$2 $ &&$1$ &\cr    \tablerule
    &$  2 E\sb{7} + D\sb{7} $ &&$4 $ &&$1$ &\cr    \tablerule
    &$  2 E\sb{7} + D\sb{6} + A\sb{1} $ &&$2 $ &&$1, 2$ &\cr    \tablerule
    &$  2 E\sb{7} + D\sb{5} +  2 A\sb{1} $ &&$4 $ &&$1, 2$ &\cr    \tablerule
    &$  2 E\sb{7} + D\sb{4} + A\sb{3} $ &&$4 $ &&$1, 2$ &\cr    \tablerule
    &$  2 E\sb{7} + D\sb{4} + A\sb{2} + A\sb{1} $ &&$6 $ &&$1, 2$ &\cr    \tablerule
    &$  2 E\sb{7} + D\sb{4} +  3 A\sb{1} $ &&$2 $ &&$1, 2, 3$ &\cr    \tablerule
    &$  2 E\sb{7} + A\sb{4} + A\sb{3} $ &&$20 $ &&$1$ &\cr    \tablerule
    &$  2 E\sb{7} + A\sb{3} +  4 A\sb{1} $ &&$4 $ &&$2, 3$ &\cr    \tablerule
    &$  2 E\sb{7} + A\sb{2} +  5 A\sb{1} $ &&$6 $ &&$2, 3$ &\cr    \tablerule
    &$  2 E\sb{7} +  7 A\sb{1} $ &&$2 $ &&$2, 3, 4$ &\cr    \tablerule
    &$ E\sb{7} + E\sb{6} + D\sb{5} +  3 A\sb{1} $ &&$12 $ &&$1, 2$ &\cr    \tablerule
    &$ E\sb{7} + E\sb{6} + D\sb{4} + A\sb{4} $ &&$30 $ &&$1$ &\cr    \tablerule
    &$ E\sb{7} + E\sb{6} + A\sb{4} +  4 A\sb{1} $ &&$30 $ &&$2$ &\cr    \tablerule
    &$ E\sb{7} + E\sb{6} + A\sb{3} +  5 A\sb{1} $ &&$12 $ &&$2, 3$ &\cr    \tablerule
    &$ E\sb{7} + D\sb{14} $ &&$2 $ &&$1$ &\cr    \tablerule
    &$ E\sb{7} + D\sb{12} + A\sb{2} $ &&$6 $ &&$1$ &\cr    \tablerule
    &$ E\sb{7} + D\sb{12} +  2 A\sb{1} $ &&$2 $ &&$1, 2$ &\cr    \tablerule
    &$ E\sb{7} + D\sb{11} +  3 A\sb{1} $ &&$4 $ &&$1, 2$ &\cr    \tablerule
    &$ E\sb{7} + D\sb{10} + D\sb{4} $ &&$2 $ &&$1, 2$ &\cr    \tablerule
  }
 }
}
%
%
\hbox{
  \vtop{\tabskip=0pt \offinterlineskip
    \halign to \colwidth {\strut\vrule#\tabskip =\ttskip plus\pttskip&#\hfil&\vrule#&\hfil#& \vrule#& #\hfil & \vrule#\tabskip=0pt\cr
    \tablerule
    &\multispan5 \; $p=2$\hfil &\cr
    \tablerule
    \tablerule
    &\hfil$R$&& \hfil$n$ && \hfil$\sigma$ &\cr
    \tablerule
    \tablerule
    &$ E\sb{7} + D\sb{10} + A\sb{3} + A\sb{1} $ &&$4 $ &&$1, 2$ &\cr    \tablerule
    &$ E\sb{7} + D\sb{10} + A\sb{2} +  2 A\sb{1} $ &&$6 $ &&$1, 2$ &\cr    \tablerule
    &$ E\sb{7} + D\sb{10} +  4 A\sb{1} $ &&$2 $ &&$1, 2, 3$ &\cr    \tablerule
    &$ E\sb{7} + D\sb{9} + A\sb{5} $ &&$12 $ &&$1$ &\cr    \tablerule
    &$ E\sb{7} + D\sb{9} +  5 A\sb{1} $ &&$4 $ &&$2, 3$ &\cr    \tablerule
    &$ E\sb{7} + D\sb{8} + D\sb{6} $ &&$2 $ &&$1, 2$ &\cr    \tablerule
    &$ E\sb{7} + D\sb{8} + D\sb{4} + A\sb{2} $ &&$6 $ &&$2$ &\cr    \tablerule
    &$ E\sb{7} + D\sb{8} + D\sb{4} +  2 A\sb{1} $ &&$2 $ &&$1, 2, 3$ &\cr    \tablerule
    &$ E\sb{7} + D\sb{8} + A\sb{3} +  3 A\sb{1} $ &&$4 $ &&$1, 2, 3$ &\cr    \tablerule
    &$ E\sb{7} + D\sb{8} + A\sb{2} +  4 A\sb{1} $ &&$6 $ &&$1, 2, 3$ &\cr    \tablerule
    &$ E\sb{7} + D\sb{8} +  6 A\sb{1} $ &&$2 $ &&$2, 3, 4$ &\cr    \tablerule
    &$ E\sb{7} + D\sb{7} + D\sb{6} + A\sb{1} $ &&$4 $ &&$1, 2$ &\cr    \tablerule
    &$ E\sb{7} + D\sb{7} + D\sb{4} +  3 A\sb{1} $ &&$4 $ &&$2, 3$ &\cr    \tablerule
    &$ E\sb{7} + D\sb{7} + A\sb{5} +  2 A\sb{1} $ &&$12 $ &&$1, 2$ &\cr    \tablerule
    &$ E\sb{7} + D\sb{7} + A\sb{4} +  3 A\sb{1} $ &&$20 $ &&$1, 2$ &\cr    \tablerule
    &$ E\sb{7} + D\sb{7} +  7 A\sb{1} $ &&$4 $ &&$3, 4$ &\cr    \tablerule
    &$ E\sb{7} +  2 D\sb{6} + A\sb{2} $ &&$6 $ &&$1, 2$ &\cr    \tablerule
    &$ E\sb{7} +  2 D\sb{6} +  2 A\sb{1} $ &&$2 $ &&$1, 2, 3$ &\cr    \tablerule
    &$ E\sb{7} + D\sb{6} + D\sb{5} +  3 A\sb{1} $ &&$4 $ &&$1, 2, 3$ &\cr    \tablerule
    &$ E\sb{7} + D\sb{6} +  2 D\sb{4} $ &&$2 $ &&$2, 3$ &\cr    \tablerule
    &$ E\sb{7} + D\sb{6} + D\sb{4} + A\sb{3} + A\sb{1} $ &&$4 $ &&$1, 2, 3$ &\cr    \tablerule
    &$ E\sb{7} + D\sb{6} + D\sb{4} + A\sb{2} +  2 A\sb{1} $ &&$6 $ &&$2, 3$ &\cr    \tablerule
    &$ E\sb{7} + D\sb{6} + D\sb{4} +  4 A\sb{1} $ &&$2 $ &&$2, 3, 4$ &\cr    \tablerule
    &$ E\sb{7} + D\sb{6} + A\sb{8} $ &&$18 $ &&$1$ &\cr    \tablerule
    &$ E\sb{7} + D\sb{6} + A\sb{6} + A\sb{2} $ &&$42 $ &&$1$ &\cr    \tablerule
    &$ E\sb{7} + D\sb{6} + A\sb{5} + A\sb{3} $ &&$12 $ &&$1, 2$ &\cr    \tablerule
    &$ E\sb{7} + D\sb{6} + A\sb{4} + A\sb{3} + A\sb{1} $ &&$20 $ &&$1, 2$ &\cr    \tablerule
    &$ E\sb{7} + D\sb{6} + A\sb{3} +  5 A\sb{1} $ &&$4 $ &&$2, 3, 4$ &\cr    \tablerule
    &$ E\sb{7} + D\sb{6} + A\sb{2} +  6 A\sb{1} $ &&$6 $ &&$2, 3, 4$ &\cr    \tablerule
    &$ E\sb{7} + D\sb{6} +  8 A\sb{1} $ &&$2 $ &&$2, 3, 4, 5$ &\cr    \tablerule
    &$ E\sb{7} + D\sb{5} + D\sb{4} + A\sb{5} $ &&$12 $ &&$1, 2$ &\cr    \tablerule
    &$ E\sb{7} + D\sb{5} + D\sb{4} +  5 A\sb{1} $ &&$4 $ &&$2, 3, 4$ &\cr    \tablerule
    &$ E\sb{7} + D\sb{5} + A\sb{5} +  4 A\sb{1} $ &&$12 $ &&$2, 3$ &\cr    \tablerule
    &$ E\sb{7} + D\sb{5} + A\sb{4} +  5 A\sb{1} $ &&$20 $ &&$2, 3$ &\cr    \tablerule
    &$ E\sb{7} + D\sb{5} +  9 A\sb{1} $ &&$4 $ &&$3, 4, 5$ &\cr    \tablerule
    &$ E\sb{7} +  3 D\sb{4} + A\sb{2} $ &&$6 $ &&$3$ &\cr    \tablerule
    &$ E\sb{7} +  3 D\sb{4} +  2 A\sb{1} $ &&$2 $ &&$2, 3, 4$ &\cr    \tablerule
    &$ E\sb{7} +  2 D\sb{4} + A\sb{3} +  3 A\sb{1} $ &&$4 $ &&$2, 3, 4$ &\cr    \tablerule
    &$ E\sb{7} +  2 D\sb{4} + A\sb{2} +  4 A\sb{1} $ &&$6 $ &&$2, 3, 4$ &\cr    \tablerule
    &$ E\sb{7} +  2 D\sb{4} +  6 A\sb{1} $ &&$2 $ &&$2, 3, 4, 5$ &\cr    \tablerule
    &$ E\sb{7} + D\sb{4} + A\sb{10} $ &&$22 $ &&$1$ &\cr    \tablerule
    &$ E\sb{7} + D\sb{4} + A\sb{9} + A\sb{1} $ &&$10 $ &&$1, 2$ &\cr    \tablerule
    &$ E\sb{7} + D\sb{4} + A\sb{8} +  2 A\sb{1} $ &&$18 $ &&$2$ &\cr    \tablerule
    &$ E\sb{7} + D\sb{4} + A\sb{7} + A\sb{2} + A\sb{1} $ &&$24 $ &&$1, 2$ &\cr    \tablerule
    &$ E\sb{7} + D\sb{4} + A\sb{6} + A\sb{4} $ &&$70 $ &&$1$ &\cr    \tablerule
    &$ E\sb{7} + D\sb{4} + A\sb{6} + A\sb{2} +  2 A\sb{1} $ &&$42 $ &&$2$ &\cr    \tablerule
    &$ E\sb{7} + D\sb{4} + A\sb{5} + A\sb{4} + A\sb{1} $ &&$30 $ &&$1, 2$ &\cr    \tablerule
    &$ E\sb{7} + D\sb{4} + A\sb{5} + A\sb{3} +  2 A\sb{1} $ &&$12 $ &&$1, 2, 3$ &\cr    \tablerule
    &$ E\sb{7} + D\sb{4} + A\sb{4} + A\sb{3} +  3 A\sb{1} $ &&$20 $ &&$2, 3$ &\cr    \tablerule
    &$ E\sb{7} + D\sb{4} + A\sb{3} +  7 A\sb{1} $ &&$4 $ &&$2, 3, 4, 5$ &\cr    \tablerule
  }
 }
  \vtop{\tabskip=0pt \offinterlineskip
    \halign to \colwidth {\strut\vrule#\tabskip =\ttskip plus\pttskip&#\hfil&\vrule#&\hfil#& \vrule#& #\hfil & \vrule#\tabskip=0pt\cr
    \tablerule
    &\multispan5 \; $p=2$\hfil &\cr
    \tablerule
    \tablerule
    &\hfil$R$&& \hfil$n$ && \hfil$\sigma$ &\cr
    \tablerule
    \tablerule
    &$ E\sb{7} + D\sb{4} + A\sb{2} +  8 A\sb{1} $ &&$6 $ &&$3, 4, 5$ &\cr    \tablerule
    &$ E\sb{7} + D\sb{4} +  10 A\sb{1} $ &&$2 $ &&$3, 4, 5, 6$ &\cr    \tablerule
    &$ E\sb{7} + A\sb{11} + A\sb{3} $ &&$6 $ &&$1$ &\cr    \tablerule
    &$ E\sb{7} + A\sb{10} +  4 A\sb{1} $ &&$22 $ &&$2$ &\cr    \tablerule
    &$ E\sb{7} + A\sb{9} + A\sb{3} + A\sb{2} $ &&$60 $ &&$1$ &\cr    \tablerule
    &$ E\sb{7} + A\sb{9} +  5 A\sb{1} $ &&$10 $ &&$2, 3$ &\cr    \tablerule
    &$ E\sb{7} + A\sb{8} +  6 A\sb{1} $ &&$18 $ &&$2, 3$ &\cr    \tablerule
    &$ E\sb{7} + A\sb{7} +  2 A\sb{3} + A\sb{1} $ &&$8 $ &&$1, 2$ &\cr    \tablerule
    &$ E\sb{7} + A\sb{7} + A\sb{2} +  5 A\sb{1} $ &&$24 $ &&$2, 3$ &\cr    \tablerule
    &$ E\sb{7} + A\sb{6} + A\sb{5} + A\sb{3} $ &&$84 $ &&$1$ &\cr    \tablerule
    &$ E\sb{7} + A\sb{6} + A\sb{4} +  4 A\sb{1} $ &&$70 $ &&$2$ &\cr    \tablerule
    &$ E\sb{7} + A\sb{6} + A\sb{2} +  6 A\sb{1} $ &&$42 $ &&$2, 3$ &\cr    \tablerule
    &$ E\sb{7} + A\sb{5} + A\sb{4} +  5 A\sb{1} $ &&$30 $ &&$2, 3$ &\cr    \tablerule
    &$ E\sb{7} + A\sb{5} + A\sb{3} +  6 A\sb{1} $ &&$12 $ &&$2, 3, 4$ &\cr    \tablerule
    &$ E\sb{7} + A\sb{4} + A\sb{3} +  7 A\sb{1} $ &&$20 $ &&$3, 4$ &\cr    \tablerule
    &$ E\sb{7} + A\sb{3} +  11 A\sb{1} $ &&$4 $ &&$3, 4, 5, 6$ &\cr    \tablerule
    &$ E\sb{7} + A\sb{2} +  12 A\sb{1} $ &&$6 $ &&$3, 4, 5, 6$ &\cr    \tablerule
    &$ E\sb{7} +  14 A\sb{1} $ &&$2 $ &&$3, 4, 5, 6, 7$ &\cr    \tablerule
    &$  3 E\sb{6} +  3 A\sb{1} $ &&$6 $ &&$1$ &\cr    \tablerule
    &$  2 E\sb{6} + D\sb{4} + A\sb{5} $ &&$6 $ &&$1$ &\cr    \tablerule
    &$  2 E\sb{6} + A\sb{5} +  4 A\sb{1} $ &&$6 $ &&$2$ &\cr    \tablerule
    &$ E\sb{6} + D\sb{15} $ &&$12 $ &&$1$ &\cr    \tablerule
    &$ E\sb{6} + D\sb{12} + A\sb{3} $ &&$12 $ &&$1, 2$ &\cr    \tablerule
    &$ E\sb{6} + D\sb{11} + D\sb{4} $ &&$12 $ &&$1, 2$ &\cr    \tablerule
    &$ E\sb{6} + D\sb{10} + A\sb{4} + A\sb{1} $ &&$30 $ &&$1$ &\cr    \tablerule
    &$ E\sb{6} + D\sb{9} + A\sb{6} $ &&$84 $ &&$1$ &\cr    \tablerule
    &$ E\sb{6} + D\sb{9} +  6 A\sb{1} $ &&$12 $ &&$2, 3$ &\cr    \tablerule
    &$ E\sb{6} + D\sb{8} + D\sb{7} $ &&$12 $ &&$1, 2$ &\cr    \tablerule
    &$ E\sb{6} + D\sb{8} + D\sb{4} + A\sb{3} $ &&$12 $ &&$1, 2, 3$ &\cr    \tablerule
    &$ E\sb{6} + D\sb{8} + A\sb{4} +  3 A\sb{1} $ &&$30 $ &&$2$ &\cr    \tablerule
    &$ E\sb{6} + D\sb{7} +  2 D\sb{4} $ &&$12 $ &&$2, 3$ &\cr    \tablerule
    &$ E\sb{6} + D\sb{7} +  8 A\sb{1} $ &&$12 $ &&$3, 4$ &\cr    \tablerule
    &$ E\sb{6} + D\sb{6} + D\sb{5} +  4 A\sb{1} $ &&$12 $ &&$2, 3$ &\cr    \tablerule
    &$ E\sb{6} + D\sb{6} + D\sb{4} + A\sb{4} + A\sb{1} $ &&$30 $ &&$2$ &\cr    \tablerule
    &$ E\sb{6} + D\sb{6} + A\sb{4} +  5 A\sb{1} $ &&$30 $ &&$2, 3$ &\cr    \tablerule
    &$ E\sb{6} + D\sb{6} + A\sb{3} +  6 A\sb{1} $ &&$12 $ &&$2, 3, 4$ &\cr    \tablerule
    &$ E\sb{6} + D\sb{5} + D\sb{4} + A\sb{6} $ &&$84 $ &&$1, 2$ &\cr    \tablerule
    &$ E\sb{6} + D\sb{5} + D\sb{4} +  6 A\sb{1} $ &&$12 $ &&$3, 4$ &\cr    \tablerule
    &$ E\sb{6} + D\sb{5} +  10 A\sb{1} $ &&$12 $ &&$4, 5$ &\cr    \tablerule
    &$ E\sb{6} +  3 D\sb{4} + A\sb{3} $ &&$12 $ &&$2, 3, 4$ &\cr    \tablerule
    &$ E\sb{6} +  2 D\sb{4} + A\sb{4} +  3 A\sb{1} $ &&$30 $ &&$3$ &\cr    \tablerule
    &$ E\sb{6} + D\sb{4} + A\sb{11} $ &&$4 $ &&$1, 2$ &\cr    \tablerule
    &$ E\sb{6} + D\sb{4} + A\sb{10} + A\sb{1} $ &&$66 $ &&$1$ &\cr    \tablerule
    &$ E\sb{6} + D\sb{4} + A\sb{8} + A\sb{2} + A\sb{1} $ &&$18 $ &&$1$ &\cr    \tablerule
    &$ E\sb{6} + D\sb{4} + A\sb{7} + A\sb{4} $ &&$120 $ &&$1, 2$ &\cr    \tablerule
    &$ E\sb{6} + D\sb{4} +  2 A\sb{5} + A\sb{1} $ &&$6 $ &&$1, 2$ &\cr    \tablerule
    &$ E\sb{6} + D\sb{4} + A\sb{4} +  7 A\sb{1} $ &&$30 $ &&$3, 4$ &\cr    \tablerule
    &$ E\sb{6} + D\sb{4} + A\sb{3} +  8 A\sb{1} $ &&$12 $ &&$3, 4, 5$ &\cr    \tablerule
    &$ E\sb{6} + A\sb{15} $ &&$12 $ &&$1$ &\cr    \tablerule
    &$ E\sb{6} + A\sb{12} + A\sb{3} $ &&$156 $ &&$1$ &\cr    \tablerule
  }
 }
}
%
%
\hbox{
  \vtop{\tabskip=0pt \offinterlineskip
    \halign to \colwidth {\strut\vrule#\tabskip =\ttskip plus\pttskip&#\hfil&\vrule#&\hfil#& \vrule#& #\hfil & \vrule#\tabskip=0pt\cr
    \tablerule
    &\multispan5 \; $p=2$\hfil &\cr
    \tablerule
    \tablerule
    &\hfil$R$&& \hfil$n$ && \hfil$\sigma$ &\cr
    \tablerule
    \tablerule
    &$ E\sb{6} + A\sb{11} + A\sb{4} $ &&$20 $ &&$1$ &\cr    \tablerule
    &$ E\sb{6} + A\sb{11} + A\sb{3} + A\sb{1} $ &&$2 $ &&$1$ &\cr    \tablerule
    &$ E\sb{6} + A\sb{10} +  5 A\sb{1} $ &&$66 $ &&$2$ &\cr    \tablerule
    &$ E\sb{6} + A\sb{8} + A\sb{2} +  5 A\sb{1} $ &&$18 $ &&$2$ &\cr    \tablerule
    &$ E\sb{6} + A\sb{6} + A\sb{3} +  6 A\sb{1} $ &&$84 $ &&$2, 3$ &\cr    \tablerule
    &$ E\sb{6} +  3 A\sb{5} $ &&$2 $ &&$1$ &\cr    \tablerule
    &$ E\sb{6} +  2 A\sb{5} +  5 A\sb{1} $ &&$6 $ &&$2, 3$ &\cr    \tablerule
    &$ E\sb{6} + A\sb{4} +  11 A\sb{1} $ &&$30 $ &&$4, 5$ &\cr    \tablerule
    &$ E\sb{6} + A\sb{3} +  12 A\sb{1} $ &&$12 $ &&$4, 5, 6$ &\cr    \tablerule
    &$ D\sb{21} $ &&$4 $ &&$1$ &\cr    \tablerule
    &$ D\sb{20} + A\sb{1} $ &&$2 $ &&$1$ &\cr    \tablerule
    &$ D\sb{18} + A\sb{2} + A\sb{1} $ &&$6 $ &&$1$ &\cr    \tablerule
    &$ D\sb{18} +  3 A\sb{1} $ &&$2 $ &&$1, 2$ &\cr    \tablerule
    &$ D\sb{17} + D\sb{4} $ &&$4 $ &&$1, 2$ &\cr    \tablerule
    &$ D\sb{17} + A\sb{4} $ &&$20 $ &&$1$ &\cr    \tablerule
    &$ D\sb{16} + D\sb{5} $ &&$4 $ &&$1, 2$ &\cr    \tablerule
    &$ D\sb{16} + D\sb{4} + A\sb{1} $ &&$2 $ &&$1, 2$ &\cr    \tablerule
    &$ D\sb{16} + A\sb{2} +  3 A\sb{1} $ &&$6 $ &&$1, 2$ &\cr    \tablerule
    &$ D\sb{16} +  5 A\sb{1} $ &&$2 $ &&$2, 3$ &\cr    \tablerule
    &$ D\sb{15} +  6 A\sb{1} $ &&$4 $ &&$2, 3$ &\cr    \tablerule
    &$ D\sb{14} + D\sb{6} + A\sb{1} $ &&$2 $ &&$1, 2$ &\cr    \tablerule
    &$ D\sb{14} + D\sb{4} + A\sb{2} + A\sb{1} $ &&$6 $ &&$1, 2$ &\cr    \tablerule
    &$ D\sb{14} + D\sb{4} +  3 A\sb{1} $ &&$2 $ &&$2, 3$ &\cr    \tablerule
    &$ D\sb{14} + A\sb{3} +  4 A\sb{1} $ &&$4 $ &&$1, 2, 3$ &\cr    \tablerule
    &$ D\sb{14} + A\sb{2} +  5 A\sb{1} $ &&$6 $ &&$2, 3$ &\cr    \tablerule
    &$ D\sb{14} +  7 A\sb{1} $ &&$2 $ &&$2, 3, 4$ &\cr    \tablerule
    &$ D\sb{13} + D\sb{8} $ &&$4 $ &&$1, 2$ &\cr    \tablerule
    &$ D\sb{13} +  2 D\sb{4} $ &&$4 $ &&$2, 3$ &\cr    \tablerule
    &$ D\sb{13} + D\sb{4} + A\sb{4} $ &&$20 $ &&$1, 2$ &\cr    \tablerule
    &$ D\sb{13} + A\sb{5} +  3 A\sb{1} $ &&$12 $ &&$1, 2$ &\cr    \tablerule
    &$ D\sb{13} +  8 A\sb{1} $ &&$4 $ &&$3, 4$ &\cr    \tablerule
    &$ D\sb{12} + D\sb{9} $ &&$4 $ &&$1, 2$ &\cr    \tablerule
    &$ D\sb{12} + D\sb{8} + A\sb{1} $ &&$2 $ &&$1, 2$ &\cr    \tablerule
    &$ D\sb{12} + D\sb{6} + A\sb{2} + A\sb{1} $ &&$6 $ &&$1, 2$ &\cr    \tablerule
    &$ D\sb{12} + D\sb{6} +  3 A\sb{1} $ &&$2 $ &&$1, 2, 3$ &\cr    \tablerule
    &$ D\sb{12} + D\sb{5} + D\sb{4} $ &&$4 $ &&$1, 2, 3$ &\cr    \tablerule
    &$ D\sb{12} + D\sb{5} + A\sb{4} $ &&$20 $ &&$1, 2$ &\cr    \tablerule
    &$ D\sb{12} +  2 D\sb{4} + A\sb{1} $ &&$2 $ &&$2, 3$ &\cr    \tablerule
    &$ D\sb{12} + D\sb{4} + A\sb{2} +  3 A\sb{1} $ &&$6 $ &&$2, 3$ &\cr    \tablerule
    &$ D\sb{12} + D\sb{4} +  5 A\sb{1} $ &&$2 $ &&$2, 3, 4$ &\cr    \tablerule
    &$ D\sb{12} + A\sb{9} $ &&$10 $ &&$1$ &\cr    \tablerule
    &$ D\sb{12} + A\sb{8} + A\sb{1} $ &&$18 $ &&$1$ &\cr    \tablerule
    &$ D\sb{12} + A\sb{6} + A\sb{2} + A\sb{1} $ &&$42 $ &&$1$ &\cr    \tablerule
    &$ D\sb{12} + A\sb{5} + A\sb{4} $ &&$30 $ &&$1$ &\cr    \tablerule
    &$ D\sb{12} + A\sb{3} +  6 A\sb{1} $ &&$4 $ &&$2, 3, 4$ &\cr    \tablerule
    &$ D\sb{12} + A\sb{2} +  7 A\sb{1} $ &&$6 $ &&$2, 3, 4$ &\cr    \tablerule
    &$ D\sb{12} +  9 A\sb{1} $ &&$2 $ &&$3, 4, 5$ &\cr    \tablerule
    &$ D\sb{11} + D\sb{6} +  4 A\sb{1} $ &&$4 $ &&$2, 3$ &\cr    \tablerule
    &$ D\sb{11} + D\sb{4} +  6 A\sb{1} $ &&$4 $ &&$3, 4$ &\cr    \tablerule
    &$ D\sb{11} + A\sb{5} +  5 A\sb{1} $ &&$12 $ &&$2, 3$ &\cr    \tablerule
  }
 }
  \vtop{\tabskip=0pt \offinterlineskip
    \halign to \colwidth {\strut\vrule#\tabskip =\ttskip plus\pttskip&#\hfil&\vrule#&\hfil#& \vrule#& #\hfil & \vrule#\tabskip=0pt\cr
    \tablerule
    &\multispan5 \; $p=2$\hfil &\cr
    \tablerule
    \tablerule
    &\hfil$R$&& \hfil$n$ && \hfil$\sigma$ &\cr
    \tablerule
    \tablerule
    &$ D\sb{11} + A\sb{4} +  6 A\sb{1} $ &&$20 $ &&$2, 3$ &\cr    \tablerule
    &$ D\sb{11} +  10 A\sb{1} $ &&$4 $ &&$4, 5$ &\cr    \tablerule
    &$  2 D\sb{10} + A\sb{1} $ &&$2 $ &&$1, 2$ &\cr    \tablerule
    &$ D\sb{10} + D\sb{8} + A\sb{2} + A\sb{1} $ &&$6 $ &&$1, 2$ &\cr    \tablerule
    &$ D\sb{10} + D\sb{8} +  3 A\sb{1} $ &&$2 $ &&$1, 2, 3$ &\cr    \tablerule
    &$ D\sb{10} + D\sb{7} +  4 A\sb{1} $ &&$4 $ &&$1, 2, 3$ &\cr    \tablerule
    &$ D\sb{10} + D\sb{6} + D\sb{4} + A\sb{1} $ &&$2 $ &&$1, 2, 3$ &\cr    \tablerule
    &$ D\sb{10} + D\sb{6} + A\sb{3} +  2 A\sb{1} $ &&$4 $ &&$1, 2, 3$ &\cr    \tablerule
    &$ D\sb{10} + D\sb{6} + A\sb{2} +  3 A\sb{1} $ &&$6 $ &&$1, 2, 3$ &\cr    \tablerule
    &$ D\sb{10} + D\sb{6} +  5 A\sb{1} $ &&$2 $ &&$2, 3, 4$ &\cr    \tablerule
    &$ D\sb{10} + D\sb{5} + A\sb{5} + A\sb{1} $ &&$12 $ &&$1, 2$ &\cr    \tablerule
    &$ D\sb{10} + D\sb{5} +  6 A\sb{1} $ &&$4 $ &&$2, 3, 4$ &\cr    \tablerule
    &$ D\sb{10} +  2 D\sb{4} + A\sb{2} + A\sb{1} $ &&$6 $ &&$2, 3$ &\cr    \tablerule
    &$ D\sb{10} +  2 D\sb{4} +  3 A\sb{1} $ &&$2 $ &&$2, 3, 4$ &\cr    \tablerule
    &$ D\sb{10} + D\sb{4} + A\sb{3} +  4 A\sb{1} $ &&$4 $ &&$2, 3, 4$ &\cr    \tablerule
    &$ D\sb{10} + D\sb{4} + A\sb{2} +  5 A\sb{1} $ &&$6 $ &&$2, 3, 4$ &\cr    \tablerule
    &$ D\sb{10} + D\sb{4} +  7 A\sb{1} $ &&$2 $ &&$2, 3, 4, 5$ &\cr    \tablerule
    &$ D\sb{10} + A\sb{10} + A\sb{1} $ &&$22 $ &&$1$ &\cr    \tablerule
    &$ D\sb{10} + A\sb{9} +  2 A\sb{1} $ &&$10 $ &&$1, 2$ &\cr    \tablerule
    &$ D\sb{10} + A\sb{8} +  3 A\sb{1} $ &&$18 $ &&$1, 2$ &\cr    \tablerule
    &$ D\sb{10} + A\sb{7} + A\sb{2} +  2 A\sb{1} $ &&$24 $ &&$1, 2$ &\cr    \tablerule
    &$ D\sb{10} + A\sb{6} + A\sb{4} + A\sb{1} $ &&$70 $ &&$1$ &\cr    \tablerule
    &$ D\sb{10} + A\sb{6} + A\sb{2} +  3 A\sb{1} $ &&$42 $ &&$1, 2$ &\cr    \tablerule
    &$ D\sb{10} + A\sb{5} + A\sb{4} +  2 A\sb{1} $ &&$30 $ &&$1, 2$ &\cr    \tablerule
    &$ D\sb{10} + A\sb{5} + A\sb{3} +  3 A\sb{1} $ &&$12 $ &&$1, 2, 3$ &\cr    \tablerule
    &$ D\sb{10} + A\sb{4} + A\sb{3} +  4 A\sb{1} $ &&$20 $ &&$1, 2, 3$ &\cr    \tablerule
    &$ D\sb{10} + A\sb{3} +  8 A\sb{1} $ &&$4 $ &&$3, 4, 5$ &\cr    \tablerule
    &$ D\sb{10} + A\sb{2} +  9 A\sb{1} $ &&$6 $ &&$3, 4, 5$ &\cr    \tablerule
    &$ D\sb{10} +  11 A\sb{1} $ &&$2 $ &&$3, 4, 5, 6$ &\cr    \tablerule
    &$ D\sb{9} + D\sb{8} + D\sb{4} $ &&$4 $ &&$1, 2, 3$ &\cr    \tablerule
    &$ D\sb{9} + D\sb{8} + A\sb{4} $ &&$20 $ &&$1, 2$ &\cr    \tablerule
    &$ D\sb{9} + D\sb{6} + A\sb{5} + A\sb{1} $ &&$12 $ &&$1, 2$ &\cr    \tablerule
    &$ D\sb{9} + D\sb{6} +  6 A\sb{1} $ &&$4 $ &&$2, 3, 4$ &\cr    \tablerule
    &$ D\sb{9} + D\sb{5} + A\sb{7} $ &&$8 $ &&$1, 2$ &\cr    \tablerule
    &$ D\sb{9} +  3 D\sb{4} $ &&$4 $ &&$2, 3, 4$ &\cr    \tablerule
    &$ D\sb{9} +  2 D\sb{4} + A\sb{4} $ &&$20 $ &&$2, 3$ &\cr    \tablerule
    &$ D\sb{9} + D\sb{4} + A\sb{5} +  3 A\sb{1} $ &&$12 $ &&$2, 3$ &\cr    \tablerule
    &$ D\sb{9} + D\sb{4} +  8 A\sb{1} $ &&$4 $ &&$3, 4, 5$ &\cr    \tablerule
    &$ D\sb{9} + A\sb{12} $ &&$52 $ &&$1$ &\cr    \tablerule
    &$ D\sb{9} + A\sb{11} + A\sb{1} $ &&$6 $ &&$1$ &\cr    \tablerule
    &$ D\sb{9} + A\sb{9} + A\sb{2} + A\sb{1} $ &&$60 $ &&$1$ &\cr    \tablerule
    &$ D\sb{9} + A\sb{8} + A\sb{4} $ &&$180 $ &&$1$ &\cr    \tablerule
    &$ D\sb{9} + A\sb{5} +  7 A\sb{1} $ &&$12 $ &&$3, 4$ &\cr    \tablerule
    &$ D\sb{9} + A\sb{4} +  8 A\sb{1} $ &&$20 $ &&$3, 4$ &\cr    \tablerule
    &$ D\sb{9} +  12 A\sb{1} $ &&$4 $ &&$4, 5, 6$ &\cr    \tablerule
    &$  2 D\sb{8} + D\sb{5} $ &&$4 $ &&$1, 2, 3$ &\cr    \tablerule
    &$  2 D\sb{8} + D\sb{4} + A\sb{1} $ &&$2 $ &&$1, 2, 3$ &\cr    \tablerule
    &$  2 D\sb{8} + A\sb{2} +  3 A\sb{1} $ &&$6 $ &&$1, 2, 3$ &\cr    \tablerule
    &$  2 D\sb{8} +  5 A\sb{1} $ &&$2 $ &&$1, 2, 3, 4$ &\cr    \tablerule
    &$ D\sb{8} + D\sb{7} +  6 A\sb{1} $ &&$4 $ &&$2, 3, 4$ &\cr    \tablerule
  }
 }
}
%
%
\hbox{
  \vtop{\tabskip=0pt \offinterlineskip
    \halign to \colwidth {\strut\vrule#\tabskip =\ttskip plus\pttskip&#\hfil&\vrule#&\hfil#& \vrule#& #\hfil & \vrule#\tabskip=0pt\cr
    \tablerule
    &\multispan5 \; $p=2$\hfil &\cr
    \tablerule
    \tablerule
    &\hfil$R$&& \hfil$n$ && \hfil$\sigma$ &\cr
    \tablerule
    \tablerule
    &$ D\sb{8} +  2 D\sb{6} + A\sb{1} $ &&$2 $ &&$1, 2, 3$ &\cr    \tablerule
    &$ D\sb{8} + D\sb{6} + D\sb{4} + A\sb{2} + A\sb{1} $ &&$6 $ &&$1, 2, 3$ &\cr    \tablerule
    &$ D\sb{8} + D\sb{6} + D\sb{4} +  3 A\sb{1} $ &&$2 $ &&$1, 2, 3, 4$ &\cr    \tablerule
    &$ D\sb{8} + D\sb{6} + A\sb{3} +  4 A\sb{1} $ &&$4 $ &&$1, 2, 3, 4$ &\cr    \tablerule
    &$ D\sb{8} + D\sb{6} + A\sb{2} +  5 A\sb{1} $ &&$6 $ &&$2, 3, 4$ &\cr    \tablerule
    &$ D\sb{8} + D\sb{6} +  7 A\sb{1} $ &&$2 $ &&$2, 3, 4, 5$ &\cr    \tablerule
    &$ D\sb{8} + D\sb{5} +  2 D\sb{4} $ &&$4 $ &&$1, 2, 3, 4$ &\cr    \tablerule
    &$ D\sb{8} + D\sb{5} + D\sb{4} + A\sb{4} $ &&$20 $ &&$1, 2, 3$ &\cr    \tablerule
    &$ D\sb{8} + D\sb{5} + A\sb{5} +  3 A\sb{1} $ &&$12 $ &&$1, 2, 3$ &\cr    \tablerule
    &$ D\sb{8} + D\sb{5} +  8 A\sb{1} $ &&$4 $ &&$2, 3, 4, 5$ &\cr    \tablerule
    &$ D\sb{8} +  3 D\sb{4} + A\sb{1} $ &&$2 $ &&$2, 3, 4$ &\cr    \tablerule
    &$ D\sb{8} +  2 D\sb{4} + A\sb{2} +  3 A\sb{1} $ &&$6 $ &&$2, 3, 4$ &\cr    \tablerule
    &$ D\sb{8} +  2 D\sb{4} +  5 A\sb{1} $ &&$2 $ &&$2, 3, 4, 5$ &\cr    \tablerule
    &$ D\sb{8} + D\sb{4} + A\sb{9} $ &&$10 $ &&$2$ &\cr    \tablerule
    &$ D\sb{8} + D\sb{4} + A\sb{8} + A\sb{1} $ &&$18 $ &&$2$ &\cr    \tablerule
    &$ D\sb{8} + D\sb{4} + A\sb{6} + A\sb{2} + A\sb{1} $ &&$42 $ &&$2$ &\cr    \tablerule
    &$ D\sb{8} + D\sb{4} + A\sb{5} + A\sb{4} $ &&$30 $ &&$2$ &\cr    \tablerule
    &$ D\sb{8} + D\sb{4} + A\sb{3} +  6 A\sb{1} $ &&$4 $ &&$2, 3, 4, 5$ &\cr    \tablerule
    &$ D\sb{8} + D\sb{4} + A\sb{2} +  7 A\sb{1} $ &&$6 $ &&$2, 3, 4, 5$ &\cr    \tablerule
    &$ D\sb{8} + D\sb{4} +  9 A\sb{1} $ &&$2 $ &&$2, 3, 4, 5, 6$ &\cr    \tablerule
    &$ D\sb{8} + A\sb{10} +  3 A\sb{1} $ &&$22 $ &&$2$ &\cr    \tablerule
    &$ D\sb{8} + A\sb{9} +  4 A\sb{1} $ &&$10 $ &&$1, 2, 3$ &\cr    \tablerule
    &$ D\sb{8} + A\sb{8} +  5 A\sb{1} $ &&$18 $ &&$2, 3$ &\cr    \tablerule
    &$ D\sb{8} + A\sb{7} + A\sb{2} +  4 A\sb{1} $ &&$24 $ &&$1, 2, 3$ &\cr    \tablerule
    &$ D\sb{8} + A\sb{6} + A\sb{4} +  3 A\sb{1} $ &&$70 $ &&$2$ &\cr    \tablerule
    &$ D\sb{8} + A\sb{6} + A\sb{2} +  5 A\sb{1} $ &&$42 $ &&$2, 3$ &\cr    \tablerule
    &$ D\sb{8} + A\sb{5} + A\sb{4} +  4 A\sb{1} $ &&$30 $ &&$1, 2, 3$ &\cr    \tablerule
    &$ D\sb{8} + A\sb{5} + A\sb{3} +  5 A\sb{1} $ &&$12 $ &&$1, 2, 3, 4$ &\cr    \tablerule
    &$ D\sb{8} + A\sb{4} + A\sb{3} +  6 A\sb{1} $ &&$20 $ &&$2, 3, 4$ &\cr    \tablerule
    &$ D\sb{8} + A\sb{3} +  10 A\sb{1} $ &&$4 $ &&$3, 4, 5, 6$ &\cr    \tablerule
    &$ D\sb{8} + A\sb{2} +  11 A\sb{1} $ &&$6 $ &&$3, 4, 5, 6$ &\cr    \tablerule
    &$ D\sb{8} +  13 A\sb{1} $ &&$2 $ &&$3, 4, 5, 6, 7$ &\cr    \tablerule
    &$ D\sb{7} +  2 D\sb{6} +  2 A\sb{1} $ &&$4 $ &&$1, 2, 3$ &\cr    \tablerule
    &$ D\sb{7} + D\sb{6} + D\sb{4} +  4 A\sb{1} $ &&$4 $ &&$2, 3, 4$ &\cr    \tablerule
    &$ D\sb{7} + D\sb{6} + A\sb{5} +  3 A\sb{1} $ &&$12 $ &&$1, 2, 3$ &\cr    \tablerule
    &$ D\sb{7} + D\sb{6} + A\sb{4} +  4 A\sb{1} $ &&$20 $ &&$2, 3$ &\cr    \tablerule
    &$ D\sb{7} + D\sb{6} +  8 A\sb{1} $ &&$4 $ &&$3, 4, 5$ &\cr    \tablerule
    &$ D\sb{7} +  2 D\sb{4} +  6 A\sb{1} $ &&$4 $ &&$2, 3, 4, 5$ &\cr    \tablerule
    &$ D\sb{7} + D\sb{4} + A\sb{5} +  5 A\sb{1} $ &&$12 $ &&$2, 3, 4$ &\cr    \tablerule
    &$ D\sb{7} + D\sb{4} + A\sb{4} +  6 A\sb{1} $ &&$20 $ &&$3, 4$ &\cr    \tablerule
    &$ D\sb{7} + D\sb{4} +  10 A\sb{1} $ &&$4 $ &&$3, 4, 5, 6$ &\cr    \tablerule
    &$ D\sb{7} + A\sb{11} +  3 A\sb{1} $ &&$6 $ &&$1, 2$ &\cr    \tablerule
    &$ D\sb{7} + A\sb{9} + A\sb{2} +  3 A\sb{1} $ &&$60 $ &&$1, 2$ &\cr    \tablerule
    &$ D\sb{7} + A\sb{7} + A\sb{3} +  4 A\sb{1} $ &&$8 $ &&$1, 2, 3$ &\cr    \tablerule
    &$ D\sb{7} + A\sb{6} + A\sb{5} +  3 A\sb{1} $ &&$84 $ &&$1, 2$ &\cr    \tablerule
    &$ D\sb{7} + A\sb{5} +  9 A\sb{1} $ &&$12 $ &&$3, 4, 5$ &\cr    \tablerule
    &$ D\sb{7} + A\sb{4} +  10 A\sb{1} $ &&$20 $ &&$4, 5$ &\cr    \tablerule
    &$ D\sb{7} +  14 A\sb{1} $ &&$4 $ &&$3, 4, 5, 6, 7$ &\cr    \tablerule
    &$  3 D\sb{6} + A\sb{3} $ &&$4 $ &&$1, 2, 3$ &\cr    \tablerule
    &$  3 D\sb{6} + A\sb{2} + A\sb{1} $ &&$6 $ &&$1, 2, 3$ &\cr    \tablerule
  }
 }
  \vtop{\tabskip=0pt \offinterlineskip
    \halign to \colwidth {\strut\vrule#\tabskip =\ttskip plus\pttskip&#\hfil&\vrule#&\hfil#& \vrule#& #\hfil & \vrule#\tabskip=0pt\cr
    \tablerule
    &\multispan5 \; $p=2$\hfil &\cr
    \tablerule
    \tablerule
    &\hfil$R$&& \hfil$n$ && \hfil$\sigma$ &\cr
    \tablerule
    \tablerule
    &$  3 D\sb{6} +  3 A\sb{1} $ &&$2 $ &&$1, 2, 3, 4$ &\cr    \tablerule
    &$  2 D\sb{6} + D\sb{5} +  4 A\sb{1} $ &&$4 $ &&$1, 2, 3, 4$ &\cr    \tablerule
    &$  2 D\sb{6} +  2 D\sb{4} + A\sb{1} $ &&$2 $ &&$1, 2, 3, 4$ &\cr    \tablerule
    &$  2 D\sb{6} + D\sb{4} + A\sb{3} +  2 A\sb{1} $ &&$4 $ &&$1, 2, 3, 4$ &\cr    \tablerule
    &$  2 D\sb{6} + D\sb{4} + A\sb{2} +  3 A\sb{1} $ &&$6 $ &&$1, 2, 3, 4$ &\cr    \tablerule
    &$  2 D\sb{6} + D\sb{4} +  5 A\sb{1} $ &&$2 $ &&$1, 2, 3, 4, 5$ &\cr    \tablerule
    &$  2 D\sb{6} + A\sb{9} $ &&$10 $ &&$1, 2$ &\cr    \tablerule
    &$  2 D\sb{6} + A\sb{8} + A\sb{1} $ &&$18 $ &&$2$ &\cr    \tablerule
    &$  2 D\sb{6} + A\sb{7} + A\sb{2} $ &&$24 $ &&$1, 2$ &\cr    \tablerule
    &$  2 D\sb{6} + A\sb{6} + A\sb{2} + A\sb{1} $ &&$42 $ &&$2$ &\cr    \tablerule
    &$  2 D\sb{6} + A\sb{5} + A\sb{4} $ &&$30 $ &&$1, 2$ &\cr    \tablerule
    &$  2 D\sb{6} + A\sb{5} + A\sb{3} + A\sb{1} $ &&$12 $ &&$1, 2, 3$ &\cr    \tablerule
    &$  2 D\sb{6} + A\sb{4} + A\sb{3} +  2 A\sb{1} $ &&$20 $ &&$1, 2, 3$ &\cr    \tablerule
    &$  2 D\sb{6} + A\sb{3} +  6 A\sb{1} $ &&$4 $ &&$2, 3, 4, 5$ &\cr    \tablerule
    &$  2 D\sb{6} + A\sb{2} +  7 A\sb{1} $ &&$6 $ &&$2, 3, 4, 5$ &\cr    \tablerule
    &$  2 D\sb{6} +  9 A\sb{1} $ &&$2 $ &&$2, 3, 4, 5, 6$ &\cr    \tablerule
    &$ D\sb{6} + D\sb{5} + D\sb{4} + A\sb{5} + A\sb{1} $ &&$12 $ &&$1, 2, 3$ &\cr    \tablerule
    &$ D\sb{6} + D\sb{5} + D\sb{4} +  6 A\sb{1} $ &&$4 $ &&$2, 3, 4, 5$ &\cr    \tablerule
    &$ D\sb{6} + D\sb{5} + A\sb{5} +  5 A\sb{1} $ &&$12 $ &&$2, 3, 4$ &\cr    \tablerule
    &$ D\sb{6} + D\sb{5} + A\sb{4} +  6 A\sb{1} $ &&$20 $ &&$2, 3, 4$ &\cr    \tablerule
    &$ D\sb{6} + D\sb{5} +  10 A\sb{1} $ &&$4 $ &&$3, 4, 5, 6$ &\cr    \tablerule
    &$ D\sb{6} +  3 D\sb{4} + A\sb{2} + A\sb{1} $ &&$6 $ &&$2, 3, 4$ &\cr    \tablerule
    &$ D\sb{6} +  3 D\sb{4} +  3 A\sb{1} $ &&$2 $ &&$1, 2, 3, 4, 5$ &\cr    \tablerule
    &$ D\sb{6} +  2 D\sb{4} + A\sb{3} +  4 A\sb{1} $ &&$4 $ &&$1, 2, 3, 4, 5$ &\cr    \tablerule
    &$ D\sb{6} +  2 D\sb{4} + A\sb{2} +  5 A\sb{1} $ &&$6 $ &&$2, 3, 4, 5$ &\cr    \tablerule
    &$ D\sb{6} +  2 D\sb{4} +  7 A\sb{1} $ &&$2 $ &&$2, 3, 4, 5, 6$ &\cr    \tablerule
    &$ D\sb{6} + D\sb{4} + A\sb{10} + A\sb{1} $ &&$22 $ &&$2$ &\cr    \tablerule
    &$ D\sb{6} + D\sb{4} + A\sb{9} +  2 A\sb{1} $ &&$10 $ &&$2, 3$ &\cr    \tablerule
    &$ D\sb{6} + D\sb{4} + A\sb{8} +  3 A\sb{1} $ &&$18 $ &&$2, 3$ &\cr    \tablerule
    &$ D\sb{6} + D\sb{4} + A\sb{7} + A\sb{2} +  2 A\sb{1} $ &&$24 $ &&$2, 3$ &\cr    \tablerule
    &$ D\sb{6} + D\sb{4} + A\sb{6} + A\sb{4} + A\sb{1} $ &&$70 $ &&$2$ &\cr    \tablerule
    &$ D\sb{6} + D\sb{4} + A\sb{6} + A\sb{2} +  3 A\sb{1} $ &&$42 $ &&$2, 3$ &\cr    \tablerule
    &$ D\sb{6} + D\sb{4} + A\sb{5} + A\sb{4} +  2 A\sb{1} $ &&$30 $ &&$2, 3$ &\cr    \tablerule
    &$ D\sb{6} + D\sb{4} + A\sb{5} + A\sb{3} +  3 A\sb{1} $ &&$12 $ &&$1, 2, 3, 4$ &\cr    \tablerule
    &$ D\sb{6} + D\sb{4} + A\sb{4} + A\sb{3} +  4 A\sb{1} $ &&$20 $ &&$2, 3, 4$ &\cr    \tablerule
    &$ D\sb{6} + D\sb{4} + A\sb{3} +  8 A\sb{1} $ &&$4 $ &&$2, 3, 4, 5, 6$ &\cr    \tablerule
    &$ D\sb{6} + D\sb{4} + A\sb{2} +  9 A\sb{1} $ &&$6 $ &&$2, 3, 4, 5, 6$ &\cr    \tablerule
    &$ D\sb{6} + D\sb{4} +  11 A\sb{1} $ &&$2 $ &&$2, 3, 4, 5, 6, $ &\cr    &$  $ && &&$7$ &\cr    \tablerule
    &$ D\sb{6} + A\sb{15} $ &&$4 $ &&$1$ &\cr    \tablerule
    &$ D\sb{6} + A\sb{13} + A\sb{2} $ &&$42 $ &&$1$ &\cr    \tablerule
    &$ D\sb{6} + A\sb{11} + A\sb{3} + A\sb{1} $ &&$6 $ &&$1, 2$ &\cr    \tablerule
    &$ D\sb{6} + A\sb{11} +  2 A\sb{2} $ &&$12 $ &&$1$ &\cr    \tablerule
    &$ D\sb{6} + A\sb{10} + A\sb{5} $ &&$66 $ &&$1$ &\cr    \tablerule
    &$ D\sb{6} + A\sb{10} +  5 A\sb{1} $ &&$22 $ &&$2, 3$ &\cr    \tablerule
    &$ D\sb{6} + A\sb{9} + A\sb{3} + A\sb{2} + A\sb{1} $ &&$60 $ &&$1, 2$ &\cr    \tablerule
    &$ D\sb{6} + A\sb{9} +  6 A\sb{1} $ &&$10 $ &&$2, 3, 4$ &\cr    \tablerule
    &$ D\sb{6} + A\sb{8} + A\sb{5} + A\sb{2} $ &&$18 $ &&$1$ &\cr    \tablerule
    &$ D\sb{6} + A\sb{8} +  7 A\sb{1} $ &&$18 $ &&$3, 4$ &\cr    \tablerule
    &$ D\sb{6} + A\sb{7} +  2 A\sb{3} +  2 A\sb{1} $ &&$8 $ &&$1, 2, 3$ &\cr    \tablerule
  }
 }
}
%
%
\hbox{
  \vtop{\tabskip=0pt \offinterlineskip
    \halign to \colwidth {\strut\vrule#\tabskip =\ttskip plus\pttskip&#\hfil&\vrule#&\hfil#& \vrule#& #\hfil & \vrule#\tabskip=0pt\cr
    \tablerule
    &\multispan5 \; $p=2$\hfil &\cr
    \tablerule
    \tablerule
    &\hfil$R$&& \hfil$n$ && \hfil$\sigma$ &\cr
    \tablerule
    \tablerule
    &$ D\sb{6} + A\sb{7} + A\sb{2} +  6 A\sb{1} $ &&$24 $ &&$2, 3, 4$ &\cr    \tablerule
    &$ D\sb{6} + A\sb{6} + A\sb{5} + A\sb{3} + A\sb{1} $ &&$84 $ &&$1, 2$ &\cr    \tablerule
    &$ D\sb{6} + A\sb{6} + A\sb{4} +  5 A\sb{1} $ &&$70 $ &&$2, 3$ &\cr    \tablerule
    &$ D\sb{6} + A\sb{6} + A\sb{2} +  7 A\sb{1} $ &&$42 $ &&$3, 4$ &\cr    \tablerule
    &$ D\sb{6} +  3 A\sb{5} $ &&$6 $ &&$1, 2$ &\cr    \tablerule
    &$ D\sb{6} + A\sb{5} + A\sb{4} +  6 A\sb{1} $ &&$30 $ &&$2, 3, 4$ &\cr    \tablerule
    &$ D\sb{6} + A\sb{5} + A\sb{3} +  7 A\sb{1} $ &&$12 $ &&$2, 3, 4, 5$ &\cr    \tablerule
    &$ D\sb{6} + A\sb{4} + A\sb{3} +  8 A\sb{1} $ &&$20 $ &&$3, 4, 5$ &\cr    \tablerule
    &$ D\sb{6} + A\sb{3} +  12 A\sb{1} $ &&$4 $ &&$2, 3, 4, 5, 6, $ &\cr    &$  $ && &&$7$ &\cr    \tablerule
    &$ D\sb{6} + A\sb{2} +  13 A\sb{1} $ &&$6 $ &&$3, 4, 5, 6, 7$ &\cr    \tablerule
    &$ D\sb{6} +  15 A\sb{1} $ &&$2 $ &&$2, 3, 4, 5, 6, $ &\cr    &$  $ && &&$7, 8$ &\cr    \tablerule
    &$  2 D\sb{5} + D\sb{4} + A\sb{7} $ &&$8 $ &&$1, 2, 3$ &\cr    \tablerule
    &$ D\sb{5} +  4 D\sb{4} $ &&$4 $ &&$1, 2, 3, 4, 5$ &\cr    \tablerule
    &$ D\sb{5} +  3 D\sb{4} + A\sb{4} $ &&$20 $ &&$2, 3, 4$ &\cr    \tablerule
    &$ D\sb{5} +  2 D\sb{4} + A\sb{5} +  3 A\sb{1} $ &&$12 $ &&$2, 3, 4$ &\cr    \tablerule
    &$ D\sb{5} +  2 D\sb{4} +  8 A\sb{1} $ &&$4 $ &&$2, 3, 4, 5, 6$ &\cr    \tablerule
    &$ D\sb{5} + D\sb{4} + A\sb{12} $ &&$52 $ &&$1, 2$ &\cr    \tablerule
    &$ D\sb{5} + D\sb{4} + A\sb{11} + A\sb{1} $ &&$6 $ &&$1, 2$ &\cr    \tablerule
    &$ D\sb{5} + D\sb{4} + A\sb{9} + A\sb{2} + A\sb{1} $ &&$60 $ &&$1, 2$ &\cr    \tablerule
    &$ D\sb{5} + D\sb{4} + A\sb{8} + A\sb{4} $ &&$180 $ &&$1, 2$ &\cr    \tablerule
    &$ D\sb{5} + D\sb{4} + A\sb{5} +  7 A\sb{1} $ &&$12 $ &&$2, 3, 4, 5$ &\cr    \tablerule
    &$ D\sb{5} + D\sb{4} + A\sb{4} +  8 A\sb{1} $ &&$20 $ &&$3, 4, 5$ &\cr    \tablerule
    &$ D\sb{5} + D\sb{4} +  12 A\sb{1} $ &&$4 $ &&$3, 4, 5, 6, 7$ &\cr    \tablerule
    &$ D\sb{5} + A\sb{16} $ &&$68 $ &&$1$ &\cr    \tablerule
    &$ D\sb{5} + A\sb{15} + A\sb{1} $ &&$2 $ &&$1$ &\cr    \tablerule
    &$ D\sb{5} + A\sb{11} + A\sb{5} $ &&$2 $ &&$1$ &\cr    \tablerule
    &$ D\sb{5} + A\sb{11} +  5 A\sb{1} $ &&$6 $ &&$2, 3$ &\cr    \tablerule
    &$ D\sb{5} + A\sb{9} + A\sb{6} + A\sb{1} $ &&$140 $ &&$1$ &\cr    \tablerule
    &$ D\sb{5} + A\sb{9} + A\sb{2} +  5 A\sb{1} $ &&$60 $ &&$2, 3$ &\cr    \tablerule
    &$ D\sb{5} +  2 A\sb{8} $ &&$36 $ &&$1$ &\cr    \tablerule
    &$ D\sb{5} +  2 A\sb{7} +  2 A\sb{1} $ &&$4 $ &&$1, 2$ &\cr    \tablerule
    &$ D\sb{5} + A\sb{7} + A\sb{3} +  6 A\sb{1} $ &&$8 $ &&$2, 3, 4$ &\cr    \tablerule
    &$ D\sb{5} + A\sb{6} + A\sb{5} +  5 A\sb{1} $ &&$84 $ &&$2, 3$ &\cr    \tablerule
    &$ D\sb{5} + A\sb{5} +  11 A\sb{1} $ &&$12 $ &&$3, 4, 5, 6$ &\cr    \tablerule
    &$ D\sb{5} + A\sb{4} +  12 A\sb{1} $ &&$20 $ &&$4, 5, 6$ &\cr    \tablerule
    &$ D\sb{5} +  16 A\sb{1} $ &&$4 $ &&$2, 3, 4, 5, 6, $ &\cr    &$  $ && &&$7, 8$ &\cr    \tablerule
    &$  5 D\sb{4} + A\sb{1} $ &&$2 $ &&$1, 2, 3, 4, 5$ &\cr    \tablerule
    &$  4 D\sb{4} + A\sb{2} +  3 A\sb{1} $ &&$6 $ &&$1, 2, 3, 4, 5$ &\cr    \tablerule
    &$  4 D\sb{4} +  5 A\sb{1} $ &&$2 $ &&$1, 2, 3, 4, 5, $ &\cr    &$  $ && &&$6$ &\cr    \tablerule
    &$  3 D\sb{4} + A\sb{9} $ &&$10 $ &&$3$ &\cr    \tablerule
    &$  3 D\sb{4} + A\sb{8} + A\sb{1} $ &&$18 $ &&$3$ &\cr    \tablerule
    &$  3 D\sb{4} + A\sb{6} + A\sb{2} + A\sb{1} $ &&$42 $ &&$3$ &\cr    \tablerule
    &$  3 D\sb{4} + A\sb{5} + A\sb{4} $ &&$30 $ &&$3$ &\cr    \tablerule
    &$  3 D\sb{4} + A\sb{3} +  6 A\sb{1} $ &&$4 $ &&$1, 2, 3, 4, 5, $ &\cr    &$  $ && &&$6$ &\cr    \tablerule
    &$  3 D\sb{4} + A\sb{2} +  7 A\sb{1} $ &&$6 $ &&$2, 3, 4, 5, 6$ &\cr    \tablerule
  }
 }
  \vtop{\tabskip=0pt \offinterlineskip
    \halign to \colwidth {\strut\vrule#\tabskip =\ttskip plus\pttskip&#\hfil&\vrule#&\hfil#& \vrule#& #\hfil & \vrule#\tabskip=0pt\cr
    \tablerule
    &\multispan5 \; $p=2$\hfil &\cr
    \tablerule
    \tablerule
    &\hfil$R$&& \hfil$n$ && \hfil$\sigma$ &\cr
    \tablerule
    \tablerule
    &$  3 D\sb{4} +  9 A\sb{1} $ &&$2 $ &&$1, 2, 3, 4, 5, $ &\cr    &$  $ && &&$6, 7$ &\cr    \tablerule
    &$  2 D\sb{4} + A\sb{10} +  3 A\sb{1} $ &&$22 $ &&$3$ &\cr    \tablerule
    &$  2 D\sb{4} + A\sb{9} +  4 A\sb{1} $ &&$10 $ &&$2, 3, 4$ &\cr    \tablerule
    &$  2 D\sb{4} + A\sb{8} +  5 A\sb{1} $ &&$18 $ &&$3, 4$ &\cr    \tablerule
    &$  2 D\sb{4} + A\sb{7} + A\sb{2} +  4 A\sb{1} $ &&$24 $ &&$2, 3, 4$ &\cr    \tablerule
    &$  2 D\sb{4} + A\sb{6} + A\sb{4} +  3 A\sb{1} $ &&$70 $ &&$3$ &\cr    \tablerule
    &$  2 D\sb{4} + A\sb{6} + A\sb{2} +  5 A\sb{1} $ &&$42 $ &&$3, 4$ &\cr    \tablerule
    &$  2 D\sb{4} + A\sb{5} + A\sb{4} +  4 A\sb{1} $ &&$30 $ &&$2, 3, 4$ &\cr    \tablerule
    &$  2 D\sb{4} + A\sb{5} + A\sb{3} +  5 A\sb{1} $ &&$12 $ &&$2, 3, 4, 5$ &\cr    \tablerule
    &$  2 D\sb{4} + A\sb{4} + A\sb{3} +  6 A\sb{1} $ &&$20 $ &&$2, 3, 4, 5$ &\cr    \tablerule
    &$  2 D\sb{4} + A\sb{3} +  10 A\sb{1} $ &&$4 $ &&$2, 3, 4, 5, 6, $ &\cr    &$  $ && &&$7$ &\cr    \tablerule
    &$  2 D\sb{4} + A\sb{2} +  11 A\sb{1} $ &&$6 $ &&$2, 3, 4, 5, 6, $ &\cr    &$  $ && &&$7$ &\cr    \tablerule
    &$  2 D\sb{4} +  13 A\sb{1} $ &&$2 $ &&$2, 3, 4, 5, 6, $ &\cr    &$  $ && &&$7, 8$ &\cr    \tablerule
    &$ D\sb{4} + A\sb{17} $ &&$2 $ &&$1$ &\cr    \tablerule
    &$ D\sb{4} + A\sb{16} + A\sb{1} $ &&$34 $ &&$1$ &\cr    \tablerule
    &$ D\sb{4} + A\sb{15} +  2 A\sb{1} $ &&$4 $ &&$1, 2$ &\cr    \tablerule
    &$ D\sb{4} + A\sb{14} + A\sb{2} + A\sb{1} $ &&$10 $ &&$1$ &\cr    \tablerule
    &$ D\sb{4} + A\sb{13} + A\sb{4} $ &&$70 $ &&$1$ &\cr    \tablerule
    &$ D\sb{4} + A\sb{13} + A\sb{2} +  2 A\sb{1} $ &&$42 $ &&$1, 2$ &\cr    \tablerule
    &$ D\sb{4} + A\sb{12} + A\sb{4} + A\sb{1} $ &&$130 $ &&$1$ &\cr    \tablerule
    &$ D\sb{4} + A\sb{11} + A\sb{3} +  3 A\sb{1} $ &&$6 $ &&$2, 3$ &\cr    \tablerule
    &$ D\sb{4} + A\sb{11} +  2 A\sb{2} +  2 A\sb{1} $ &&$12 $ &&$1, 2$ &\cr    \tablerule
    &$ D\sb{4} + A\sb{10} + A\sb{5} +  2 A\sb{1} $ &&$66 $ &&$2$ &\cr    \tablerule
    &$ D\sb{4} + A\sb{10} + A\sb{4} + A\sb{2} + A\sb{1} $ &&$330 $ &&$1$ &\cr    \tablerule
    &$ D\sb{4} + A\sb{10} +  7 A\sb{1} $ &&$22 $ &&$3, 4$ &\cr    \tablerule
    &$ D\sb{4} + A\sb{9} +  2 A\sb{4} $ &&$10 $ &&$1$ &\cr    \tablerule
    &$ D\sb{4} + A\sb{9} + A\sb{3} + A\sb{2} +  3 A\sb{1} $ &&$60 $ &&$2, 3$ &\cr    \tablerule
    &$ D\sb{4} + A\sb{9} +  8 A\sb{1} $ &&$10 $ &&$3, 4, 5$ &\cr    \tablerule
    &$ D\sb{4} + A\sb{8} + A\sb{5} + A\sb{2} +  2 A\sb{1} $ &&$18 $ &&$2$ &\cr    \tablerule
    &$ D\sb{4} + A\sb{8} +  9 A\sb{1} $ &&$18 $ &&$3, 4, 5$ &\cr    \tablerule
    &$ D\sb{4} + A\sb{7} +  2 A\sb{3} +  4 A\sb{1} $ &&$8 $ &&$2, 3, 4$ &\cr    \tablerule
    &$ D\sb{4} + A\sb{7} + A\sb{2} +  8 A\sb{1} $ &&$24 $ &&$3, 4, 5$ &\cr    \tablerule
    &$ D\sb{4} + A\sb{6} + A\sb{5} + A\sb{3} +  3 A\sb{1} $ &&$84 $ &&$2, 3$ &\cr    \tablerule
    &$ D\sb{4} + A\sb{6} + A\sb{4} +  7 A\sb{1} $ &&$70 $ &&$3, 4$ &\cr    \tablerule
    &$ D\sb{4} + A\sb{6} + A\sb{2} +  9 A\sb{1} $ &&$42 $ &&$3, 4, 5$ &\cr    \tablerule
    &$ D\sb{4} +  3 A\sb{5} +  2 A\sb{1} $ &&$6 $ &&$1, 2, 3$ &\cr    \tablerule
    &$ D\sb{4} + A\sb{5} + A\sb{4} +  8 A\sb{1} $ &&$30 $ &&$3, 4, 5$ &\cr    \tablerule
    &$ D\sb{4} + A\sb{5} + A\sb{3} +  9 A\sb{1} $ &&$12 $ &&$2, 3, 4, 5, 6$ &\cr    \tablerule
    &$ D\sb{4} + A\sb{4} + A\sb{3} +  10 A\sb{1} $ &&$20 $ &&$3, 4, 5, 6$ &\cr    \tablerule
    &$ D\sb{4} + A\sb{3} +  14 A\sb{1} $ &&$4 $ &&$2, 3, 4, 5, 6, $ &\cr    &$  $ && &&$7, 8$ &\cr    \tablerule
    &$ D\sb{4} + A\sb{2} +  15 A\sb{1} $ &&$6 $ &&$2, 3, 4, 5, 6, $ &\cr    &$  $ && &&$7, 8$ &\cr    \tablerule
    &$ D\sb{4} +  17 A\sb{1} $ &&$2 $ &&$2, 3, 4, 5, 6, $ &\cr    &$  $ && &&$7, 8, 9$ &\cr    \tablerule
    &$ A\sb{19} +  2 A\sb{1} $ &&$20 $ &&$1$ &\cr    \tablerule
  }
 }
}
%
%
\hbox{
  \vtop{\tabskip=0pt \offinterlineskip
    \halign to \colwidth {\strut\vrule#\tabskip =\ttskip plus\pttskip&#\hfil&\vrule#&\hfil#& \vrule#& #\hfil & \vrule#\tabskip=0pt\cr
    \tablerule
    &\multispan5 \; $p=2$\hfil &\cr
    \tablerule
    \tablerule
    &\hfil$R$&& \hfil$n$ && \hfil$\sigma$ &\cr
    \tablerule
    \tablerule
    &$ A\sb{18} +  3 A\sb{1} $ &&$38 $ &&$1$ &\cr    \tablerule
    &$ A\sb{17} + A\sb{3} + A\sb{1} $ &&$36 $ &&$1$ &\cr    \tablerule
    &$ A\sb{17} + A\sb{3} + A\sb{1} $ &&$4 $ &&$1$ &\cr    \tablerule
    &$ A\sb{17} + A\sb{2} +  2 A\sb{1} $ &&$6 $ &&$1$ &\cr    \tablerule
    &$ A\sb{17} +  4 A\sb{1} $ &&$2 $ &&$1, 2$ &\cr    \tablerule
    &$ A\sb{16} +  5 A\sb{1} $ &&$34 $ &&$2$ &\cr    \tablerule
    &$ A\sb{15} + A\sb{4} +  2 A\sb{1} $ &&$20 $ &&$1$ &\cr    \tablerule
    &$ A\sb{15} + A\sb{3} + A\sb{2} + A\sb{1} $ &&$6 $ &&$1$ &\cr    \tablerule
    &$ A\sb{15} + A\sb{3} +  3 A\sb{1} $ &&$2 $ &&$1, 2$ &\cr    \tablerule
    &$ A\sb{15} +  6 A\sb{1} $ &&$4 $ &&$2, 3$ &\cr    \tablerule
    &$ A\sb{14} + A\sb{3} +  2 A\sb{2} $ &&$60 $ &&$1$ &\cr    \tablerule
    &$ A\sb{14} + A\sb{2} +  5 A\sb{1} $ &&$10 $ &&$2$ &\cr    \tablerule
    &$ A\sb{13} + A\sb{5} + A\sb{3} $ &&$84 $ &&$1$ &\cr    \tablerule
    &$ A\sb{13} + A\sb{4} +  4 A\sb{1} $ &&$70 $ &&$1, 2$ &\cr    \tablerule
    &$ A\sb{13} + A\sb{2} +  6 A\sb{1} $ &&$42 $ &&$2, 3$ &\cr    \tablerule
    &$ A\sb{12} + A\sb{6} +  3 A\sb{1} $ &&$182 $ &&$1$ &\cr    \tablerule
    &$ A\sb{12} + A\sb{4} +  5 A\sb{1} $ &&$130 $ &&$2$ &\cr    \tablerule
    &$ A\sb{12} + A\sb{3} +  6 A\sb{1} $ &&$52 $ &&$2, 3$ &\cr    \tablerule
    &$ A\sb{11} + A\sb{9} + A\sb{1} $ &&$60 $ &&$1$ &\cr    \tablerule
    &$ A\sb{11} + A\sb{6} + A\sb{3} + A\sb{1} $ &&$42 $ &&$1$ &\cr    \tablerule
    &$ A\sb{11} + A\sb{6} +  2 A\sb{2} $ &&$84 $ &&$1$ &\cr    \tablerule
    &$ A\sb{11} +  2 A\sb{5} $ &&$12 $ &&$1$ &\cr    \tablerule
    &$ A\sb{11} + A\sb{5} + A\sb{3} + A\sb{2} $ &&$6 $ &&$1$ &\cr    \tablerule
    &$ A\sb{11} + A\sb{5} + A\sb{3} +  2 A\sb{1} $ &&$2 $ &&$1, 2$ &\cr    \tablerule
    &$ A\sb{11} + A\sb{5} +  5 A\sb{1} $ &&$4 $ &&$1, 2, 3$ &\cr    \tablerule
    &$ A\sb{11} +  2 A\sb{3} +  2 A\sb{2} $ &&$12 $ &&$1, 2$ &\cr    \tablerule
    &$ A\sb{11} + A\sb{3} +  7 A\sb{1} $ &&$6 $ &&$2, 3, 4$ &\cr    \tablerule
    &$ A\sb{11} +  2 A\sb{2} +  6 A\sb{1} $ &&$12 $ &&$2, 3$ &\cr    \tablerule
    &$ A\sb{10} + A\sb{9} +  2 A\sb{1} $ &&$110 $ &&$1$ &\cr    \tablerule
    &$ A\sb{10} + A\sb{7} +  4 A\sb{1} $ &&$88 $ &&$1, 2$ &\cr    \tablerule
    &$ A\sb{10} + A\sb{6} + A\sb{3} + A\sb{2} $ &&$924 $ &&$1$ &\cr    \tablerule
    &$ A\sb{10} + A\sb{5} +  6 A\sb{1} $ &&$66 $ &&$2, 3$ &\cr    \tablerule
    &$ A\sb{10} + A\sb{4} + A\sb{2} +  5 A\sb{1} $ &&$330 $ &&$2$ &\cr    \tablerule
    &$ A\sb{10} +  11 A\sb{1} $ &&$22 $ &&$4, 5$ &\cr    \tablerule
    &$  2 A\sb{9} + A\sb{3} $ &&$4 $ &&$1$ &\cr    \tablerule
  }
 }
  \vtop{\tabskip=0pt \offinterlineskip
    \halign to \colwidth {\strut\vrule#\tabskip =\ttskip plus\pttskip&#\hfil&\vrule#&\hfil#& \vrule#& #\hfil & \vrule#\tabskip=0pt\cr
    \tablerule
    &\multispan5 \; $p=2$\hfil &\cr
    \tablerule
    \tablerule
    &\hfil$R$&& \hfil$n$ && \hfil$\sigma$ &\cr
    \tablerule
    \tablerule
    &$  2 A\sb{9} + A\sb{2} + A\sb{1} $ &&$6 $ &&$1$ &\cr    \tablerule
    &$  2 A\sb{9} +  3 A\sb{1} $ &&$2 $ &&$1, 2$ &\cr    \tablerule
    &$ A\sb{9} + A\sb{6} + A\sb{3} +  3 A\sb{1} $ &&$140 $ &&$1, 2$ &\cr    \tablerule
    &$ A\sb{9} +  2 A\sb{4} +  4 A\sb{1} $ &&$10 $ &&$1, 2$ &\cr    \tablerule
    &$ A\sb{9} + A\sb{3} + A\sb{2} +  7 A\sb{1} $ &&$60 $ &&$2, 3, 4$ &\cr    \tablerule
    &$ A\sb{9} +  12 A\sb{1} $ &&$10 $ &&$3, 4, 5, 6$ &\cr    \tablerule
    &$ A\sb{8} + A\sb{5} + A\sb{2} +  6 A\sb{1} $ &&$18 $ &&$2, 3$ &\cr    \tablerule
    &$ A\sb{8} + A\sb{4} + A\sb{3} +  6 A\sb{1} $ &&$180 $ &&$2, 3$ &\cr    \tablerule
    &$ A\sb{8} +  13 A\sb{1} $ &&$18 $ &&$4, 5, 6$ &\cr    \tablerule
    &$  2 A\sb{7} +  2 A\sb{3} + A\sb{1} $ &&$2 $ &&$1, 2$ &\cr    \tablerule
    &$  2 A\sb{7} + A\sb{3} +  4 A\sb{1} $ &&$4 $ &&$1, 2, 3$ &\cr    \tablerule
    &$ A\sb{7} + A\sb{6} + A\sb{4} +  4 A\sb{1} $ &&$280 $ &&$1, 2$ &\cr    \tablerule
    &$ A\sb{7} + A\sb{6} +  2 A\sb{3} + A\sb{2} $ &&$168 $ &&$1, 2$ &\cr    \tablerule
    &$ A\sb{7} + A\sb{5} + A\sb{4} +  5 A\sb{1} $ &&$120 $ &&$1, 2, 3$ &\cr    \tablerule
    &$ A\sb{7} +  4 A\sb{3} + A\sb{2} $ &&$24 $ &&$1, 2, 3$ &\cr    \tablerule
    &$ A\sb{7} +  2 A\sb{3} +  8 A\sb{1} $ &&$8 $ &&$3, 4, 5$ &\cr    \tablerule
    &$ A\sb{7} + A\sb{2} +  12 A\sb{1} $ &&$24 $ &&$3, 4, 5, 6$ &\cr    \tablerule
    &$  3 A\sb{6} + A\sb{3} $ &&$28 $ &&$1$ &\cr    \tablerule
    &$  3 A\sb{6} +  3 A\sb{1} $ &&$14 $ &&$1$ &\cr    \tablerule
    &$ A\sb{6} +  3 A\sb{5} $ &&$42 $ &&$1$ &\cr    \tablerule
    &$ A\sb{6} + A\sb{5} + A\sb{3} +  7 A\sb{1} $ &&$84 $ &&$3, 4$ &\cr    \tablerule
    &$ A\sb{6} + A\sb{4} +  11 A\sb{1} $ &&$70 $ &&$4, 5$ &\cr    \tablerule
    &$ A\sb{6} + A\sb{2} +  13 A\sb{1} $ &&$42 $ &&$4, 5, 6$ &\cr    \tablerule
    &$  4 A\sb{5} + A\sb{1} $ &&$2 $ &&$1, 2$ &\cr    \tablerule
    &$  3 A\sb{5} +  6 A\sb{1} $ &&$6 $ &&$1, 2, 3, 4$ &\cr    \tablerule
    &$ A\sb{5} + A\sb{4} +  12 A\sb{1} $ &&$30 $ &&$3, 4, 5, 6$ &\cr    \tablerule
    &$ A\sb{5} + A\sb{3} +  13 A\sb{1} $ &&$12 $ &&$3, 4, 5, 6, 7$ &\cr    \tablerule
    &$ A\sb{4} + A\sb{3} +  14 A\sb{1} $ &&$20 $ &&$3, 4, 5, 6, 7$ &\cr    \tablerule
    &$  7 A\sb{3} $ &&$4 $ &&$1, 2, 3, 4$ &\cr    \tablerule
    &$ A\sb{3} +  18 A\sb{1} $ &&$4 $ &&$2, 3, 4, 5, 6, $ &\cr    &$  $ && &&$7, 8, 9$ &\cr    \tablerule
    &$ A\sb{2} +  19 A\sb{1} $ &&$6 $ &&$2, 3, 4, 5, 6, $ &\cr    &$  $ && &&$7, 8, 9$ &\cr    \tablerule
    &$  21 A\sb{1} $ &&$2 $ &&$1, 2, 3, 4, 5, $ &\cr    &$  $ && &&$6, 7, 8, 9, 10$ &\cr    \tablerule
  }
 }
}
%
%
%
%
\bibliographystyle{amsplain}

\providecommand{\bysame}{\leavevmode\hbox to3em{\hrulefill}\thinspace}

\end{document}